\documentclass[11pt,reqno]{amsart}
\usepackage{amssymb,mathrsfs,graphicx,subfigure, enumerate}
\usepackage{amsmath,amsfonts,amssymb,amscd,amsthm,bbm}
\usepackage{extpfeil}
\usepackage{graphicx,colortbl}
\usepackage{float}
\usepackage{epsfig}
\usepackage{caption}
\usepackage{bm}
\graphicspath{{Figure/}}

\usepackage[margin=1in]{geometry}

\title[Stability of composite wave of shocks for 3D Navier-Stokes]{Long-time behavior toward composite wave of shocks for 3D barotropic Navier-Stokes system}

\author[Kang]{Moon-Jin Kang}
\address[Moon-Jin Kang]{\newline Department of Mathematical Sciences \newline Korea Advanced Institute of Science and Technology, Daejeon  34141, Republic of Korea}
\email{moonjinkang@kaist.ac.kr}

\author[Lee]{Hobin Lee}
\address[Hobin Lee]{\newline Department of Mathematical Sciences \newline Korea Advanced Institute of Science and Technology, Daejeon  34141, Republic of Korea}
\email{lcuh11@kaist.ac.kr}

\newtheorem{theorem}{Theorem}[section]
\newtheorem{lemma}{Lemma}[section]

\newtheorem{proposition}{Proposition}[section]
\newtheorem{remark}{Remark}[section]

\newcommand{\beq}{\begin{equation}}
\newcommand{\eeq}{\end{equation}}

\newcommand{\bbr}{\mathbb R}

\newcommand{\bbt}{\mathbb T}
\newcommand{\eps}{\varepsilon }

\newcommand{\R}{\mathbb{R}}

\newcommand{\bq}{\begin{equation}}
\newcommand{\eq}{\end{equation}}
\newcommand{\e}{\varepsilon}

\newcommand{\pa}{\partial}

\newcommand{\cB}{\mathcal{B}}
\newcommand{\cG}{\mathcal{G}}

\newcommand{\cD}{\mathcal{D}}

\newcommand{\tU}{\widetilde{U}}
\newcommand{\tu}{\widetilde{u}}

\newcommand{\tv}{\widetilde{v}}
\newcommand{\tvi}{\widetilde{v}_i}

\newcommand{\tr}{\widetilde{\rho}}

\newcommand{\pv}{p(v)}

\newcommand{\tpv}{p(\widetilde{v})}

\newcommand{\norm}[1]{\left\lVert#1\right\rVert}
\newcommand{\di}{\displaystyle}

\newcommand{\bu}{\bold{u}}

\begin{document}
\bibliographystyle{plain}

\date{\today}


\thanks{\textbf{Acknowledgment.} 
The authors are partially supported by  the National Research Foundation of Korea  (NRF-2019R1C1C1009355, NRF-2019R1A5A1028324 and NRF-RS-2024-00361663)}

\begin{abstract} 
We consider the barotropic Navier-Stokes system in three space dimensions with periodic boundary condition in the transversal direction. 
	We show the long-time behavior of the 3D barotropic Navier-Stokes flow perturbed from a composition of two shock waves with suitably small amplitudes. We prove that the perturbed Navier-Stokes flow converges, uniformly in space, towards a composition of two planar viscous shock waves  as time goes to infinity, up to dynamical shifts. This is the first result on time-asymptotic stability of composite wave of two shocks for multi-D Navier-Stokes system.
	The main part of proof is based on the method of $a$-contraction with shifts.  
	\end{abstract}

\maketitle

\tableofcontents

\section{Introduction}\label{sec:1}
\setcounter{equation}{0}

We consider the 3D barotropic compressible Navier-Stokes equations :
\begin{equation}  \label{main}
\begin{cases} \partial_t\rho+div_x(\rho \bold{u})=0,&  \\
\partial_t(\rho u)+div_x(\rho \bu \otimes \bu)+\nabla_xp(\rho)=\mu\Delta_x\bu+(\mu+\lambda)\nabla_x div_x\bu.&
\end{cases}
\end{equation}
Here, $\rho=\rho(t,x)$, $\bold{u}=\bold{u}(t,x)=(u_1,u_2,u_3)^T(t,x)$ denote the mass density and the velocity field of a fluid, respectively. The pressure $p(\rho)=b\rho^{\gamma}$ $(b>0, \gamma>1)$ follows the $\gamma$-law, and the two constants $\mu$ and $\lambda$ represent the viscosity coefficients satisfying the physical constraints 
\[
\mu>0, \quad 2\mu+3\lambda\ge 0.
\]
We handle the system \eqref{main} for $ x=(x_1,x_2,x_3)\in \Omega :=\mathbb{R}\times \mathbb{T}^2$ with $\mathbb{T}^2:=\mathbb{R}^2/\mathbb{Z}^2$, where the periodic boundary condition is considered for the transversal direction. \\
Consider an initial datum of the system \eqref{main} that connects prescribed far-field constant states :
\begin{equation} \label{initial data}
(\rho, \bold{u})|_{t=0}=(\rho_0, \bold{u}_0) \to (\rho_{\pm}, \bold{u}_{\pm}), \quad as \quad x_1 \to \pm\infty , 
\end{equation}
where $\rho_{\pm}>0$ and $\bold{u}_{\pm}=(u_{\pm},0,0)$. 

It is known from a heuristic argument (c.f. \cite{MatsumuraBook}) that the asymptotic profile of large-time behavior of Navier-Stokes solutions is related to the Riemann solution to the associated Euler system with the Riemann data composed of the above states:
\begin{equation}   \label{eq:Euler}
\begin{cases}
\partial_t\rho+div_x(\rho \bold{u})=0, \\
\partial_t(\rho \bold{u})+div_x(\rho \bold{u} \otimes \bold{u})+\nabla_x p(\rho)=0, \\
(\rho, \bold{u})(0,x)=
\begin{cases}
(\rho_-, \bold{u}_-), \quad x_1<0 \\
(\rho_+, \bold{u}_+), \quad x_1>0.
\end{cases}
\end{cases}
\end{equation}

In this paper, we are interested in the long-time behavior of solutions to the Navier-Stokes system \eqref{main}-\eqref{initial data} when the Riemann data with the end states $(v_\pm,u_\pm)$ that generate the Riemann solution composed of two shock waves.
The Riemann solution is an one-dimensional (self-similar entropy) solution to the Riemann problem of \eqref{eq:Euler} as a solution to
\begin{equation}\label{viscous shock}
\begin{cases}
\partial_t\rho+\partial_{x_1}(\rho u_1)=0, \\
\partial_t(\rho u_1)+\partial_{x_1}(\rho u_1^2+p(\rho))=0,\\
(\rho, u_1)(0,x_1)=
\begin{cases}
(\rho_-, u_-), \quad x_1<0 \\
(\rho_+, u_+), \quad x_1>0.
\end{cases}
\end{cases}
\end{equation}

To consider the Riemann solution composed of two shock waves, we assume that there exists a unique intermediate state $(\rho_m,u_m)$ which is connected with $(\rho_-,u_-)$ by 1-shock curve and with $(\rho_+,u_+)$ by 2-shock curve. That is, there exists a unique $(\rho_m,u_m)$ such that the following Rankine-Hugoniot condition and Lax entropy condition hold:
\begin{equation}\label{RHCLEC}
\begin{aligned}
&\begin{cases}
-\sigma_1(\rho_m-\rho_-)+(\rho_{m}u_{m}-\rho_-u_{-})=0, \\
-\sigma_1(\rho_mu_{m}-\rho_-u_{-})+(\rho_mu_{m}^2-\rho_-u_{-}^2)+(p(\rho_m)-p(\rho_-))=0, \quad \rho_-< \rho_m,\quad u_->u_m;\\
\end{cases} \\
&\begin{cases}
-\sigma_2(\rho_+-\rho_m)+(\rho_{+}u_{+}-\rho_mu_{m})=0, \\
-\sigma_2(\rho_+u_{+}-\rho_mu_{m})+(\rho_+u_{+}^2-\rho_mu_{m}^2)+(p(\rho_+)-p(\rho_m))=0,\quad \rho_m>\rho_+,\quad u_m>u_+.
\end{cases} 
\end{aligned}
\end{equation}
For such a Riemann data, the associated Riemann solution is the superposition $(\bar{\rho},\bar{u})=(\rho^s_1,u^s_1)+(\rho^s_2,u^s_2)-(\rho_m,u_m)$ of 1-shock wave $(\rho^s_1,u^s_1)$ and 2-shock wave $(\rho^s_2,u^s_2)$ defined as
\[(\rho^s_1,u^s_1)(t,x_1) = \begin{cases}
(\rho_-,u_-),\quad x_1<\sigma_1t,\\
(\rho_m,u_m),\quad x_1>\sigma_1t,
\end{cases},\quad (\rho^s_2,u^s_2)(t,x_1) = \begin{cases}
(\rho_m,u_m),\quad x_1<\sigma_2t,\\
(\rho_+,u_+),\quad x_1>\sigma_2t.
\end{cases}\]
The viscous counterpart of the Riemann solution $(\bar{v},\bar{u})$ is given by the composite wave:
\beq\label{com2}
(\tr(t,x),\tu(t,x)):=\Big(\tr_1(x_1-\sigma_1 t),\tu_1(x_1-\sigma_1t)\Big)+\Big(\tr_2(x_1-\sigma_2 t),\tu_2(x_1-\sigma_2 t)\Big)-(v_m,u_m),
\eeq
which is composed of 1-viscous shock $(\tr_1,\tu_1)(x_1-\sigma_1t)$ and 2-viscous shock $(\tr_2,\tu_2)(x_1-\sigma_2t)$ satisfying: for each $i=1,2,$
\begin{equation}\label{viscous-shock-u}
\begin{cases}
-\sigma_i(\widetilde{\rho}_i)'+(\widetilde{\rho}_i\widetilde{u}_{i})'=0, \\
-\sigma_i(\widetilde{\rho}_i\widetilde{u}_{i})'+(\widetilde{\rho}_i(\widetilde{u}_{i})^2)'+p(\widetilde{\rho}_i)'=(2\mu+\lambda)(\widetilde{u}_{i})'' ,\\
(\tr_1,\tu_1)(-\infty) = (\rho_-,u_-),\quad (\tr_1,\tu_1)(+\infty) = (\rho_m,u_m),\\
(\tr_2,\tu_2)(-\infty) = (\rho_m,u_m),\quad (\tr_2,\tu_2)(+\infty) = (\rho_+,u_+).
\end{cases}
\end{equation} 
Notice that each viscous shock $(\widetilde{v}_i(x_1-\sigma_it), \widetilde{\bold{u}}_{i}(x_1-\sigma_it))$ is a traveling wave solution to \eqref{main}, where
we use the notations $\widetilde{\bold{u}}_i(x_1-\sigma_it):=(\widetilde{u}_i(x_1-\sigma_it), 0, 0)^T$ and $\widetilde{\bold{u}}(t,x):=(\widetilde{u}(t,x),0,0)^T$.\\

For previous results on time-asymptotic stability for the barotropic Navier-Stokes system when the far-field states are connected by a single shock, we refer to the first result by Matsumura-Nishihara  \cite{MN-S} in 1D, and its improvements 
\cite{MN-S,Liu,LZ,SX,HLZ10,M-Zumbrun}.  
Recently, Wang-Wang \cite{WW} show the time-asymptotic stability towards a single planar viscous (weak) shock for \eqref{main}-\eqref{initial data} under 3D perturbations with periodic boundary condition in transversal direction. We also refer to \cite{HLZ17} for stability result of a planar viscous shock with large amplitude under multi-dimensional perturbations on whole space $\bbr^n$.  
We aims to extend the result \cite{WW} to more general case where the far-field states are connected by two shocks. More precisely, we prove that solutions to the 3D Navier-Stokes system \eqref{main}-\eqref{initial data} with \eqref{RHCLEC} converge to the composite wave $(\tr,\tu)$ up to shifts, uniformly in $x$ as $t\to\infty$.

For the 1D results on long-time behavior towards two shocks, we refer to Huang-Matsumura \cite{HM} and Han-Kang-Kim \cite{HKK23}.
We also refer to \cite{KVW-2022,KVW-NSF} for generic superposition case of shock and rarefaction wave possibly with contact discontinuity. For the results on stability of planar rarefaction waves for 3D Navier-Stokes, we refer to \cite{LWW2,LWW-1}. \\

We now state the main result as follows.
\begin{theorem}\label{thm:main}
	For a given constant state $(\rho_+,u_+)\in\R_+\times\R$, there exist positive constants $\delta^0,\e^0$ such that the following holds.\\	
	For any constant states $(\rho_m,u_m)$ and $(\rho_-,u_-)$ satisfying \eqref{RHCLEC}  with
	\beq\label{assind}
	|\rho_+-\rho_m|+|\rho_m-\rho_-|<\delta^0,
	\eeq
	let $(\tr_i,\tu_i)(x_1-\sigma_it)$ be the $i$-viscous shock wave satisfying \eqref{viscous-shock-u} for each $i=1,2$. In addition, let $(\rho_0,\bu_0)$ be any initial data satisfying
	\[\sum_{\pm} \left(\|\rho_0-\rho_\pm\|_{L^2(\mathbb{R}_{\pm}\times\mathbb{T}^2)}+\|\bu_0- \bu_\pm\|_{L^2(\mathbb{R}_{\pm}\times\mathbb{T}^2)}\right)+ \|\nabla_x\rho_0\|_{L^2(\Omega)}+\|\nabla_x \bu_0\|_{L^2(\Omega)}<\e^0,\]
	where $\R_+:=(0,+\infty)$ and $\R_-:=(-\infty,0)$. Then, the compressible Navier-Stokes system  \eqref{main}-\eqref{initial data} with \eqref{RHCLEC}  admits a unique global-in-time solution $(\rho, \bu)$ in the following sense: there exist Lipschitz continuous shift functions $X_1(t), X_2(t)$ such that
	\begin{align*}
	&\rho(t,x)-(\tr_1(x_1-\sigma_1t-X_1(t))+\tr_2(x_1-\sigma_2t-X_2(t))-\rho_m)\in C(0,+\infty;H^2(\Omega)),\\
	&\bu(t,x)-\Big(\widetilde{\bu}_1(x_1-\sigma_1t-X_1(t))+\widetilde{\bu}_2(x_1-\sigma_2t-X_2(t))- \bu_m \Big) \in C(0,+\infty;H^2(\Omega)),
	\end{align*}
	where $\bu_m:=(u_m,0,0)^T$.
	Moreover, we have the long-time behavior:
	\begin{align*}
	&\lim_{t\to+\infty}\sup_{x\in \Omega}\Big|\rho(t,x)-\big(\tr_1(x_1-\sigma_1t-X_1(t))+\tr_2(x_1-\sigma_2t-X_2(t))-\rho_m\big)\Big|=0,\\
	&\lim_{t\to+\infty}\sup_{x\in \Omega}\Big| \bu(t,x)-\Big(\widetilde{\bu}_1(x_1-\sigma_1t-X_1(t))+\widetilde{\bu}_2(x_1-\sigma_2t-X_2(t))- \bu_m \Big) \Big|=0,
	\end{align*}
	where 
	\beq\label{as1}
	\lim_{t\to+\infty}|\dot{X}_i(t)|=0,\quad\mbox{for}\quad i=1,2.
	\eeq
	Especially,
\beq\label{sx121}
	X_1(t)+\sigma_1t \le \frac{3\sigma_1+\sigma_2}{4}t < \frac{\sigma_1+\sigma_2}{2}t < \frac{\sigma_1+3\sigma_2}{4}t  \le X_2(t) +\sigma_2t ,\quad t>0.
	\eeq
\end{theorem}

\begin{remark}
	1.
	The property \eqref{as1} implies that 
	\[
	\lim_{t\to+\infty}\frac{X_i(t)}{t}=0,\quad\mbox{for}\quad i=1,2.
	\]
	Thus, the shift functions $|X_i(t)|$ grow at most sub-linearly, and so the shifted composite wave tends to the original composite wave $(\tr,\tu)$ time-asymptotically.\\
	2. The two shifts are well-separated in the sense of \eqref{sx121}, which holds from our construction on the shifts  in the proof.  Indeed,  it will be ensured that $|\dot X_i(t)|\ll 1$ with $X_i(0)=0$ and so,
	\[
	|X_i(t)| \le  \frac{\sigma_2-\sigma_1}{4} t, \quad t>0.
	\]
	
	\end{remark}
	
	\begin{remark}
	 We will use the method of $a$-contraction with shifts to prove the main theorem.
	The method of $a$-contraction with shifts was introduced in \cite{KVARMA,Vasseur-2013} for the study on stability of extreme shocks in the hyperbolic system of conservation laws (especially for the Euler system). The first extension of the method to viscous conservation laws was done in the 1D scalar case \cite{Kang-V-1}  (\cite{Kang19} for more general fluxes), and then in the multi-D case \cite{KO,KVW}. This method was extended to the 1D Navier-Stokes system to show the contraction property of any large perturbations for a single shock in \cite{KV21,KV-Inven}, and for a composite wave of two shocks in \cite{KV-2shock}. Furthermore, it was also used to show the long-time behavior of the barotropic NS system towards the composition of shock and rarefaction under the 1D perturbation in \cite{KVW-2022}, and towards a single shock under the multi-D perturbation in \cite{WW}. Recently, it was extended to the study on long-time behavior towards generic Riemann solution for 1D Navier-Stokes-Fourier system in \cite{KVW-NSF}. For extension to more general Cauchy problem of Euler and Navier-Stokes systems (and related models), we refer to \cite{CKV,CKV20,CKV24}.  
We also refer to \cite{SV-16,SV-16dcds,Serre-2015} for the extension to abstract study on hyperbolic system of conservation laws.	
\end{remark}

The paper is organized as follows. In Section \ref{sec:2}, we present a reformulation for the NS system with the specific volume, together with a new statement for the main result. 
In Section 3, we present the main proposition for the a priori estimates that completes the proof of the main result. 
 In Section 4, we use the method of $a$-contraction with shifts to obtain the zeroth order estimates. Then, the higher order estimates will be obtained  in Section \ref{sec:5}.

\section{Reformulation of the main result}\label{sec:2}
\setcounter{equation}{0}
A starting point of the main strategy for the proof is to rewrite the system by using the specific volume $v=\frac{1}{\rho}$ in the Eulerian coordinates. This allows us to use the method of $a$-contraction with shifts easily as done in \cite{KV21,KV-Inven,KV-2shock,WW}, but we prefer to use the Eulerian frame rather than the mass Lagrangian frame due to the multi dimensionality.  \\

As in \cite{WW}, we use the decomposition of the Laplacian $\Delta \bold{u}$ into the irrotational part $\nabla (\nabla\cdot\bold{u})$ and the rotational part $\nabla\times \nabla\times \bold{u}$ as
\[\Delta_{x}\bold{u}=\nabla_xdiv_x\bold{u}- \nabla_{x}\times \nabla_{x}\times \bold{u},\]
which is importantly used to control the transversal perturbations around the planar shock. \\
Together with the decomposition, we use the specific volume $v=\frac{1}{\rho}$ to rewrite the system \eqref{main} into 
(as in \cite{M01} for 1D case)
\begin{equation} \label{changed main}
\begin{cases}
\rho\left( \partial_tv+\bold{u}\cdot \nabla_{x}v \right)=div_{x}\bold{u}, \\
\rho\left( \partial_t\bold{u}+\bold{u}\nabla_{x}\bold{u} \right)+\nabla_{x}p(v)=(2\mu+\lambda)\nabla_{x}div_{x}\bold{u}-\mu \nabla_{x}\times \nabla_{x}\times \bold{u}, \\
\end{cases}
\end{equation}
with the initial datum
\beq \label{nini}
(v, \bold{u})|_{t=0}=(v_0, \bold{u}_0) \to (v_{\pm}, \bold{u}_{\pm}), \quad as \quad x_1 \to \pm\infty , 
\eeq
where $p(v)=v^{-\gamma}$ with normalized coefficient $b=1$ for simplicity. \\
Similarly, by using $\widetilde{v}_i:=\frac{1}{\widetilde{\rho}_i}$ for each $i=1,2$, the ODE system \eqref{viscous-shock-u} can be written as
\begin{equation} \label{changed viscous shock equation}
\begin{cases}
\widetilde{\rho}_i\left(-\sigma_i(\widetilde{v}_i)+\widetilde{u}_{i}(\widetilde{v}_i)'\right)=(\widetilde{u}_{i})', \\
\widetilde{\rho}_i\left(-\sigma_i(\widetilde{u}_{i})'+\widetilde{u}_{i}(\widetilde{u}_{i})'\right)+p(\widetilde{v}_i)'=(2\mu+\lambda)(\widetilde{u}_{i})'', 
\end{cases}    
\end{equation}
where $p(\widetilde{v}_i)=(\widetilde{v}_i)^{-\gamma}$. Integrating \eqref{viscous-shock-u}$_1$ on $(-\infty, x_1-\sigma_i t)$, we have
\begin{equation} \label{changed RH condition}
\begin{cases}
-\sigma_1\widetilde{\rho}_1+\widetilde{\rho}_1\widetilde{u}_{1}=-\sigma_1\rho_-+\rho_-u_{-}, \\ 
-\sigma_2\widetilde{\rho}_2+\widetilde{\rho}_2\widetilde{u}_{2}=-\sigma_2\rho_m+\rho_mu_{m} .
\end{cases}    
\end{equation}
Let 
\[
\sigma_1^* := \sigma_1\rho_- - \rho_-u_{-},\qquad \sigma_2^* := \sigma_2\rho_m - \rho_mu_{m}.
\]
Then, by \eqref{changed RH condition}, the system \eqref{changed viscous shock equation} and the far-field conditions can be rewritten as for each $i=1,2,$
\begin{equation} \label{change shock equation 1}
\begin{cases}
-\sigma_{i}^*(\widetilde{v}_i)'=(\widetilde{u}_{i})',\\ 
-\sigma_{i}^*(\widetilde{u}_{i})'+p(\widetilde{v}_i)'=(2\mu+\lambda)(\widetilde{u}_{i})'',
\end{cases}    
\end{equation}
and
\begin{equation} \label{change shock equation 2}
\begin{cases}
(\widetilde{v}_1,\widetilde{u}_{1})(-\infty)=(v_-,u_{-}), \quad (\widetilde{v}_1,\widetilde{u}_{1})(+\infty)=(v_m,u_{m}), \quad v_-=\frac{1}{\rho_-}, \  v_m=\frac{1}{\rho_m}, \\ 
(\widetilde{v}_2,\widetilde{u}_{2})(-\infty)=(v_m,u_{m}), \quad (\widetilde{v}_2,\widetilde{u}_{2})(+\infty)=(v_+,u_{+}), \quad v_+=\frac{1}{\rho_+}.
\end{cases}
\end{equation}
Thus, we have
\begin{equation} \label{changed RH condition 2}    
\begin{cases} 
-\sigma_1^*(v_m-v_-)=u_{m}-u_{-}. \\
-\sigma_1^*(u_{m}-u_{-})+p(v_m)-p(v_-)=0,
\end{cases}
\begin{cases}
-\sigma_2^*(v_+-v_m)=u_{+}-u_{m}. \\
-\sigma_2^*(u_{+}-u_{m})+p(v_+)-p(v_m)=0,
\end{cases}
\end{equation}
from which,
 \beq\label{sspeed}
 \sigma_1^*=-\sqrt{-\frac{p(v_m)-p(v_-)}{v_m-v_-}},\qquad \sigma_2^*=\sqrt{-\frac{p(v_+)-p(v_m)}{v_+-v_m}}.
\eeq
\subsection{Estimates on viscous shocks}
Since the ODE system \eqref{change shock equation 1} with \eqref{change shock equation 2}-\eqref{sspeed} is the same as the one for viscous shocks of Navier-Stokes system in the mass Lagrangian coordinates, the following lemma on estimates of viscous shocks is verified by the previous results \cite{HKK,Kang-V-2shocks,Kang-V-NS17,MN-S}. 
\begin{lemma}\label{lem:shock-est}
	For a given constant $U_{\ast}:=(v_{\ast},u_{\ast})\in \mathbb{R}_+\times \mathbb{R}$, there exist positive constants $\delta_0$, $C$, $C_1$, and $C_2$ such that the following holds. 
 Let $U_-:=(v_-, u_{-})$, $U_m:=(v_m,u_{m})$, and $U_+:=(v_+,u_{+})\in \mathbb{R}_+\times \mathbb{R}$ be 
any constants such that $U_-$, $U_m$, $U_+\in B_{\delta_0}(U_{\ast})$, and $|p(v_-)-p(v_m)|=:\delta_1<\delta_0$ and 
$|p(v_m)-p(v_+)|=:\delta_2<\delta_0$. Let $\widetilde{U}_1=(\widetilde{v}_1,\widetilde{u}_{1})$ and $\widetilde{U}_2=(\widetilde{v}_2, \widetilde{u}_{2})$ be the $1-$ and 
$2-$shocks connecting from $U_-$ to $U_m$ and from $U_m$ to $U_+$ respectively, 
satisfying $\widetilde{v}_1(0)=\frac{v_-+v_m}{2}$ and $\widetilde{v}_2(0)=\frac{v_m+v_+}{2}$ without loss of generality. Let $\xi_i:=x_1-\sigma_it$ for $i=1,2$.

Then the following estimates holds : for each $i=1,2$,
\[
\widetilde{v}_i'\sim \widetilde{u}_{i}' \quad i.e., \quad C^{-1}\widetilde{v}_i'(\xi_i)\le \widetilde{u}_{i}'(\xi_i)\le C\widetilde{v}_i'(\xi_i), \quad \xi_i\in \mathbb{R},
\]
\[
C^{-1}\delta_i^2e^{-C_1\delta_i|\xi_i|}\le \widetilde{v}_i'(\xi_i)\le -C\delta_i^2e^{-C_2\delta_i|\xi_i|}, \quad \xi_i\in\mathbb{R},
\] 
in addition, 
\[
|(\widetilde{v}_i''(\xi_i),\widetilde{u}_{i}''(\xi_i))|\le C\delta_i|(\widetilde{v}_i'(\xi_i),\widetilde{u}_{i}'(\xi_i))|,
\]
\[
|(\widetilde{v}_i'''(\xi_i),\widetilde{u}_{i}'''(\xi_i))|\le C\delta_i^2|(\widetilde{v}_i'(\xi_i), \widetilde{u}_{i}'(\xi_i))|.
\]
\end{lemma}

\begin{remark}
	Throughout the paper, we will use Lemma \ref{lem:shock-est} with $U_*=U_+$. Thus, the constants $\delta_0$, $C$, $C_1$ and $C_2$ depend only on $U_+$. 
\end{remark}

\subsection{Statement for the main result}

In the remaining part of the paper, we will prove the following theorem that is stated in terms of the volume variable. Then, Theorem \ref{thm} obviously implies Theorem \ref{thm:main}, since we are considering small perturbations in $L^\infty ((0,\infty)\times\Omega)$ and so the values of $v$ and $\rho$ stay near the reference point.

\begin{theorem} \label{thm}
For a given constant state $(v_+, u_+)\in \mathbb{R}_+\times \mathbb{R}$, there exists constants $\delta_0, \varepsilon_0>0$ such that the following holds true.

\noindent For any constant states $(v_m, u_m)$ and $(v_-,u_-)$ satisfying \eqref{changed RH condition 2} with
\begin{equation}
|v_+-v_m|+|v_m-v_-|<\delta_0,    
\end{equation}
let $(\widetilde{v}_i, \widetilde{u}_i)(x_1-\sigma_it)$ be the $i$-viscous shock wave satisfying \eqref{changed viscous shock equation}. In addition, let $(v_0, \bu_0)$ be any initial data satisfying
\[\sum\limits_{\pm}\left(\lVert v_0-v_{\pm} \rVert_{L^2(\mathbb{R}_{\pm}\times\mathbb{T}^2)}+\lVert \bu_0-\bu_{\pm} \rVert_{L^2(\mathbb{R}_{\pm}\times\mathbb{T}^2)}\right)+\lVert \nabla_xv_0 \rVert_{H^1(\Omega)}+\lVert \nabla_x\bu_0 \rVert_{H^1(\Omega)}<\varepsilon_0,\]
where $\mathbb{R}_+:=(0, +\infty)$ and $\mathbb{R}_-:=(-\infty, 0)$. Then the compressible Navier-Stokes system \eqref{changed main}-\eqref{nini} admits a unique global-in-time solution $(v, \bu)$ in the following sense: there exist Lipschitz continuous functions $X_1(t), X_2(t)$ such that

\begin{align}
\begin{aligned} \label{thmstate1}  
&v(t, x)-\left(\widetilde{v}_1(x_1-\sigma_1t-X_1(t))+\widetilde{v}_2(x_1-\sigma_2t-X_2(t))-v_m\right)\in C\left(0, +\infty; H^2(\Omega)\right),\\
&\bold{u}(t, x)-\left(\bold{\widetilde{u}}_1(x_1-\sigma_1t-X_1(t))+\bold{\widetilde{u}}_2(x_1-\sigma_2t-X_2(t))-\bold{u}_m\right)\in C\left(0, +\infty; H^2(\Omega)\right),\\
&\nabla_x\left(v(t, x)-\left(\widetilde{v}_1(x_1-\sigma_1t-X_1(t))+\widetilde{v}_2(x_1-\sigma_2t-X_2(t))-v_m\right)\right)\in L^2\left(0, +\infty; H^1(\Omega)\right), \\
&\nabla_x\left(\bold{u}(t, x)-\left(\widetilde{\bold{u}}_1(x_1-\sigma_1t-X_1(t))+\widetilde{\bold{u}}_2(x_1-\sigma_2t-X_2(t))-\bold{u}_m\right)\right)\in L^2\left(0, +\infty; H^2(\Omega)\right).
\end{aligned}
\end{align}

\noindent Moreover, we have the large-time behavior:
\begin{align}
\begin{aligned} \label{thmstate2}
&\lim\limits_{t\rightarrow +\infty}\sup\limits_{x\in \Omega}\left|v(t, x)-\left(\widetilde{v}_1(x_1-\sigma_1t-X_1(t))+\widetilde{v}_2(x_1-\sigma_2t-X_2(t))-v_m\right)\right|=0,\\
&\lim\limits_{t\rightarrow +\infty}\sup\limits_{x\in \Omega}\left|\bold{u}(t, x)-\left(\widetilde{\bold{u}}_1(x_1-\sigma_1t-X_1(t))+\widetilde{\bold{u}}_2(x_1-\sigma_2t-X_2(t))-\bold{u}_m\right)\right|=0,
\end{aligned}
\end{align}
where
\begin{equation} \label{as}
\lim\limits_{t\rightarrow +\infty}\frac{X_i(t)}{t}=0, \quad for \quad i=1, 2,  
\end{equation}
and
\beq
	X_1(t)+\sigma_1t \le \frac{3\sigma_1+\sigma_2}{4}t < \frac{\sigma_1+\sigma_2}{2}t < \frac{\sigma_1+3\sigma_2}{4}t  \le X_2(t) +\sigma_2t ,\quad t>0.
	\eeq
\end{theorem}

\section{Main proposition for proof of Theorem \ref{thm}}\label{sec:3}
\setcounter{equation}{0}
In this section, we present a main proposition on the a priori estimates for the proof of Theorem \ref{thm}. 

\subsection{Local existence}
We first recall the classical result on local-in-time existence of $H^2$ solutions to \eqref{main} connecting the two different constant states at far-fields.

\begin{proposition}
Let $\underline{v}$ and $\underline{u}$ be some smooth monotone functions in $\mathbb{R}$ such that 
\begin{equation*}
\underline{v}=v_{\pm}, \quad \underline{u}=u_{\pm}, \quad for \quad \pm x\ge 1.    
\end{equation*}
In addition, let $\underline{\bold{u}}=(\underline{u}, 0, 0)$.

\noindent Then for any constants $M_0$, $M_1$, $\underline{\kappa}_0$, $\overline{\kappa}_0$, $\underline{\kappa}_1$, $\overline{\kappa}_1$ with 
\begin{equation*}
0<M_0<M_1, \quad and \ \ 0<\underline{\kappa}_1<\underline{\kappa}_0<\overline{\kappa}_0<\overline{\kappa}_1,    
\end{equation*}
there exists a constant $T_0>0$ such that if the initial data $(v_0,\bold{u}_0)$ satisfy 
\begin{equation*}
\lVert v_0-\underline{v} \rVert_{H^2(\Omega)}+\lVert \bold{u}_0-\underline{\bold{u}}\rVert_{H^2(\Omega)}\le M_0, \quad and \ \ \underline{\kappa}_0\le v_0(x)\le \overline{\kappa}_0, \ \ \forall x\in \Omega,
\end{equation*}
then the Navier-Stokes equations \eqref{main} admit a unique soluition $(v, \bold{u})$ on $\left[0, T_0\right]$ satisfying
\begin{align*}
&v-\underline{v}\in C(\left[0, T_0\right] ; H^2(\Omega)), \\ 
&\bold{u}-\underline{\bold{u}}\in C(\left[0, T_0\right] ; H^2(\Omega)\cap L^2(\left[0, T_0\right] ; H^3(\Omega)),
\end{align*}
together with
\begin{equation*}
\lVert v-\underline{v} \rVert_{L^\infty(\left[0, T_0\right] ; H^2(\Omega))}+\lVert \bold{u}-\underline{\bold{u}} \rVert_{L^\infty(\left[0, T_0\right] ; H^2(\Omega))}\le M_1,
\end{equation*}
and
\beq\label{bddab}
\underline{\kappa}_1\le v(t,x)\le \overline{\kappa}_1, \quad \forall(t, x)\in \left[0, T_0\right]\times \Omega.    
\eeq
\end{proposition}

\subsection{Construction of shifts}
As desired, we will show the orbital stability of a composite wave of viscous shocks. More precisely, we will prove that a small $H^2$-perturbation of a composite wave of two viscous shocks is stable and uniformly converges to the composite wave up to shifts where each shock is shifted by $X_i(t)$ as follows:
\begin{align}
	\begin{aligned}\label{composite_wave}
	&(\tv^{X_1,X_2},\widetilde{u}^{X_1,X_2})(t,x) \\
	&\qquad := \left(\tv_1^{X_1}(x_1-\sigma_1t)+\tv_2^{X_2}(x_1-\sigma_2t)-v_m,\tu_{1}^{X_1}(x_1-\sigma_1t)+\tu_{2}^{X_2}(x_1-\sigma_2t)-u_{m}\right),
	\end{aligned}
\end{align}
where $f^{X_i}$ denotes a function $f$ shifted by $X_i$, that is, $f^{X_i}(x_1) := f(x_1-X_i(t))$ for any function $f$. This notation will be used throughout the paper. \\
We here introduce the shift functions explicitly, from which we could obtain a bound of derivative of shifts (at least locally-in-time) in Lemma \ref{lem:xex}, and obtain the desired a priori estimates  in Proposition \ref{prop:main}. 
We define a pair of shifts $(X_1, X_2)$ as a solution to the system of ODEs:
\begin{equation}\label{X(t)}
\left\{
\begin{array}{ll}
\di \dot{X}_1(t)=-\frac{M}{\delta_1}\bigg[
 \int_{\Omega}\frac{a^{X_1,X_2}}{\sigma_1^*}(\widetilde{h}_{1})_x^{X_1}(p(v)-p(\widetilde{v}^{X_1,X_2}))\,dx  -\int_{\Omega}a^{X_1,X_2} (p(\widetilde{v}^{X_1}_1))_x(v-\widetilde{v}^{X_1,X_2})\,dx\bigg],\\ [4mm]
 \di \dot{X}_2(t)=-\frac{M}{\delta_2}\bigg[
 \int_{\Omega}\frac{a^{X_1,X_2}}{\sigma_2^*}(\widetilde{h}_{2})_x^{X_2}(p(v)-p(\widetilde{v}^{X_1,X_2}))\,dx  -\int_{\Omega}a^{X_1,X_2} (p(\widetilde{v}^{X_2}_2))_x(v-\widetilde{v}^{X_1,X_2})\,dx\bigg],\\ [4mm]
\di  X_1(0)=X_2(0)=0,
\end{array}
\right.
\end{equation}

\noindent where $a^{X_1,X_2}$ is the shifted weight function as defined in \eqref{weight}, $\widetilde{h}_{i}:=\widetilde{u}_{i}-(2\mu+\lambda)\partial_{x_1}\widetilde{v}_i$, and $M$ is the specific constant chosen as $M:=\frac{5}{4}\sigma_m^4v_m^2\alpha_m$ with $\sigma_m=\sqrt{-p'(v_m)}$ and $\alpha_m=\frac{\gamma+1}{2\gamma \sigma_m p(v_m)}$. 
For well-posedness of the above ODEs in Lemma \ref{lem:xex}, we only need the following assumption for the shifted weight function at this point: $a^{X_1,X_2}$ is a $C^1$-function of $(x_1, X_1, X_2)$ with finite $C^1$-norm. This is verified by the explicit one defined in \eqref{weight}.\\
The following lemma ensures that \eqref{X(t)} has a unique Lipschitz continuous solution at least for the lifespan $[0,T_0]$ of solution $v$ satisfying \eqref{bddab}.
 We refer to \cite[Lemma 3.1]{HKK} for its proof.
\begin{lemma}\label{lem:xex}
For any $c_1,c_2>0$, there exists a constant $C>0$ such that the following is true.  For any $T>0$, and any   function $v\in L^\infty ((0,T)\times \R)$ 
verifying
\beq\label{odes}
c_1 \le v(t,x)\le c_2,\qquad \forall (t,x)\in [0,T]\times \bbr,
\eeq
 the system \eqref{X(t)} has a unique Lipschitz continuous solution $(X_1,X_2)$ on $[0,T]$. Moreover, 
\beq\label{roughx}
|X_1(t)| + |X_2(t)|  \le Ct,\quad \forall t\le T.
\eeq
\end{lemma}

\subsection{Main proposition for a priori estimates} We here present the main proposition for a priori estimates.

\begin{proposition} \label{prop:main}
For a given constant $U_+:=(v_+,\bold{u}_+)\in \mathbb{R_+}\times (\mathbb{R}\times\left\{ 0 \right\}\times\left\{ 0 \right\})$, there exist positive constants $\delta_0$, $C_0$, and $\varepsilon_1$ such that the following holds:   
For any constant states $U_m:=(v_m, u_m)$ and $U_-:=(v_-,u_-)$ satisfying \eqref{changed RH condition 2} 
satisfying $|p(v_-)-p(v_m)|=:\delta_1<\delta_0$ and $|p(v_m)-p(v_+)|=:\delta_2<\delta_0$, let $(\widetilde{v},\widetilde{u})$ denote the composite wave of two shifted shocks as in \eqref{com2}, where $(X_1,X_2)$ solves \eqref{X(t)}. Suppose that $(v,u)$ is the solution to \eqref{main} on $\left[ 0,T \right]$ for some $T>0$, and satisfy 
\begin{equation*}
\begin{aligned}
v-\widetilde{v}^{X_1,X_2}&\in C(\left[ 0,T \right]; H^2(\Omega)), \quad \nabla_x(v-\widetilde{v}^{X_1,X_2})\in L^2(\left[ 0,T \right]; H^1(\Omega)), \\    
\bold{u}-\widetilde{\bold{u}}^{X_1,X_2}&\in C(\left[ 0,T \right]; H^2(\Omega)), \quad \nabla_x(\bold{u}-\widetilde{\bold{u}}^{X_1,X_2})\in L^2(\left[ 0,T \right]; H^2(\Omega)),
\end{aligned}    
\end{equation*}
with
\begin{equation} \label{perturbation_small}
\lVert v-\widetilde{v}^{X_1,X_2} \rVert_{L^\infty(0,T;H^2(\Omega))}+\lVert \bold{u}-\widetilde{\bold{u}}^{X_1,X_2} \rVert_{L^\infty(0,T;H^2(\Omega))}\le \varepsilon_1,
\end{equation}
\end{proposition}
then for all $t\in\left[ 0,T \right]$,
\begin{equation}
\begin{aligned} \label{C-3}
&\lVert v-\widetilde{v}^{X_1,X_2} \rVert_{H^2(\Omega)}+\lVert \bold{u}-\widetilde{\bold{u}}^{X_1,X_2} \rVert_{H^2(\Omega)}+\sqrt{\int_{0}^{t}\sum\limits_{i=1}^2\delta_i|\dot{X}_i|^2\, d\tau} \\
&\quad \ \ \ +\sqrt{\int_{0}^{t}(\mathcal{G}^S(U)+\bold{D}(U)+\mathbf{D}_1(U)+\mathbf{D}_2(U))+\bold{D}_3(U))\, d\tau}\\
&\quad \le C_0(\lVert v_0-\widetilde{v}(0,\cdot) \rVert_{H^2(\Omega)}+\lVert \bold{u}_0-\widetilde{\bold{u}}(0,\cdot) \rVert_{H^2(\Omega)})+C_0\delta_0^{\frac{1}{4}}.
\end{aligned}
\end{equation}
In particular, for all $t\le T$,
\begin{equation} \label{bddx12}
|\dot{X}_1(t)|+|\dot{X}_2(t)|\le C_0\lVert (v-\widetilde{v}^{X_1,X_2})(t,\cdot) \rVert_{L^\infty(\Omega)}, 
\end{equation}
and
\beq\label{sx12}
	X_1(t)+\sigma_1t \le \frac{3\sigma_1+\sigma_2}{4}t < \frac{\sigma_1+\sigma_2}{2}t < \frac{\sigma_1+3\sigma_2}{4}t  \le X_2(t) +\sigma_2t.
 \eeq
\noindent Here, the constant $C_0$ is independent of $T$ and
\begin{equation}
\begin{aligned}
\mathcal{G}^S(U)&:=\sum\limits_{i=1}^2 \int_{\mathbb{T}^2} \int_{\mathbb{R}} |(\widetilde{v}_i)_{x_1}^{X_i}| |\Phi_i(v-\widetilde{v}^{X_1,X_2})|^2 \, dx_1 \, dx', \\
\mathbf{D}(U)&:=\int_{\mathbb{T}^2} \int_{\mathbb{R}} |\nabla_x(p(v)-p(\widetilde{v}^{X_1,X_2}))|^2 \, dx_1 \, dx' , \\
\mathbf{D}_1(U)&:=\int_{\mathbb{T}^2} \int_{\mathbb{R}} |\nabla_x(\bold{u}-\widetilde{\bold{u}}^{X_1,X_2})|^2 \, dx_1 \, dx',  \\
\mathbf{D}_2(U)&:=\int_{\mathbb{T}^2} \int_{\mathbb{R}} |\nabla^2_x(\bold{u}-\widetilde{\bold{u}}^{X_1,X_2})|^2 \, dx_1 \, dx' , \\
\mathbf{D}_3(U)&:=\int_{\mathbb{T}^2} \int_{\mathbb{R}} |\nabla^3_x(\bold{u}-\widetilde{\bold{u}}^{X_1,X_2})|^2 \, dx_1 \, dx', 
\end{aligned}
\end{equation}
where $\Phi_1$, $\Phi_2$ are cutoff functions defined by
 \begin{align} \label{def cutoff}
   \Phi_1(t,x_1):=
\begin{cases}
1, & \mbox{if } x_1<\frac{3\sigma_1+\sigma_2}{4}t, \\
0, & \mbox{if } x_1>\frac{\sigma_1+3\sigma_2}{4}t, \\
\text{linearly} \ \text{decreasing} \ 1 \ \text{to} \ 0, & \mbox{if } \frac{3\sigma_1+\sigma_2}{4}t\le x_1 \le \frac{\sigma_1+3\sigma_2}{4}t,
\end{cases} \Phi_2(t,x_1):=1-\Phi_1(t,x_1).
\end{align}

\begin{remark}\label{rem-phi}
Because the shifted composite wave $(\tv^{X_1,X_2}, \tu^{X_1,X_2})$ is not a solution to the Navier-Stokes equations, we have to control the interaction term of the 1- and 2-waves to have the desired results of Proposition \ref{prop:main}. For that, we localize the perturbation near each wave by applying the cutoff functions above. More precisely, $\phi_1$ (resp. $\phi_2$) localizes the perturbation near the 1-wave (resp. 2-wave) shifted by $X_1$ (resp. $X_2$). 
\end{remark}

\subsection{Global-in-time existence for perturbation}
Based on Proposition 3.1 and 3.2, we use the continuation argument to prove \eqref{thmstate1} for the global-in-time existence of perturbations. We can also use Proposition 3.2 to prove \eqref{thmstate2} for the long-time behavior. Those proofs are typical and use the same arguments as in \cite{HKK}. Therefore, we omit those details, and complete the proof of Theorem \ref{thm}. \\

Hence, the remaining part of the paper is dedicated to the proof of Proposition 3.2.\\

\subsection{Useful estimates}
We here present the useful estimates for the proofs.

\subsubsection{Sobolev inequalities} 
The following inequality is an extension of the 1D Poincare-type inequality \cite[Lemma 2.9]{Kang-V-NS17} to the multi-D domain $[0,1]\times\bbt^2$. 
\begin{lemma} \cite[Lemma 3.1]{WW}  \label{lem-poin}
	For any $f:\left[ 0, 1 \right]\times \mathbb{T}^2 \longrightarrow \mathbb{R}$ satisfying
\[\int_{\mathbb{T}^2} \int_{0}^{1} \left[ y_1(1-y_1)\left| \partial_{y_1}f \right|^2+\frac{\left| \nabla_{y'}f \right|^2}{y_1(1-y_1)} \right]  \, dy_1 \, dy'<\infty,\]
we have
\begin{equation} \label{poincare inequality}
\begin{aligned}
&\int_{\mathbb{T}^2} \int_{0}^{1} \left| f-\bar{f} \right|^2 \, dy_1 \, dy'\le \frac{1}{2} \int_{\mathbb{T}^2} \int_{0}^{1}  y_1(1-y_1)\left| \partial_{y_1}f \right|^2 \, dy_1 \, dy' \\
& \quad \quad \quad \quad \quad \quad \quad \quad \quad \quad \quad \ +\frac{1}{16\pi^2}\int_{\mathbb{T}^2} \int_{0}^{1} \frac{\left| \nabla_{y'}f \right|^2}{y_1(1-y_1)} \, dy_1 \, dy'.
\end{aligned}
\end{equation}
\end{lemma}

We will use the following interpolation inequality. We refer to \cite{adams2003sobolev} for its proof. 
\begin{lemma} \label{infty interpolation inequality lemma}
For any $g \in H^2(\Omega)$ where $ \Omega=\mathbb{R}\times{\mathbb{T}^2}$, we have
\[ 
\lVert g \rVert_{L^{\infty}(\Omega)}\le C\lVert g \rVert^{\frac{1}{2}}_{L^2{(\Omega})}\lVert \partial_{x_1}g \rVert^{\frac{1}{2}}_{L^2{(\Omega)}}+C\lVert \nabla_xg \rVert^{\frac{1}{2}}_{L^2{(\Omega)}}\lVert \nabla_x^2g \rVert^{\frac{1}{2}}_{L^2{(\Omega)}}.    
\]  
\end{lemma}

\subsubsection{Estimate on the relative quantities}
Our approach is based on the relative entropy method, which was first employed by Dafermos \cite{Dafermos1} and Diperna \cite{DiPerna} fore the $L^2$ stability and uniqueness of Lipschitz solution to the hyperbolic conservation laws.
We here present some useful estimates for the relative quantities, such as
\[p(v|w) = p(v)-p(w)-p'(w)(v-w),\quad Q(v|w) = Q(v)-Q(w)-Q'(w)(v-w).\]
To use the Taylor-Expansion, these must be almost locally quadratic quantities. The estimates below show the precise estimate on this local behavior of the relative quantities. The proof of the lemma can be found in \cite{Kang-V-NS17}.

\begin{lemma}\label{lem-rel-quant}
	Let $\gamma>1$ and $v_+$ be given constants. Then there exist constants $C,\delta_*>0$ such that the following assertions hold:
	\begin{enumerate}
		\item For any $v,w$ satisfying $0<w<2v_+$ and $0<v<3v_+$,
		\begin{equation}\label{est-rel-1}
		|v-w|^2\le CQ(v|w),\quad |v-w|^2\le Cp(v|w).
		\end{equation}
		\item For any $v,w$ satisfying $v,w>v_+/2$,
		\begin{equation}\label{est-rel-2}
		|p(v)-p(w)|\le C|v-w|.
		\end{equation}
		\item For any $0<\delta<\delta_*$ and any $(v,w)\in\bbr_+^2$ satisfying $|p(v)-p(w)|<\delta$ and $|p(w)-p(v_+)|<\delta$,
		\begin{align}
		\begin{aligned}\label{est-rel-3}
		&p(v|w) \le \left(\frac{\gamma+1}{2\gamma}\frac{1}{p(w)}+C\delta\right)|p(v)-p(w)|^2,\\
		&Q(v|w) \ge \frac{p(w)^{-\frac{1}{\gamma}-1}}{2\gamma}|p(v)-p(w)|^2-\frac{1+\gamma}{3\gamma^2}p(w)^{-\frac{1}{\gamma}-2}(p(v)-p(w))^3,\\ 
		&Q(v|w)\le \left(\frac{p(w)^{-\frac{1}{\gamma}-1}}{2\gamma}+C\delta\right)|p(v)-p(w)|^2.
		\end{aligned}
		\end{align}
	\end{enumerate}
\end{lemma}
\subsubsection{Estimate on shifts} 
Here, we prove the estimates \eqref{bddx12} and \eqref{sx12}. To use the Lemma \ref{lem-rel-quant}, choose $\varepsilon_1$ and $\delta_0$ such that $\varepsilon_1, \delta_0 \in (0,\delta_*)$, where $\delta_*$ is the constant in Lemma \ref{lem-rel-quant}. \\
Using the assumption \eqref{perturbation_small}, Sobolev inequality, and \eqref{est-rel-2}, we have
\beq\label{smp1}
\|p(v)-p(\tv)\|_{L^\infty((0,T)\times\Omega)}\le C\|v-\tv\|_{L^\infty((0,T)\times\Omega)} \le C\varepsilon_1.
\eeq
Note that from \eqref{viscous-shock-h}, we know for each $i=1, 2$,
\begin{equation}\label{hp}
\sigma_i^*(\widetilde{h}_{i})' =p(\tv_i)'.    
\end{equation} 
Thus, applying \eqref{smp1} and \eqref{hp} to the definition of $X_i$ \eqref{X(t)}, we get for each $i=1, 2$,
\beq\label{dxbound}
|\dot{X}_i(t)| \le \frac{C}{\delta_i}  \left(\|p(v)-p(\tv) \|_{L^\infty(\Omega)}+\|v-\tv \|_{L^\infty(\Omega)}\right) \int_\bbr |(\tv_i)_{x_1}^{X_i} | dx_1 \le C\|v-\tv\|_{L^\infty(\Omega)},\quad t\le T.
\eeq
This completes the proof of \eqref{bddx12}. 

\noindent On the other hand, since $\sigma_2-\sigma_1>0$, taking $\varepsilon_1$ small enough together with \eqref{dxbound} and \eqref{smp1} we can get
\begin{equation} \label{bddx13}
|\dot{X}_i(t)|\le \frac{\sigma_2-\sigma_1}{8}, \quad i=1,2, \quad t\le T.    
\end{equation}
Integrating \eqref{bddx13} over $\left[0, t\right]$ for any $t\le T$, we can get 
\begin{equation} \label{bddx14}
X_1(t) \le \frac{\sigma_2-\sigma_1}{4}t, \quad X_2(t)\ge \frac{\sigma_1-\sigma_2}{4}t,
\end{equation}
and especially,
\beq\label{fsx12}
X_1(t)+\sigma_1t\le \frac{3\sigma_1+\sigma_2}{4}t, \quad X_2(t)+\sigma_2t\ge \frac{\sigma_1+3\sigma_2}{4}t.
\eeq
Because of \eqref{fsx12}, this completes the proof of \eqref{sx12}. Thanks to the estimates above, we can prove Lemmas \ref{shock interaction lemma1} and \ref{shock interaction lemma2} on the interaction estimates.
\subsubsection{Interaction estimates}
We present Lemmas \ref{shock interaction lemma1} and \ref{shock interaction lemma2} on interaction estimates. Those proofs are similar to \cite[Lemma 3.2]{HKK23} and \cite[Lemma 3.3]{HKK23}.
\begin{lemma} \cite[Lemma 3.2]{HKK23} \label{shock interaction lemma1}
Given $v_+>0$, there exist positive constant $\delta_0$, $C$ such that for any $\delta_1,\delta_2\in (0, \delta_0)$, the following estimates hold. For each $i=1, 2$,
\[|(\widetilde{v}_i)^{X_i}_{x_1}||\widetilde{v}^{X_1, X_2}-\widetilde{v}_i^{X_i}|
\le C\delta_i\delta_1\delta_2 \exp(-C\min \left\{ \delta_1,\delta_2 \right\}t), \quad t>0, \quad x_1\in \mathbb{R}, \]
\[\int_{\mathbb{R}}|(\widetilde{v}_i)^{X_i}_{x_1}||\widetilde{v}^{X_1, X_2}-\widetilde{v}_i^{X_i}|dx_1
\le C\delta_1\delta_2 \exp(-C\min \left\{ \delta_1,\delta_2 \right\}t), \quad t>0,\]
\[\int_{\mathbb{R}}|(\widetilde{v}_1)^{X_1}_{x_1}||(\widetilde{v}_2)^{X_2}_{x_1}|dx_1
\le C\delta_1\delta_2 \exp(-C\min \left\{ \delta_1,\delta_2 \right\}t), \quad t>0.\]
\end{lemma}
\begin{lemma} \cite[Lemma 3.3]{HKK23} \label{shock interaction lemma2}
Let $\Phi_i$ be the functions defined in \eqref{def cutoff}. Given $v_+>0$, there exist positive constant $\delta_0$, $C$ such that for any $\delta_1,\delta_2\in (0, \delta_0)$, the following estimates hold.
\[\Phi_2|(\widetilde{v}^{X_1}_1)_{x_1}|\le C\delta_1^2\exp(-C\delta_1t), \quad \Phi_1|(\widetilde{v}^{X_2}_2)_{x_1}|\le C\delta_2^2\exp(-C\delta_2t), \quad t>0, \quad x_1\in \mathbb{R},\]
\[\int_{\mathbb{R}}\Phi_2|(\widetilde{v}^{X_1}_1)_{x_1}|dx_1\le C\delta_1\exp(-C\delta_1t), \quad \int_{\mathbb{R}}\Phi_1|(\widetilde{v}^{X_2}_2)_{x_1}|\le C\delta_2\exp(-C\delta_2t), \quad t>0.\]
\end{lemma}

\subsection{Notations} In what follows, we use the following notations for simplicity. \\
1. $C$ denotes a positive $O(1)$-constant that may change from line to line, but is independent of the small constants $\delta_0, \eps_1, \delta_1,\delta_2$, $\lambda_1,\lambda_2$ (to be introduced below) and the time $T$.\\
2. We omit the dependence on the pair of shifts $(X_1,X_2)$ of the composite wave \eqref{composite_wave} without confusion as:
\[
(\tv, \tu) (t,x) := (\tv^{X_1,X_2}, \tu^{X_1,X_2}) (t,x).
\]
\section{Estimates on weighted relative entropy with shifts}\label{sec:4}
\setcounter{equation}{0}
Using the $a$-contraction method with shifts, we can get the bounds of perturbations in \eqref{C-3}. For simplicity of our analysis, we employ the effective velocity $\bold{h}:=\bold{u}-\left( 2\mu+\lambda \right)\nabla_xv$ as in \cite{KVW-2022,Kang-V-NS17,KV-Inven}, which is related to the BD entropy \cite{BDL, BD_03,BD_06}).
With this effective velocity, we transform the Navier-Stokes equations \eqref{changed main} to the system:
\begin{equation} \label{4.1}
\begin{cases}
\rho\left(\partial_tv+\bold{u}\cdot\nabla_xv\right)-div_x\bold{h}=(2\mu+\lambda)\Delta_xv, \\
\rho\left(\partial_t\bold{h}+\bold{u}\nabla_x\bold{h}\right)+\nabla_xp(v)=R,
\end{cases}
\end{equation}
where
\begin{equation} \label{presentR}
R=\frac{2\mu+\lambda}{v}(\nabla_x\bold{u} \nabla_xv - div_x\bold{u}\nabla_xv)-\mu \nabla_x \times \nabla_x \times \bold{u}.
\end{equation}
Then it follows from \eqref{changed viscous shock equation} - \eqref{change shock equation 2} that the associated viscous shocks $\widetilde{v}_i$ and
\ $\widetilde{h}_{i}:=\widetilde{u}_{i}-(2\mu+\lambda)\partial_{x_1}\widetilde{v}_i$ satisfy the following ODEs:
\begin{equation}\label{viscous-shock-h}
\begin{cases}
\widetilde{\rho}_i^{X_i}(-\sigma_i(\widetilde{v}_i)_{x_1}^{X_i}+\widetilde{u}_{i}^{X_i}(\widetilde{v}_i)_{x_1}^{X_i})-(\widetilde{h}_{i})_{x_1}^{X_i}=(2\mu+\lambda)(\widetilde{v}_i)_{x_1x_1}^{X_i}, \quad (i=1,2)\\
\widetilde{\rho}_i^{X_i}((-\sigma_i+\widetilde{u}_{i}^{X_i})(\widetilde{h}_{i})_{x_1}^{X_i})+p(\widetilde{v}_i^{X_i})_{x_1}=0, \quad (i=1,2)\\
(\widetilde{v}_1,\widetilde{h}_{1})(-\infty) = (v_-,u_{-}),\quad (\tv_1,\widetilde{h}_{1})(+\infty) = (v_m,u_{m}),\\
(\tv_2,\widetilde{h}_{2})(-\infty) = (v_m,u_{m}),\quad (\tv_2,\widetilde{h}_{2})(+\infty) = (v_+,u_{+}).
\end{cases}
\end{equation} 
Let $(\tv,\widetilde{h})$ denote the shifted composite wave such that
\[
(\tv,\widetilde{h})(t,x) := \left(\tv_1^{X_1}(x_1-\sigma_1t)+\tv_2^{X_2}(x_1-\sigma_2t)-v_m,\widetilde{h}_{1}^{X_1}(x_1-\sigma_1t)+\widetilde{h}_{2}^{X_2}(x_1-\sigma_2t)-u_{m}\right).
\]
In addition, we denote
\[\widetilde{\bold{h}}_1(t,x):=(\widetilde{h}_{1}^{X_1}(x_1-\sigma_1t), 0, 0), \quad \widetilde{\bold{h}}_2(t,x):=(\widetilde{h}_{2}^{X_2}(x_1-\sigma_2t), 0, 0).\]
It follows from \eqref{4.1} and \eqref{viscous-shock-h} that
\begin{equation} \label{mixed equations}
\begin{cases}
\rho(v-\widetilde{v})_t+\rho \bold{u}\cdot \nabla_x(v-\widetilde{v})-div_x(\bold{h}-\widetilde{\bold{h}})-\rho\sum\limits_{i=1}^2 \dot{X}_i(t)(\widetilde{v}_i)_{x_1}^{X_i}+\sum\limits_{i=1}^2 F_i(\widetilde{v}_i)_{x_1}^{X_i}=(2\mu+\lambda)\Delta_x(v-\widetilde{v}), \\
\rho(\bold{h}-\widetilde{\bold{h}})_t+\rho \bold{u} \nabla_x(\bold{h}-\widetilde{\bold{h}})+\nabla_x(p(v)-p(\widetilde{v}))-\rho\sum\limits_{i=1}^2 \dot{X}_i(t)(\widetilde{\bold{h}}_i)_{x_1}+\sum\limits_{i=1}^2F_i(\widetilde{\bold{h}}_i)_{x_1}=R^*,
\end{cases}
\end{equation}
where for each $i=1,2$,
\begin{align}
\begin{aligned} \label{presentF}
F_i&=-\sigma_i(\rho-\widetilde{\rho}_i^{X_i})+\rho u_1 - \widetilde{\rho}_i^{X_i} \widetilde{u}_{i}^{X_i}\\
&=-\frac{\sigma_i^*}{\widetilde{\rho}_i^{X_i}}(\rho-\widetilde{\rho}_i^{X_i})+\rho(u_1-\widetilde{u}_{i}^{X_i}),\\
\end{aligned}
\end{align}
and, 
\begin{equation} \label{presentR^*}
R^*=R-\nabla_x( p(\widetilde{v})-p(\widetilde{v}_1^{X_1})-p(\widetilde{v}_2^{X_2}) ).
\end{equation}

\begin{lemma} \label{lemma4.1}
There exists a constant $C>0$ such that for any $t\in\left[0, T\right]$,
\begin{align}
\begin{aligned} \label{lemma4.1eq}
&\int_{\mathbb{T}^2} \int_{\mathbb{R}} \left( \frac{|\bold{h}-\widetilde{\bold{h}}|^2}{2}+Q(v|\widetilde{v}) \right) \, dx_1 \, dx'+\int_{0}^{t} \left( \sum\limits_{i=1}^2\delta_i|\dot{X}_i(\tau)|^2+\mathbf{G}_1(\tau)+\mathbf{G}_3(\tau)+\mathcal{G}^S(\tau)+\mathbf{D}(\tau) \right) \, d\tau \\
&\le C\int_{\mathbb{T}^2}\int_{\mathbb{R}}\left( \frac{|\bold{h}_0(x)-\widetilde{\bold{h}}(0,x)|^2}{2}+Q(v_0(x) | \widetilde{v}(0,x)) \right) \,dx_1 \, dx'+C(\delta_0+\varepsilon_1)\int_{0}^{t}||\nabla_x(\bold{u}-\widetilde{\bold{u}})||^2_{H^1}\, d\tau+C\delta_0,
\end{aligned}
\end{align}
where
\begin{equation*}
\begin{aligned}
\mathbf{G}_1&:=\sum\limits_{i=1}^2 \int_{\mathbb{T}^2}  \int_{\mathbb{R}} |(a_i)_{x_1}^{X_i}| \left| h_1-\widetilde{h}-\frac{p(v)-p(\widetilde{v})}{\sigma_i^*} \right|^2\, dx_1 \, dx',\\  
\mathbf{G}_3&:=\sum\limits_{i=1}^2 \int_{\mathbb{T}^2}  \int_{\mathbb{R}} |(a_i)_{x_1}^{X_i}| (h_2^2+h_3^2)\, dx_1 \, dx',\\
\mathcal{G}^S&:=\sum_{i=1}^2 \int_{\mathbb{T}^2}\int_\R |(\widetilde{v}_i)^{X_i}_{x_1}| |\Phi_i(p(v)-p(\widetilde{v})) |^2 \,dx_1\,dx',\\
\mathbf{D}&:=\int_{\mathbb{T}^2}  \int_{\mathbb{R}} |\nabla_x(p(v)-p(\widetilde{v}))|^2\, dx_1 \, dx'.
\end{aligned}    
\end{equation*}
\end{lemma}

\subsection{Construction of weight functions}
In order to deal with the two shock waves, we first introduce two weight functions $a_1, a_2$ associated with 1-shock and 2-shock respectively: for each $i=1,2,$ we define
\begin{equation} \label{weightdef}
a_i(x_1-\sigma_i t) = 1+\frac{\nu_i(p(v_m)-p(\widetilde{v_i}(x_1-\sigma_it)))}{\delta_i},
\end{equation}
where $\delta_i$ is the $i$th-shock strength. At this time, for each $i=1, 2$, $\nu_i$ is a small constant satisfying $\delta_i\ll\nu_i<C\sqrt{\delta_i}$ for so that it is large enough compared to the $i$th-shock strength. 
In addition, we define $\nu:=\nu_1+\nu_2$, which satisfies $\nu \le 2\max(\nu_1,\nu_2)\le C\sqrt{\delta_0}$.\\
Note that $\frac{3}{4} \le 1-\nu \le a_i \le 1+\nu \le \frac{5}{4}$ and
\beq\label{avd}
(a_i)_{x_1} = -\frac{\nu_i}{\delta_i}(p(\widetilde{v_i}))_{x_1},
\eeq
from which we have
\beq\label{inta}
|(a_i)_{x_1}| \sim \frac{\nu_i}{\delta_i} |(\widetilde{v}_i)_{x_1}|,\quad \mbox{and so,}\quad  \|(a_i)_{x_1}\|_{L^\infty(\bbr)}\le \nu_i\delta_i,\quad  \|(a_i)_{x_1}\|_{L^1(\bbr)} =\nu_i.
\eeq
To deal with the shifted composite wave, we think about the following composition of shifted weight functions as
\begin{equation}\label{weight}
	a^{X_1,X_2}(t,x_1):=a_1^{X_1}(x_1-\sigma_1 t) +a_2^{X_2}(x_1-\sigma_2t)-1=a_1(x_1-\sigma_1t-X_1(t))+a_2(x_1-\sigma_2t-X_2(t))-1.
\end{equation}
For notational simplicity as before, we will omit the dependence on shifts:
\[
a(t,x_1) := a^{X_1,X_2}(t,x_1).
\]

\subsection{Evolution of the weighted relative entropy} 
We present the representation for the evolution of the relative entropy. The computation is typical but we present its proof in Appendix for the readers' convenience.   

\begin{lemma}\label{lem:rel-ent}
	Let $a$ be the weighted function defined by \eqref{weightdef} and \eqref{weight}. It holds
\begin{align}
	\begin{aligned}\label{est-rel-ent}
	&\frac{d}{dt}\int_{\mathbb{T}^2}\int_{\mathbb{R}} a \rho \left(Q(v|\tv)+\frac{|\bold{h}-\widetilde{\bold{h}}|^2}{2}\right) \,dx_1 \,dx'=\sum\limits_{i=1}^2\left( \dot{X}_i(t)Y_i(t) \right)+\mathcal{J}^{\rm bad}(U)-\mathcal{J}^{\rm good}(U),
	\end{aligned}
\end{align}
where
\begin{align}
\begin{aligned} \label{before diffusionterm}
&\sum\limits_{i=1}^2\left( \dot{X}_i(t)Y_i(t) \right)=-\int_{\mathbb{T}^2}\int_{\mathbb{R}}\rho Q(v|\widetilde{v})\sum\limits_{i=1}^2\dot{X}_i(t)(a_i)_{x_1}^{X_i}dx_1dx'-\int_{\mathbb{T}^2}\int_{\mathbb{R}}a\rho p'(\widetilde{v})(v-\widetilde{v})\sum\limits_{i=1}^2\dot{X}_i(t)(\widetilde{v}_i)_{x_1}^{X_i}dx_1dx', \\ 
&\qquad \quad \mathcal{J}^{\rm bad}(U)=\int_{\mathbb{T}^2}\int_{\mathbb{R}}Q(v|\widetilde{v})\sum\limits_{i=1}^2 F_i (a_i)_{x_1}^{X_i}dx_1dx'-\int_{\mathbb{T}^2}\int_{\mathbb{R}}ap(v|\widetilde{v})\sum\limits_{i=1}^2 F_i (\widetilde{v}_i)_{x_1}^{X_i}dx_1dx'\\
&\qquad \qquad \qquad \qquad +\int_{\mathbb{T}^2}\int_{\mathbb{R}}ap(v|\widetilde{v})\sum\limits_{i=1}^2 \sigma_i^* (\widetilde{v}_i)_{x_1}^{X_i}dx_1dx'+\int_{\mathbb{T}^2}\int_{\mathbb{R}}a(p(v)-p(\widetilde{v}))\sum\limits_{i=1}^2 F_i(\widetilde{v}_i)_{x_1}^{X_i}dx_1dx' \\ 
&\qquad \qquad \qquad \qquad -(2\mu+\lambda)\int_{\mathbb{T}^2}\int_{\mathbb{R}}a\partial_{x_1}(p(v)-p(\tilde{v})) \partial_{x_1}p(\widetilde{v})\left( \frac{1}{\gamma p(v)^{1+\frac{1}{\gamma}}}- \frac{1}{\gamma p(\widetilde{v})^{{1+\frac{1}{\gamma}}}}\right)dx_1dx'\\
&\qquad \qquad \qquad \qquad -(2\mu+\lambda)\int_{\mathbb{T}^2}\int_{\mathbb{R}}(p(v)-p(\widetilde{v}))\frac{\partial_{x_1}(p(v)-p(\widetilde{v}))}{\gamma p(v)^{1+\frac{1}{\gamma}}}\sum\limits_{i=1}^2 (a_i)_{x_1}^{X_i}dx_1dx'\\ 
&\qquad \qquad \qquad \qquad -(2\mu+\lambda)\int_{\mathbb{T}^2}\int_{\mathbb{R}}(p(v)-p(\widetilde{v}))\partial_{x_1}p(\widetilde{v})\left( \frac{1}{\gamma p(v)^{1+\frac{1}{\gamma}}}- \frac{1}{\gamma p(\widetilde{v})^{{1+\frac{1}{\gamma}}}}\right) \sum\limits_{i=1}^2 (a_i)_{x_1}^{X_i}dx_1dx', \\ 
&\qquad \ \mathcal{J}^{\rm good}(U)=\int_{\mathbb{T}^2}\int_{\mathbb{R}}Q(v|\widetilde{v})\sum\limits_{i=1}^2 \sigma_i^* (a_i)_{x_1}^{X_i}dx_1dx'+(2\mu+\lambda)\int_{\mathbb{T}^2}\int_{\mathbb{R}}a\frac{|\nabla_x(p(v)-p(\widetilde{v}))|^2}{\gamma p(v)^{1+\frac{1}{\gamma}}}dx_1dx'
\end{aligned}
\end{align}
\end{lemma}
\begin{remark}
Since $\sigma_i^{*}(a_i^{X_i})_{x_1}>0$ and $a>\frac{1}{2}$, $-\mathcal{G}(t)$ consists of five terms with good sign, while $\mathcal{B}(t)$ consists of bad terms.
\end{remark}
\begin{remark}
We note that the definition of the weight $a_i$ implies
\[\sigma_i^* (a_i)^{X_i}_{x_1} = -\frac{\nu_i\sigma_i^*}{\delta_i}(p(\widetilde{v}_i^{X_i}))_{x_1} = \frac{\gamma\nu_i\sigma_i^*}{\delta_i}(\widetilde{v}_i^{X_i})^{-\gamma-1}(\widetilde{v}_i)_{x_1}^{X_i}>0,\]
from which we observe that $\mathcal{J}^{\rm good}$ consists of good terms.
\end{remark}

\subsection{Maximization on $\bold{h}-\widetilde{\bold{h}}$}
Among the terms in $\mathcal{J}^{\rm bad}(U)$, a main bad term is
\[\int_{\mathbb{T}^2}\int_{\bbr}(a_i)^{X_i}_{x_1}(p(v)-p(\widetilde{v}))(h_1-\widetilde{h})\,dx_1\,dx'\] 
where the perturbations for $p(v)$ and $h_1$ are multiplied. To exploit the parabolic term on $v$-variable and so use the Poincar\'e-type inequality, we separate $h_1-\widetilde{h}$ from $p(v)-p(\widetilde{v})$ by using the quadratic structure of $\bold{h}-\widetilde{\bold{h}}$. Exactly, using
\begin{align*}
	(a_i)^{X_i}_{x_1}&(p(v)-p(\widetilde{v}))(h_1-\widetilde{h})-\frac{\sigma_i^*}{2}(a_i)^{X_i}_{x_1}|\bold{h}-\widetilde{\bold{h}}|^2\\
	&=-\frac{\sigma_i^*(a_i)^{X_i}_{x_1}}{2}\left|h_1-\widetilde{h}-\frac{p(v)-p(\widetilde{v})}{\sigma_i^*}\right|^2+\frac{(a_i)_{x_1}^{X_i}}{2\sigma_i^*}|p(v)-p(\widetilde{v})|^2-\sigma_i^* (a_i)_{x_1}^{X_i} \left( \frac{h_2^2+h_3^2}{2} \right),
\end{align*}
we can rewrite the terms $\mathcal{J}^{\rm bad}(U)-\mathcal{J}^{\rm good}(U)$ in the above lemma \eqref{est-rel-ent} as 
\[\mathcal{J}^{\text{bad}}(U) - \mathcal{J}^{\text{good}}(U) = \mathcal{B}(U)-\mathcal{G}(U),\]
where
\begin{align*}
\mathcal{B}(U)& =\sum\limits_{i=1}^2  \Bigg[ \frac{1}{2\sigma_i^*} \int_{\mathbb{T}^2} \int_{\mathbb{R}} (a_i)_{x_1}^{X_i} |p(v)-p(\widetilde{v})|^2 \, dx_1 \, dx'+\sigma_i^* \int_{\mathbb{T}^2}  \int_{\mathbb{R}}  ap(v|\widetilde{v}) (\widetilde{v}_i)_{x_1}^{X_i} \, dx_1 \, dx'   \,\\
 & \qquad \quad \ +\int_{\mathbb{T}^2}  \int_{\mathbb{R}} ap'(\widetilde{v})(v-\widetilde{v})F_i (\widetilde{v}_i)_{x_1}^{X_i}-a(h_1-\widetilde{h})F_i(\widetilde{h}_{i})_{x_1}^{X_i} \, dx_1 \, dx' \\
 & \qquad \quad \ +\int_{\mathbb{T}^2}  \int_{\mathbb{R}}  \left( Q(v|\widetilde{v})+\frac{|\bold{h}-\widetilde{\bold{h}}|^2}{2} \right) F_i (a_i)_{x_1}^{X_i}\, dx_1 \, dx' \\
 & \qquad \quad \ -\frac{(2\mu+\lambda)}{2}  \int_{\mathbb{T}^2}  \int_{\mathbb{R}} a\partial_{x_1}(p(v)-p(\widetilde{v})) \partial_{x_1}p(\widetilde{v})\left( \frac{1}{\gamma p(v)^{1+\frac{1}{\gamma}}}- \frac{1}{\gamma p(\widetilde{v})^{{1+\frac{1}{\gamma}}}}\right) \, dx_1 \, dx' \\ 
 & \qquad \quad \ -(2\mu+\lambda)  \int_{\mathbb{T}^2}  \int_{\mathbb{R}} (a_i)_{x_1}^{X_i}(p(v)-p(\widetilde{v}))\frac{\partial_{x_1}(p(v)-p(\widetilde{v}))}{\gamma p(v)^{1+\frac{1}{\gamma}}} \, dx_1 \, dx' \\ 
 & \qquad \quad \ -(2\mu+\lambda)\int_{\mathbb{T}^2}  \int_{\mathbb{R}} (a_i)_{x_1}^{X_i}(p(v)-p(\widetilde{v}))\partial_{x_1}p(\widetilde{v})\left( \frac{1}{\gamma p(v)^{1+\frac{1}{\gamma}}}- \frac{1}{\gamma p(\widetilde{v})^{{1+\frac{1}{\gamma}}}}\right)  \, dx_1 \, dx' \Bigg], \\ 
 &
\end{align*} 
\begin{align*}
\mathcal{G}(U)& =\sum\limits_{i=1}^2 \Bigg[ \frac{\sigma_i^*}{2} \int_{\mathbb{T}^2}  \int_{\mathbb{R}} (a_i)_{x_1}^{X_i} \left| h_1-\widetilde{h}-\frac{p(v)-p(\widetilde{v})}{\sigma_i^*} \right|^2\, dx_1 \, dx'   \,\\
&\qquad \quad \ +\sigma_i^* \int_{\mathbb{T}^2}  \int_{\mathbb{R}} Q(v|\widetilde{v}) (a_i)_{x_1}^{X_i}\, dx_1 \, dx' \\
&\qquad \quad \ +\frac{\sigma_i^*}{2} \int_{\mathbb{T}^2}  \int_{\mathbb{R}} (a_i)_{x_1}^{X_i} (h_2^2+h_3^2)\, dx_1 \, dx' \Bigg], \\
&\qquad \quad \ +(2\mu+\lambda) \int_{\mathbb{T}^2}  \int_{\mathbb{R}} a\frac{|\nabla_x(p(v)-p(\widetilde{v}))|^2}{\gamma p(v)^{1+\frac{1}{\gamma}}}\, dx_1 \, dx'.
\end{align*} 
Thus, we estimate the right-hand side of the below equation: 
\begin{equation}\label{est-1}
	\frac{d}{dt}\int_{\mathbb{T}^2}\int_{\bbr}a\rho\eta(U|\widetilde{U})\,dx_1\,dx' = \sum_{i=1}^2(\dot{X}_iY_i(U)) + \mathcal{B}(U)-\mathcal{G}(U).
\end{equation}

\subsection{Decompositions}
First, we name each term of $\mathcal{B}(U)$ and $\mathcal{G}(U)$ as follows:
\begin{align*}
	&\mathcal{B}(U) =\sum_{i=1}^8 \mathcal{B}_i(U),\\
	&\mathcal{G}(U) = \mathcal{G}_1(U) + \mathcal{G}_2(U) + \mathcal{G}_3(U)  + \mathcal{D}(U),
\end{align*}
where
\begin{align*}
 & \mathcal{B}_1(U):=\sum\limits_{i=1}^2 \frac{1}{2\sigma_i^*} \int_{\mathbb{T}^2} \int_{\mathbb{R}} (a_i)_{x_1}^{X_i} |p(v)-p(\widetilde{v})|^2 \, dx_1 \, dx'\\ 
 & \mathcal{B}_2(U):=\sum\limits_{i=1}^2 \sigma_i^* \int_{\mathbb{T}^2}  \int_{\mathbb{R}}  ap(v|\widetilde{v}) (\widetilde{v}_i)_{x_1}^{X_i} \, dx_1 \, dx'\\
 & \mathcal{B}_3(U):=\sum\limits_{i=1}^2  \int_{\mathbb{T}^2}  \int_{\mathbb{R}} ap'(\widetilde{v})(v-\widetilde{v})F_i (\widetilde{v}_i)_{x_1}^{X_i}-a(h_1-\widetilde{h})F_i(\widetilde{h}_{i})_{x_1}^{X_i} \, dx_1 \, dx' \\
 & \mathcal{B}_4(U):=\sum\limits_{i=1}^2 \int_{\mathbb{T}^2}  \int_{\mathbb{R}}  \left( Q(v|\widetilde{v})+\frac{|\bold{h}-\widetilde{\bold{h}}|^2}{2} \right) F_i (a_i)_{x_1}^{X_i}\, dx_1 \, dx' \\
 & \mathcal{B}_5(U):=-(2\mu+\lambda)  \int_{\mathbb{T}^2}  \int_{\mathbb{R}} a\partial_{x_1}(p(v)-p(\widetilde{v})) \partial_{x_1}p(\widetilde{v})\left( \frac{1}{\gamma p(v)^{1+\frac{1}{\gamma}}}- \frac{1}{\gamma p(\widetilde{v})^{{1+\frac{1}{\gamma}}}}\right) \, dx_1 \, dx' \\ 
 & \mathcal{B}_6(U):=-(2\mu+\lambda)  \sum\limits_{i=1}^2  \int_{\mathbb{T}^2}  \int_{\mathbb{R}} (a_i)_{x_1}^{X_i}(p(v)-p(\widetilde{v}))\frac{\partial_{x_1}(p(v)-p(\widetilde{v}))}{\gamma p(v)^{1+\frac{1}{\gamma}}} \, dx_1 \, dx' \\ 
 & \mathcal{B}_7(U):=-(2\mu+\lambda)\sum\limits_{i=1}^2 \int_{\mathbb{T}^2}  \int_{\mathbb{R}} (a_i)_{x_1}^{X_i}(p(v)-p(\widetilde{v}))\partial_{x_1}p(\widetilde{v})\left( \frac{1}{\gamma p(v)^{1+\frac{1}{\gamma}}}- \frac{1}{\gamma p(\widetilde{v})^{{1+\frac{1}{\gamma}}}}\right)  \, dx_1 \, dx' \\ 
 & \mathcal{B}_8(U):= \int_{\mathbb{T}^2}  \int_{\mathbb{R}}a(\bold{h}-\widetilde{\bold{h}})\cdot R^*\, dx_1 \, dx', \quad \quad (R^* \ was \ defined \ in \ \eqref{presentR} \ and \ \eqref{presentR^*})
\end{align*}
and
\begin{align*}
& \mathcal{G}_1(U):=\sum\limits_{i=1}^2 \frac{\sigma_i^*}{2} \int_{\mathbb{T}^2}  \int_{\mathbb{R}} (a_i)_{x_1}^{X_i} \left| h_1-\widetilde{h}-\frac{p(v)-p(\widetilde{v})}{\sigma_i^*} \right|^2\, dx_1 \, dx' \\
& \mathcal{G}_2(U):=\sum\limits_{i=1}^2 \sigma_i^* \int_{\mathbb{T}^2}  \int_{\mathbb{R}} (a_i)_{x_1}^{X_i}Q(v|\widetilde{v}) \, dx_1 \, dx' \\
& \mathcal{G}_3(U):=\sum\limits_{i=1}^2 \frac{\sigma_i^*}{2} \int_{\mathbb{T}^2}  \int_{\mathbb{R}} (a_i)_{x_1}^{X_i} (h_2^2+h_3^2)\, dx_1 \, dx' \\
& \mathcal{D}(U):=(2\mu+\lambda) \int_{\mathbb{T}^2}  \int_{\mathbb{R}} a\frac{|\nabla_x(p(v)-p(\widetilde{v}))|^2}{\gamma p(v)^{1+\frac{1}{\gamma}}}\, dx_1 \, dx'.
\end{align*}

For each $Y_i$ in \eqref{before diffusionterm}, we initially write $Y_i$ more explicitly as follows:
\begin{align*}
	Y_i(U)
	& = -\int_{\mathbb{T}^2}  \int_{\mathbb{R}} \rho(a_i)_{x_1}^{X_i} \left( Q(v|\widetilde{v})+\frac{|\bold{h}-\widetilde{\bold{h}}|^2}{2} \right)  \, dx_1 \, dx'  +\int_{\mathbb{T}^2}  \int_{\mathbb{R}} a\rho (\widetilde{h}_{i})_{x_1}^{X_i}(h_1-\widetilde{h}_1) \, dx_1 \, dx'\\
	&\quad -\int_{\mathbb{T}^2}  \int_{\mathbb{R}} a\rho p'(\widetilde{v})(\widetilde{v}_i)_{x_1}^{X_i}(v-\widetilde{v}) \, dx_1 \, dx'.
\end{align*}
Because an essential idea in our analysis is to apply Poincar\'e-type inequality of Lemma \ref{lem-poin} by extracting a good term on an average of perturbation $p(v)-p(\tv)$ from the shift part $\sum_{i=1}^2(\dot{X}_iY_i(U))$, we decompose $Y_i$ as follows: for each $i=1, 2$,
\[Y_i = \sum_{j=1}^6 Y_{ij},\]
where
\begin{align*}
	&Y_{i1} :=\int_{\mathbb{T}^2}  \int_{\mathbb{R}}  \frac{a}{\sigma_i^*}\rho(\widetilde{h}_{i})_{x_1}^{X_i} (p(v)-p(\widetilde{v}))  \, dx_1 \, dx',\\
	&Y_{i2} := - \int_{\mathbb{T}^2}  \int_{\mathbb{R}} a\rho p(\widetilde{v}_i)_{x_1}^{X_i}(v-\widetilde{v})  \, dx_1 \, dx',\\
	&Y_{i3} := \int_{\mathbb{T}^2}  \int_{\mathbb{R}} a\rho(\widetilde{h}_{i})_{x_1}^{X_i} \left( h_1-\widetilde{h}-\frac{p(v)-p(\widetilde{v})}{\sigma_i^*} \right) \, dx_1 \, dx'\\
	&Y_{i4} := - \int_{\mathbb{T}^2}  \int_{\mathbb{R}} a\rho (\widetilde{v}_i)_{x_1}^{X_i} \left( p'(\widetilde{v})-p'(\widetilde{v}_i^{X_i}) \right)\left( v-\widetilde{v} \right)  \, dx_1 \, dx', \\
	&Y_{i5} :=-\frac{1}{2} \int_{\mathbb{T}^2}  \int_{\mathbb{R}} \rho ({a_i})_{x_1}^{X_i} \left( h_1-\widetilde{h} - \frac{p(v)-p(\widetilde{v})}{\sigma_i^*} \right)\left( h_1-\widetilde{h} + \frac{p(v)-p(\widetilde{v})}{\sigma_i^*} \right) \, dx_1 \, dx'\\
	&Y_{i6} :=-\int_{\mathbb{T}^2}  \int_{\mathbb{R}} \rho({a_i})_{x_1}^{X_i} \left(Q(v|\widetilde{v})+\frac{|p(v)-p(\widetilde{v})|^2}{2{\sigma_i^*}^2}+\frac{h_2^2+h_3^2}{2} \right) \, dx_1 \, dx'.
\end{align*}
Notice that it follows from our construction \eqref{X(t)} on shifts $X_i$ that
\[\dot{X}_i = -\frac{M}{\delta_i}(Y_{i1}+Y_{i2}),\]
which implies 
\begin{equation}\label{XiYi}
	\dot{X}_iY_i(U) = -\frac{\delta_i}{M}|\dot{X_i}|^2+\dot{X}_i\sum_{j=3}^6Y_{ij}.
\end{equation}
Here, the good term $-\frac{\delta_i}{M}|\dot{X_i}|^2$ would give an average of linear perturbation on $v$-variable as mentioned above, whereas the remaining part would be controlled by the good terms in $\mathcal{G}(U)$. To show it, we combine \eqref{est-1} and \eqref{XiYi} to get
\begin{align}
	\begin{aligned}\label{est}
	&\frac{d}{dt} \int_{\mathbb{T}^2}\int_\R a \eta (U| \tU) \,dx_1 \,dx' =\mathcal{R}, \quad\mbox{where}\\
	&\mathcal{R} :=-\sum_{i=1}^2\frac{\delta_i}{M}|\dot{X}_i|^2 + \sum_{i=1}^2\left(\dot{X_i}\sum_{j=3}^6Y_{ij}\right) +\sum_{i=1}^8\mathcal{B}_i -\mathcal{G}_1-\mathcal{G}_2-\mathcal{G}_3-\mathcal{D}. \\
	&\quad= \underbrace{-\sum_{i=1}^2\frac{\delta_i}{2M}|\dot{X}_i|^2 + \mathcal{B}_1 +\mathcal{B}_2 -\mathcal{G}_2-\frac{3}{4}\mathcal{D}}_{=:\mathcal{R}_1}\\
	&\quad\quad \underbrace{-\sum_{i=1}^2\frac{\delta_i}{2M}|\dot{X}_i|^2 + \sum_{i=1}^2\left(\dot{X_i}\sum_{j=3}^6Y_{ij}\right) +\sum_{i=3}^8\mathcal{B}_i -\mathcal{G}_1-\frac{1}{4}\mathcal{D}}_{=:\mathcal{R}_2}.
	\end{aligned}
\end{align}

A reason of the decomposition \eqref{est} is that the bad terms contained in $\mathcal{R}_1$ must be estimated delicately. To estimate these terms, we use the sharp Poincar\'e-type inequality of Lemma \ref{lem-poin}. The remaining terms in $\mathcal{R}_2$ can be estimated in a rather rough way. 
First, we concentrate on the estimate of $\mathcal{R}_1$.

\subsection{Estimate of the main part $\mathcal{R}_{1}$}\label{sec:4.5}
An essential idea for estimates of $\mathcal{R}_{1}$ is to apply the Poincar\'e-type inequality \eqref{poincare inequality} in Lemma \ref{lem-poin}. To apply the Poincar\'e-type inequality, we require to localize the perturbation near each wave by using the cutoff functions $\Phi_1, \Phi_2$ defined in \eqref{def cutoff} (see Remark \ref{rem-phi}), and then change of variables from whole space $\bbr$ to a bounded interval $(0,1)$ for each wave.  
For any fixed $t>0$, we will consider the following change of variables in space: 
\[
y_1:=1-\frac{p(v_m)-p(\widetilde{v}_1(x_1-\sigma_1t-X_1(t)))}{\delta_1},\qquad y_2 := \frac{p(v_m)-p(\widetilde{v}_2(x_1-\sigma_2t-X_2(t)))}{\delta_2}.
\]
Surely, for each $i=1, 2$, $y_i :\bbr\to (0,1)$ is a monotone function such that
\[\frac{d y_1}{dx_1} = \frac{1}{\delta_1}p(\widetilde{v}_1)' >0,  \quad \frac{d y_2}{dx_1} = -\frac{1}{\delta_2}p(\widetilde{v}_2)' >0, \]
and
\[\lim_{x_1\to-\infty}y_i = 0,\quad \lim_{x_1\to+\infty} y_i = 1.\]

In addition, both $|X_1(t)|$ and $|X_2(t)|$ are bounded on $\left[0, T\right]$ by \eqref{dxbound}.

With regard to the new variables, we will apply the Poincar\'e-type inequality to each perturbation $w_i$ localized by $\Phi_i$ respectively:
 \[
 \begin{aligned}
&w_1:=\Phi_1(t, x_1) \Bigg(p(v(t,x)) - p\Big(\tv_1(x_1-\sigma_1 t-X_1(t))+\tv_2(x_1-\sigma_2 t-X_2(t))-v_m \Big) \Bigg), \\
&w_2:=\Phi_2(t, x_1) \Bigg(p(v(t,x)) - p\Big(\tv_1(x_1-\sigma_1 t-X_1(t))+\tv_2(x_1-\sigma_2 t-X_2(t))-v_m \Big) \Bigg).
\end{aligned}
\]
In what follows, for simplicity, we use the following notations to denote constants of $O(1)$-scale:
\[\sigma_m:=\sqrt{-p'(v_m)},\quad \alpha_m:=\frac{\gamma+1}{2\gamma\sigma_mp(v_m)}=\frac{p''(v_m)}{2\sigma_m|p'(v_m)|^2},\]
which are surely independent of the shock strengths $\delta_i$ since $v_+/2\le v_m \le v_+$.\\
Because the $\delta_i$ are bounded by $\delta_0$ respectively, the following estimates on the $O(1)$-constants hold:
\begin{equation}\label{shock_speed_est}
	|\sigma_1^*- (-\sigma_m)| \leq C \delta_1,\quad |\sigma_2^*-\sigma_m|\le C\delta_2,
\end{equation}
and
\begin{equation}\label{shock_speed_est-2}
	\|\sigma_m^2-|p'(\tv_i)|\|_{L^\infty} \le C\delta_i,\quad \left\|\frac{1}{\sigma_m^2}-\frac{p(\tv_i)^{-\frac{1}{\gamma}-1}}{\gamma}\right\|_{L^\infty}\le C\delta_i,\quad \left\|\frac{1}{\sigma_m^2}-\frac{p(\tv)^{-\frac{1}{\gamma}-1}}{\gamma}\right\|_{L^\infty}\le C\delta_0.
\end{equation}
We are now ready to estimate the terms in $\mathcal{R}_1$. 
As mentioned, we need to extract a good term on an average of the perturbation $w_i$ from the shift part $\frac{\delta_i}{2M}|\dot{X}_i|^2$  as follows, so that we could apply Lemma \ref{lem-poin}.\\

\noindent $\bullet$  {\bf(Estimate of shift part $\frac{\delta_i}{2M}|\dot{X}_i|^2$):}
Our goal is that: for each $i=1, 2$,
		\begin{align}
		\begin{aligned} \label{L2.2}
		-\frac{\delta_i}{2M}|\dot{X}_i|^2 \le &-\frac{M\delta_i}{\sigma_m^4v_m^2} \left( \int_{\mathbb{T}^2}\int_{0}^{1} w_i \, dy_i \, dx' \right)^2 +  C\delta_i(\delta_0+\delta_i+\nu)^2 \int_{\mathbb{T}^2}  \int_{0}^{1} |w_i|^2 \, dy_i \, dx'\\
        & +C\delta_i^2e^{-C\delta_it}\int_{\mathbb {T}^2} \int_{\mathbb{R}} \eta(U | \widetilde{U}) \, dx_1 \, dx'.
	\end{aligned}
	\end{align}
	As the estimates of $\frac{\delta_1}{2M}|\dot{X}_1|^2$ and $\frac{\delta_2}{2M}|\dot{X}_2|^2$ are the same, we enough to handle the case of $X_1$. Since $\dot{X}_1=-\frac{M}{\delta_1}(Y_{11}+Y_{12})$, we first estimate $Y_{11}$ and $Y_{12}$.\\
 Using (4.3) and the relation $\Phi_1+\Phi_2=1$, we have
	\begin{align*}
		Y_{11} 
		&=\int_{\mathbb{T}^2}  \int_{\mathbb{R}} \frac{a}{{\sigma_1^*}^2v}\Phi_1p(\widetilde{v}_1^{X_1})_{x_1}(p(v)-p(\widetilde{v})) \, dx_1 \, dx'+\int_{\mathbb{T}^2}  \int_{\mathbb{R}} \frac{a}{{\sigma_1^*}^2v}\
\Phi_2p(\widetilde{v}_1^{X_1})_{x_1}(p(v)-p(\widetilde{v})) \, dx_1 \, dx'.
	\end{align*}
Then using \eqref{shock_speed_est} and $\|a-1\|_{L^\infty(\bbr_+\times\bbr)}\le \nu$, we have
	\[\left| Y_{11}-\frac{\delta_1}{\sigma_m^2 v_m} \int_{\mathbb{T}^2}  \int_{0}^{1} w_1 \, dy_1 \, dx' \right|\le C\delta_1(\delta_0+\nu)\int_{\mathbb{T}^2}  \int_{0}^{1} |w_1| \, dy_1 \, dx' + C \int_{\mathbb{T}^2}  \int_{\mathbb{R}} |(\widetilde{v}_1)_{x_1}^{X_1}| |\Phi_2| |p(v)-p(\widetilde{v})| \, dx_1 \, dx' .\]
 When we estimate $Y_{12}$, we first use Taylor expansion in terms of $v=p(v)^{-1/\gamma}$ to get
\[\left|v-\tv -\left(-\frac{p(\widetilde{v})^{-\frac{1}{\gamma}-1}}{\gamma}(p(v)-p(\widetilde{v}))\right)\right|\le C|p(v)-p(\tilde{v})|^2,\]
	which together with the estimates \eqref{shock_speed_est-2} and \eqref{smp1} implies
	\[\left|v-\tv -\left(-\frac{1}{\sigma_m^2}(p(v)-p(\widetilde{v}))\right)\right|\le C(\delta_0+\varepsilon_1)|p(v)-p(\widetilde{v})|.\]	
As above,
\begin{align*}
    \left| Y_{12}-\frac{\delta_1}{\sigma_m^2v_m}\int_{\mathbb{T}^2}  \int_{0}^{1} w_1  \, dy_1 \, dx' \right| \le&  C\delta_1(\nu+\delta_0+\varepsilon_1)\int_{\mathbb{T}^2}  \int_{0}^{1} |w_1|  \, dy_1 \, dx'\\
    &+C\int_{\mathbb{T}^2}  \int_{\mathbb{R}} |(\widetilde{v}_1)_{x_1}^{X_1}|  |\Phi_2| |p(v)-p(\widetilde{v})|  \, dy_1 \, dx'.
\end{align*}
	Combining the estimates for $Y_{11}$ and $Y_{12}$, we have
	\begin{align*}
		\left| \dot{X}_1+\frac{2M}{\sigma_m^2v_m} \int_{\mathbb{T}^2}\int_{0}^{1} w_1 \, dy_1 \, dx' \right|&\le \frac{M}{\delta_1}	\left( \left| Y_{11}-\frac{\delta_1}{\sigma_m^2 v_m} \int_{\mathbb{T}^2}  \int_{0}^{1} w_1 \, dy_1 \, dx' \right| + 	\left| Y_{12}-\frac{\delta_1}{\sigma_m^2v_m}\int_{\mathbb{T}^2}  \int_{0}^{1} w_1  \, dy_1 \, dx' \right| \right)\\
		&\le C(\nu+\delta_0+\varepsilon_1)\int_{\mathbb{T}^2}  \int_{0}^{1} |w_1| \, dy_1 \, dx' \\
  &\quad + \frac{C}{\delta_1}\int_{\mathbb{T}^2}  \int_{\mathbb{R}} |(\widetilde{v}_1)_{x_1}^{X_1}| |\Phi_2| |p(v)-p(\widetilde{v})| \, dx_1 \, dx' ,
	\end{align*}
	which implies by squaring both sides and using Young's inequality, 
 \begin{align}
 \begin{aligned} \label{est-X1}
\frac{2M^2}{\sigma_m^4v_m^2} \left( \int_{\mathbb{T}^2}\int_{0}^{1} w_1 \, dy_1 \, dx' \right)^2 -|\dot{X}_1|^2 \le & \ C(\nu+\delta_0+\varepsilon_1)^2 \int_{\mathbb{T}^2}  \int_{0}^{1} |w_1|^2 \, dy_1 \, dx'\\ 
&+\frac{C}{\delta_1^2}\left( \int_{\mathbb{T}^2}  \int_{\mathbb{R}} |(\widetilde{v}_1)_{x_1}^{X_1}| |\Phi_2| |p(v)-p(\widetilde{v})| \, dx_1 \, dx' \right)^2.
 \end{aligned}
 \end{align}
In addition, using Lemmas \ref{lem:shock-est}, \ref{lem-rel-quant}, and \ref{shock interaction lemma2}, 	
\begin{align*}
		\frac{C}{\delta_1^2}\left( \int_{\mathbb{T}^2}  \int_{\mathbb{R}} |(\widetilde{v}_1)_{x_1}^{X_1}| |\Phi_2| |p(v)-p(\widetilde{v})| \, dx_1 \, dx' \right)^2 
		&\leq \frac{C}{\delta_1^2}\int_{\mathbb{T}^2}\int_\R (|(\widetilde{v}_1)_{x_1}^{X_1} | \Phi_2)^2 dx_1 \, dx' \int_{\mathbb{T}^2}\int_\R | p(v)-p(\widetilde{v}) |^2 \, dx_1 \, dx' \\
		&\leq \frac{C}{\delta_1^2}\sup_{t,x}\big(\Phi_2^2 |(\widetilde{v}_1)_{x_1}^{X_1} |\big)  \|(\widetilde{v}_1)_{x_1}^{X_1}\|_{L^1(\bbr)}  \int_{\mathbb{T}^2} \int_\R Q(v|\tv) \, dx_1 \, dx' \\
		&\leq C \delta_1 \exp(-C \delta_1 t) \int_{\mathbb{T}^2} \int_\R \eta(U|\widetilde{U}) \, dx_1 \, dx'.
	\end{align*}	
In short, combining the estimate above, we have following estimate on $\dot{X}_1$:
	\begin{align*}
		-\frac{\delta_1}{2M}|\dot{X}_1|^2 \le & -\frac{M\delta_i}{\sigma_m^4v_m^2} \left( \int_{\mathbb{T}^2}\int_{0}^{1} w_i \, dy_i \, dx' \right)^2 +  C\delta_i(\delta_0+\delta_i+\nu)^2 \int_{\mathbb{T}^2}  \int_{0}^{1} |w_i|^2 \, dy_i \, dx', \\ 
  &+ C\delta_i^2e^{-C\delta_it}\int_{\mathbb {T}^2} \int_{\mathbb{R}} \eta(U | \widetilde{U}) \, dx_1 \, dx'.	\end{align*}
	which is the desired inequality \eqref{L2.2}.\\
	
\noindent $\bullet$ {\bf (Estimate of the bad term $\mathcal{B}_1$ and good term $\mathcal{G}_2$):}
Recall
\begin{align*}
	&\mathcal{B}_1(U) := \sum_{i=1}^2\underbrace{\frac{1}{2\sigma_i^*}\int_{\mathbb{T}^2}\int_{\bbr}(a_i)_{x_1}^{X_i} |p(v)-p(\widetilde{v})|^2\, dx_1 \,dx'}_{=:\mathcal{B}_{i1}},\\
	&\mathcal{G}_2(U):= \sum_{i=1}^2\underbrace{ \sigma_i^*\int_{\mathbb{T}^2}\int_{\bbr}(a_i)_{x_1}^{X_i}Q(v|\widetilde{v})\,dx_1 \, dx'}_{=:\mathcal{G}_{i2}}.
\end{align*}
Since the estimates for the two cases $i=1,2$ are the same, we enough to deal with the case of $i=2$ for simplicity.\\
First, we use the estimate on $Q(v|\tv)$ in Lemma \ref{lem-rel-quant} to obtain
	\begin{align*}
		\mathcal{G}_{22}&\ge \sigma_2^*\int_{\mathbb{T}^2} \int_{\bbr} (a_2)_{x_1}^{X_2} \frac{p(\widetilde{v}^{X_2}_2)^{-\frac{1}{\gamma}-1}}{2\gamma}|p(v)-p(\widetilde{v})|^2\, dx_1 \,dx'\\
		&\quad -\sigma_2^* \int_{\mathbb{T}^2} \int_{\bbr}(a_2)_{x_1}^{X_2}\frac{1+\gamma}{3\gamma^2}\widetilde{p}^{-\frac{1}{\gamma}-2}(p(v)-p(\widetilde{v}))^3\,dx_1 \,dx' \\
		& \quad + \frac{\sigma_2^*}{2\gamma}\int_{\mathbb{T}^2} \int_{\bbr}(a_2)_{x_1}^{X_2} \left(p(\widetilde{v})^{-\frac{1}{\gamma}-1}-p(\widetilde{v_2}^{X_2})^{-\frac{1}{\gamma}-1}\right)|p(v)-p(\widetilde{v})|^2\,dx_1 \,dx'.
	\end{align*}
For simplicity, let $\widehat{\mathcal{G}}_{22}$ denote the good term given as
	\[\widehat{\mathcal{G}}_{22}:=\sigma_2^* \int_{\mathbb{T}^2} \int_\R (a_2)_{x_1}^{X_2} \frac{p(\widetilde{v}^{X_2}_2)^{-\frac{1}{\gamma}-1}}{2\gamma}|p(v)-p(\widetilde{v})|^2 \, dx_1 \,dx'.
	\]
Using \eqref{shock_speed_est} and \eqref{shock_speed_est-2}, we have
\begin{align*}
\mathcal{B}_{21} \leq \frac{1}{2\sigma_m} \int_{\mathbb{T}^2} \int_\R (a_2)_{x_1}^{X_2} |p(v)-p(\widetilde{v})|^2 \, dx_1 \, dx'+ \frac{C \delta_2}{2\sigma_m} \int_\R  (a_2)_{x_1}^{X_2} |p(v)-p(\widetilde{v})|^2 \, dx_1 \, dx',
	\end{align*}
	and
	\begin{align*}
		\widehat{\mathcal{G}}_{22} \geq \frac{1}{2\sigma_m}(1-C \delta_2) \int_{\mathbb{T}^2} \int_\R (a_2)^{X_2}_{x_1} |\pv-\tpv|^2 \, dx_1 \, dx'
	\end{align*}
Then using $\Phi_1+\Phi_2=1$ and \eqref{avd}, we estimate
	\begin{align}
	\begin{aligned}\label{B21G22}
		\mathcal{B}_{21} - \widehat{\mathcal{G}}_{22}
		&\leq C\delta_2 \int_{\mathbb{T}^2} \int_\R (a_2)_{x_1}^{X_2} |p(v)-p(\widetilde{v})|^2 \, dx_1 \, dx'\\
		&\le C\delta_2\nu_2 \int_{\mathbb{T}^2} \int_{\bbr} \frac{|p(\tv_2^{X_2})_{x_1}|}{\delta_2} |\Phi_2(p(v)-p(\widetilde{v}))|^2\,dx_1 \,dx'+  C\nu_2 \int_{\mathbb{T}^2} \int_{\bbr} |(\tv_2)_{x_1}^{X_2}| \Phi_1^2|p(v)-p(\widetilde{v})|^2\,dx_1 \,dx'.
	\end{aligned}
	\end{align}
The first term of the right-hand side of \eqref{B21G22} is rewritten in the new variables $w_2, y_2$:
\[
C\delta_2\nu_2 \int_{\mathbb{T}^2} \int_{\bbr} \frac{|p(\tv_2^{X_2})_{x_1}|}{\delta_2} |\Phi_2(p(v)-p(\widetilde{v}))|^2\,dx_1 \,dx' = C\delta_2\nu_2 \int_{\mathbb{T}^2} \int_0^1 |w_2|^2\,dy_2 dx'.
\]	
Using Lemma \ref{shock interaction lemma2}, the last term of \eqref{B21G22} can be estimated as
\begin{align*}
C\nu_2 \int_{\mathbb{T}^2} \int_{\bbr} |(\tv_2)^{X_2}_{x_1}| \Phi_1^2|p(v)-p(\widetilde{v})|^2\,dx_1 \, dx' \leq C \nu_2 \delta^2_2 \exp(-C \delta_2 t) \int_{\mathbb{T}^2} \int_\R \eta(U|\widetilde{U}) \, dx_1 \, dx'.
\end{align*}
Hence, we have
	\[\mathcal{B}_1-\widehat{\mathcal{G}}_{2}\le \sum_{i=1}^2 C\nu_i \delta_i\int_{\mathbb{T}^2}\int_0^1|w_i|^2\,dy_i\,dx'+C\nu_i\delta_i^2\exp(-C\delta_it)\int_{\mathbb{T}^2}\int_{\R}\eta(U|\widetilde{U})\,dx_1 \,dx' .\]
\noindent $\bullet$ {\bf (Estimate of the bad term $\mathcal{B}_2$):}
Recall
\begin{align*}
	&\mathcal{B}_2(U) := \sum_{i=1}^2 \underbrace{\sigma_i^*\int_{\mathbb{T}^2}\int_\bbr a(\widetilde{v}_i)_{x_1}^{X_i}p(v|\widetilde{v})\,dx_1 \, dx'}_{=:\mathcal{B}_{i2}}.
\end{align*}
Also, we only deal with the case of $i=2$ for simplicity.
First, we have
\begin{align*}
	\mathcal{B}_{22}= \sigma_2^*\int_{\mathbb{T}^2}\int_{\R}a(\tv_2)_{x_1}^{X_2}\Phi_2^2p(v|\tv)\,dx_1 \,dx'+\sigma_2^*\int_{\mathbb{T}^2}\int_{\R}a(\tv_2)_{x_1}^{X_2}(1-\Phi_2^2)p(v|\tv)\,dx_1 \,dx'.
\end{align*}
Using Lemma \ref{lem-rel-quant}, \eqref{shock_speed_est}, \eqref{shock_speed_est-2}, and then the integration by substitution, we estimate the first term as
\begin{align*}	&\sigma_2^*\int_{\mathbb{T}^2}\int_{\R}a(\tv_2)^{X_2}_x\Phi_2^2p(v|\tv)\,dx_1 \, dx' =\sigma_2^*\int_{\mathbb{T}^2}\int_{\R}a \frac{p(\tv_2^{X_2})_{x_1}}{p'(\tv_2^{X_2})}\Phi_2^2p(v|\tv) \,dx_1 \, dx'\\
	&\le \sigma_2^*\int_{\mathbb{T}^2}\int_{\R}a \frac{|p(\tv_2^{X_2})_{x_1}|}{|p'(\tv_2^{X_2})|}\Phi_2^2 \left(\frac{\gamma+1}{2\gamma p(\widetilde{v})}+C\varepsilon_1\right) | p(v)-p(\tilde{v})|^2 \,dx_1 \,dx' \\
	&\le \frac{\gamma+1}{2\gamma\sigma_mp(v_m)}(1+C(\delta_0+\nu+\varepsilon_1)) \int_{\mathbb{T}^2}\int_\bbr |p(\tv_2^{X_2})_{x_1}|  | \Phi_2(p(v)-p(\widetilde{v}))|^2 \,dx_1 \,dx' \\		&\le\delta_2\alpha_m(1+C(\delta_0+\nu+\varepsilon_1))\int_{\mathbb{T}^2}\int_{0}^1|w_2|^2\,dy_2 \,dx',
\end{align*}
where we used $\alpha_m=\frac{\gamma+1}{2\gamma\sigma_mp(v_m)}$.\\
On the other hand, using Lemma \eqref{shock interaction lemma2} we estimate the last interaction term of $\cB_{22}$ as
\begin{align*}
	\sigma_2^* \int_{\mathbb{T}^2}\int_{\bbr}a(\widetilde{v}_2)_{x_1}^{X_2}(1+\Phi_2)\Phi_1p(v|\widetilde{v})\,dx_1 \, dx' &\leq C \int_{\mathbb{T}^2}\int_\R |(\widetilde{v}_2)_{x_1}^{X_2}| \Phi_1 |\pv-p(\tv)|^2 \, dx_1 \, dx' \\
	&\leq C \delta_2^2 \exp(-C \delta_2 t) \int_{\mathbb{T}^2}\int_\R \eta(U|\widetilde{U}) \, dx_1 \, dx' .
\end{align*}
Thus, we have
\begin{align*}
\mathcal{B}_{22}&\le \delta_2\alpha_m(1+C(\delta_0+\nu+\varepsilon_1))\int_{\mathbb{T}^2}\int_{0}^1|w_2|^2\,dy_2 \, dx'\\
&\qquad +C \delta_2^2 \exp(-C \delta_2 t) \int_{\mathbb{T}^2} \int_\R \eta(U|\widetilde{U}) \, dx_1 \, dx',
\end{align*}
which yields
\begin{align*}
\mathcal{B}_2&\le \sum_{i=1}^2\delta_i\alpha_m(1+C(\delta_0+\nu+\varepsilon_1))\int_{\mathbb{T}^2}\int_0^1|w_i|^2\,dy_i \\
&\qquad +C\delta_i^2 \exp(-C\delta_it)\int_{\mathbb{T}^2}\int_{\R}\eta(U|\widetilde{U})\,dx_1 \,dx'.   
\end{align*}
In short, combining the estimates above on $\mathcal{B}_1,\mathcal{G}_2$, and $\mathcal{B}_2$, we have
\begin{align}
\begin{aligned} \label{L2.3}
\mathcal{B}_1-\mathcal{G}_{2} + \mathcal{B}_2	 &\le
\sum_{i=1}^2\Bigg[\delta_i\alpha_m(1+C(\delta_0+\nu+\varepsilon_1))\int_{\mathbb{T}^2}\int_0^1|w_i|^2\,dy_i \, dx'+C\delta_i^2 \exp(-C\delta_it)\int_{\mathbb{T}^2}\int_{\R}\eta(U|\widetilde{U})\,dx_1 \, dx' \Bigg]\\
&\quad + \sum_{i=1}^2\Bigg[-\sigma_i^* \int_{\mathbb{T}^2}\int_{\R}(a_i)_{x_1}^{X_i}\frac{1+\gamma}{3\gamma^2}p(\widetilde{v})^{-\frac{1}{\gamma}-2}(p(v)-p(\widetilde{v}))^3\,dx_1 \, dx' \\
		&\hspace{2cm} +\frac{\sigma_i^*}{2\gamma}\int_{\mathbb{T}^2}\int_{\R}(a_i)_{x_1}^{X_i}\left(p(\widetilde{v})^{-\frac{1}{\gamma}-1}-p(\tv_i^{X_i})^{-\frac{1}{\gamma}-1}\right)|p(v)-p(\widetilde{v})|^2\,dx_1 \,dx' \Bigg].
		\end{aligned}
	\end{align}
\noindent $\bullet$ {\bf (Estimate of the diffusion term $\mathcal{D}(U)$):}
First of all, using the fact that $\Phi_{1}+\Phi_{2}=1$ and $1\ge \Phi_{i}\ge \Phi_{i}^2\ge0$ for each $i=1, 2$, we separate $\mathcal{D}(U)$ into
\begin{align*}    
\mathcal{D}(U) &= (2\mu+\lambda) \int_{\mathbb{T}^2}  \int_{\mathbb{R}} a\frac{|\partial_{x_1}(p(v)-p(\tilde{v}))|^2}{\gamma p(v)^{1+\frac{1}{\gamma}}}\, dx_1 \, dx'+(2\mu+\lambda) \int_{\mathbb{T}^2}  \int_{\mathbb{R}} a\frac{|\nabla_{x'}(p(v)-p(\widetilde{v}))|^2}{\gamma p(v)^{1+\frac{1}{\gamma}}}\, dx_1 \, dx' \\ 
&\ge (2\mu+\lambda)\sum\limits_{i=1}^2 \int_{\mathbb{T}^2}  \int_{\mathbb{R}} a\Phi_i^2\frac{|\partial_{x_1}(p(v)-p(\tilde{v}))|^2}{\gamma p(v)^{1+\frac{1}{\gamma}}}\, dx_1 \, dx'+(2\mu+\lambda)\sum\limits_{i=1}^2 \int_{\mathbb{T}^2}  \int_{\mathbb{R}} a\Phi_i^2\frac{|\nabla_{x'}(p(v)-p(\widetilde{v}))|^2}{\gamma p(v)^{1+\frac{1}{\gamma}}}\, dx_1 \, dx'
\end{align*}
Since Young's inequality yields: for any $1>\delta_*>0$ small enough, (to be determined below)
\begin{align*}
\frac{2\mu+\lambda}{\gamma}\int_{\mathbb{T}^2}  \int_{\mathbb{R}} a\frac{|\partial_{x_1}(\Phi_i(p(v)-p(\widetilde{v})))|^2}{ p(v)^{1+\frac{1}{\gamma}}}\, dx_1 \, dx' \le& (1+\delta_{\ast})\frac{2\mu+\lambda}{\gamma} \int_{\mathbb{T}^2}  \int_{\mathbb{R}} a\frac{\Phi_i^2|\partial_{x_1}(p(v)-p(\widetilde{v}))|^2}{ p(v)^{1+\frac{1}{\gamma}}}\, dx_1 \, dx',\\
&+\frac{C}{\delta_{\ast}} \frac{2\mu+\lambda}{\gamma}\sum\limits_{i=1}^2 \int_{\mathbb{T}^2} \int_{\mathbb{R}} a\frac{|\partial_{x_1}\Phi_i|^2|p(v)-p(\widetilde{v})|^2}{p(v)^{1+\frac{1}{\gamma}}} \, dx_1 \, dx'
\end{align*}
we have
\begin{align*}
	-\mathcal{D}_1(U) \le &  -\frac{1}{1+\delta_{\ast}}\frac{2\mu+\lambda}{\gamma}\sum\limits_{i=1}^2 \int_{\mathbb{T}^2} \int_{\mathbb{R}} a\frac{|\partial_{x_1}(\Phi_i(p(v)-p(\widetilde{v})))|^2}{p(v)^{1+\frac{1}{\gamma}}} \, dx_1 \, dx' \\
 &+ \frac{C}{\delta_{\ast}} \frac{2\mu+\lambda}{\gamma}\sum\limits_{i=1}^2 \int_{\mathbb{T}^2} \int_{\mathbb{R}} a\frac{|\partial_{x_1}\Phi_i|^2|p(v)-p(\widetilde{v})|^2}{p(v)^{1+\frac{1}{\gamma}}} \, dx_1 \, dx'\\
	=:&  J_1+J_2.
\end{align*}
We want to write $J_1$ in terms of the variables $y_i$ and $w_i$.
So, we apply the following estimates in the proof of \cite[Lemma 4.5]{KVW23} to $J_1$:
\begin{align} \label{diffusion lemma}
\left| \frac{1}{y_i(1-y_i)} \frac{2\mu+\lambda}{\gamma p^{\frac{1}{\gamma}+1}(\widetilde{v}_i^{X_i})}  \frac{dy_i}{dx_1}-\frac{\delta_i p''(v_m)}{2|p'(v_m)|^2\sigma_m}\right|\le C\delta_i^2
\end{align}
From $\left\|a\frac{p(\tv_i)}{p(v)}-1\right\|_{L^\infty}\le C(\delta_0+\varepsilon_1+\nu)$, \eqref{diffusion lemma} implies
\begin{align*}
	J_1&=-\frac{1}{1+\delta_{\ast}}\frac{2\mu+\lambda}{\gamma}\sum\limits_{i=1}^2 \int_{\mathbb{T}^2} \int_{0}^{1} a\frac{|\partial_{y_i}w_i|^2}{p^{1+\frac{1}{\gamma}}} \frac{dy_i}{dx_1} \, dy_i \, dx'\\
	&\le -\sum\limits_{i=1}^2(1-C(\delta_0+\varepsilon_1+\nu+\delta_{\ast}))\left( \frac{\delta_i p''(v_m)}{2p'(v_m)^2\sigma_m}-C\delta_i^2 \right) \int_{\mathbb{T}^2} \int_{0}^{1} 
y_i(1-y_i)|\partial_{y_i}w_i|^2 \, dy_i\, dx'\\
	&\le -\sum\limits_{i=1}^2\delta_i \alpha_m(1-C(\delta_0+\varepsilon_1+\nu+\delta_{\ast})) \int_{\mathbb{T}^2} \int_{0}^{1} 
y_i(1-y_i)|\partial_{y_i}w_i|^2 \, dy_i\, dx',
\end{align*}
where we used $\frac{p''(v_m)}{2|p'(v_m)|^2\sigma_m} = \frac{\gamma+1}{2\gamma\sigma_mp(v_m)}=\alpha_m$. \\
To deal with the term $J_2$, note that from the definition of cutoff functions \eqref{def cutoff} for each $i=1, 2$, 
\beq\label{phidecay}
|\partial_{x_1} \Phi_{i}(t,x_1)| \le \frac{2}{(\sigma_2- \sigma_1)} \frac{1}{t},\quad\forall x_1\in\bbr, \quad t\in (0, T].
\eeq
Using this, we can estimate as
\[
J_2 \le \frac{C}{\delta_*t^2} \int_{\mathbb{T}^2}\int_\R \eta(U|\widetilde{U}) \, dx_1 \, dx'.
\]
Thus, combining the estimates of $J_1$ and $J_2$ above, we have
		\begin{align}
	\begin{aligned}\label{L2.4}
		-\mathcal{D}_1(U)\le -\sum_{i=1}^2\delta_i \alpha_m(1-C(\delta_0+\nu+\varepsilon_1+\delta_*))\int_{\mathbb{T}^2}\int_0^1y_i(1-y_i)|\pa_{y_i}w_i|^2\,dy_i\,dx' + \frac{C}{\delta_* t^2}\int_{\mathbb{T}^2}\int_{\bbr}\eta(U|\tU)\,dx_1 \,dx'.
		\end{aligned}
	\end{align}
\noindent Since both $\Phi_1$ and $\Phi_2$ are independent of the variables $x'$ and $1=\Phi_1+\Phi_2\ge \Phi_1^2+\Phi_2^2$, we have
\begin{align*}    
-\mathcal{D}_2&=-(2\mu+\lambda) \int_{\mathbb{T}^2}  \int_{\mathbb{R}} a\frac{|\nabla_{x'}(p(v)-p(\widetilde{v}))|^2}{\gamma p(v)^{1+\frac{1}{\gamma}}}\, dx_1 \, dx'\\
&\le -(2\mu+\lambda)\sum\limits_{i=1}^2 \int_{\mathbb{T}^2}\int_{\mathbb{R}} a\frac{|\nabla_{x'}(\Phi_i(p(v)-p(\widetilde{v})))|^2}{\gamma p(v)^{1+\frac{1}{\gamma}}}  \, dx_1\, dx'\\
&= -(2\mu+\lambda)\sum\limits_{i=1}^2 \int_{\mathbb{T}^2}\int_{0}^{1} a\frac{|\nabla_{x'}w_i|^2}{\gamma p(v)^{1+\frac{1}{\gamma}}} \left(\frac{dx_1}{dy_i}\right) \, dy_i \, dx'.
\end{align*}
From \eqref{diffusion lemma}, 
\[y_i(1-y_i)\frac{dx_1}{dy_i}\ge \frac{2\mu+\lambda}{\gamma p(\widetilde{v}_i^{X_i})^{1+\frac{1}{\gamma}}(\alpha_m\delta_i+C\delta_i^2)}\ge \frac{2\mu+\lambda}{2\alpha_m\delta_i|p'(v_m)|}.\]
This implies 
\[\int_{0}^{1} a\frac{|\nabla_{x'}w_i|^2}{\gamma p(v)^{1+\frac{1}{\gamma}}} \left(\frac{dx_1}{dy_i}\right) \, dy_i \ge (1-C(\delta_i+\nu))\frac{\sigma_m(2\mu+\lambda)^2}{\delta_i p''(v_m)}\int_{0}^{1}\frac{|\nabla_{x'}w_i|^2}{y_i(1-y_i)}dy_i.  \]
Thus, we get 
\begin{align}
\begin{aligned}\label{L2.4'}
-\mathcal{D}_2(U)\le  -\sum_{i=1}^2(1-C(\delta_i+\nu))\frac{\sigma_m(2\mu+\lambda)^2}{\delta_ip''(v_m)} \int_{\mathbb{T}^2} \int_{0}^1 \frac{|\nabla_{x'}w_i|^2}{y_i(1-y_i)} \, dy_i \, dx'.     
\end{aligned}    
\end{align}
\noindent $\bullet$ {\bf (Conclusion):}
Combining the estimates \eqref{L2.3}, \eqref{L2.4}, and \eqref{L2.4'} we have
\begin{align}
\begin{aligned} \label{4.36}
	\mathcal{B}_1&+\mathcal{B}_2-\widehat{\mathcal{G}_{2}}-\frac{3}{4}\mathcal{D}\\
	&\le \sum_{i=1}^2\delta_i\alpha_m\Bigg(1+C(\delta_0+\nu+\varepsilon_1))\int_{\mathbb{T}^2}\int_0^1|w_i|^2\,dy_i\,dx'\\
 &\qquad \qquad \qquad -\frac{3}{4}(1-C_0(\delta_0+\nu+\varepsilon_1+\delta_*))\int_{\mathbb{T}^2}\int_0^1y_i(1-y_i)|\pa_{y_i}w_i|^2\,dy_i\,dx'\Bigg)\\
	&\ \quad +C\left(\sum_{i=1}^2 \delta_i^2 \exp(-C \delta_i t) + \frac{1}{\delta_*t^2}  \right) \int_\R \eta (U|\widetilde{U})\,dx\\
        &\ \quad -\frac{3}{4}\sum_{i=1}^2(1-C(\delta_i+\nu))\frac{\sigma_m(2\mu+\lambda)^2}{\delta_ip''(v_m)} \int_{\mathbb{T}^2} \int_{0}^1 \frac{|\nabla_{x'}w_i|^2}{y_i(1-y_i)} \, dy_i \, dx'.
\end{aligned}
\end{align}
At this time, we take $\delta_*$ as
\[
\delta_* = \frac{1}{24C_0},
\]
which together with the smallness of $\delta_0,\nu,\varepsilon_1$ yields
\[
C_0(\delta_0+\nu+\varepsilon_1+\delta_*) <\frac{1}{12}.
\]
Substituting this into \eqref{4.36}
\begin{align*}
	\mathcal{B}_1&+\mathcal{B}_2-\widehat{\mathcal{G}_{2}}-\frac{3}{4}\mathcal{D}\\
	&\le  \sum_{i=1}^2\delta_i\alpha_m\left(\frac{11}{10}\int_{\mathbb{T}^2}\int_0^1|w_i|^2\,dy_i\,dx'-\frac{2}{3}\int_{\mathbb{T}^2}\int_0^1y_i(1-y_i)|\pa_{y_i}w_i|^2\,dy_i\,dx'\right) \\
	&\quad +C\left(\sum_{i=1}^2 \delta_i \exp(-C \delta_i t) + \frac{1}{ t^2}  \right) \int_{\mathbb{T}^2}\int_\R \eta (U|\widetilde{U})\,dx_1\,dx'\\
 &\quad -\frac{5}{8}(2\mu+\lambda)^2\sum\limits_{i=1}^2 \frac{\sigma_m}{\delta_i p''(v_m)}\int_{\mathbb{T}^2}  \int_{0}^{1} \frac{|\nabla_{x'}w_i|^2}{y_i(1-y_i)} \, dy_1 \, dx'.
\end{align*}
Note that the identity:
\[
\int_0^1 |w_i-\bar {w_i}|^2 dy_i = \int_0^1w_i^2 dy_i -{\bar {w_i}}^2,\qquad \bar {w_i}:=\int_0^1 w_i dy_i, \quad \text{for} \ i=1, 2.
\]
Using Lemma \ref{lem-poin} with identity above, we have
\begin{align*}
	\mathcal{B}_1+\mathcal{B}_2-\widehat{\mathcal{G}_{2}}-\frac{3}{4}\mathcal{D}
	&\le \sum_{i=1}^2\left[-\frac{7\delta_i\alpha_m}{30}\int_{\mathbb{T}^2}\int_0^1|w_i|^2\,dy_i\,dx' +\frac{4\delta_i\alpha_m}{3}\left(\int_{\mathbb{T}^2}\int_0^1w_i\,dy_i\,dx'\right)^2\right]\\
	&\quad +C\left(\sum_{i=1}^2 \delta_i^2 \exp(-C \delta_i t) + \frac{1}{t^2}  \right) \int_{\mathbb{T}^2}\int_\R \eta (U|\widetilde{U})\,dx_1 \,dx' \\
        &\quad -\frac{5}{8}\sum_{i=1}^2\left(\frac{(2\mu+\lambda)^2\sigma_m}{\delta_ip''(v_m)}-\frac{2\delta_i\alpha_m}{15\pi}\right)\int_{\mathbb{T}^2}\int_{0}^{1}\frac{|\nabla_{x'}w_i|^2}{y_i(1-y_i)}dy_idx'.
\end{align*}
Lastly, using \eqref{L2.2} with the choice $M=\frac{4}{3}\sigma_m^4v_m^2\alpha_m$ and smallness of $\delta_i$, we have
\begin{align*}
	\begin{aligned}
	&-\sum_{i=1}^2\frac{\delta_i}{2M}|\dot{X}_i|^2+\mathcal{B}_1+\mathcal{B}_2-\mathcal{G}_2-\frac{3}{4}\mathcal{D}\\
	&\le \sum_{i=1}^2\Bigg[-\frac{\delta_i\alpha_m}{5}\int_{\mathbb{T}^2}\int_0^1|w_i|^2\,dy_i\,dx'+\sigma_i^*\int_{\mathbb{T}^2}\int_{\bbr}(a_i)_{x_1}^{X_i}\frac{1+\gamma}{3\gamma^2}p(\widetilde{v})^{-\frac{1}{\gamma}-2}(p(v)-p(\widetilde{v}))^3\,dx_1\,dx'\\
	&\hspace{2cm} -\frac{\sigma_i^*}{2\gamma}\int_{\mathbb{T}^2}\int_{\bbr}(a_i)_{x_1}^{X_i} \left(p(\widetilde{v})^{-\frac{1}{\gamma}-1}-p(\widetilde{v}^{X_i}_i)^{-\frac{1}{\gamma}-1}\right)|p(v)-p(\widetilde{v})|^2\,dx\Bigg]\\
	&\quad+C\left(\sum_{i=1}^2 \delta_i^2 \exp(-C \delta_i t) + \frac{1}{t^2}  \right) \int_\R \eta (U|\widetilde{U})\,dx,
	\end{aligned}
\end{align*}
which concludes
\begin{align}
	\begin{aligned}\label{est-I1}
	\mathcal{R}_1
	&\le \sum_{i=1}^2\Bigg[-C_1 \int_{\mathbb{T}^2}\int_\R |(\widetilde{v}_i)_{x_1}^{X_i}| |\Phi_i(p(v)-p(\widetilde{v})) |^2 \, dx_1 \,dx'
	+C \int_{\mathbb{T}^2}\int_\R |(a_i)_{x_1}^{X_i}| |p(v)-p(\widetilde{v})|^3 \, dx_1 \,dx' \\
	&\hspace{1.5cm}+C \int_{\mathbb{T}^2}\int_{\bbr} |(a_i)_{x_1}^{X_i}| |\tv-\tvi^{X_i}| |p(v)-p(\widetilde{v})|^2\,dx_1 \,dx' \Bigg]\\
	&\quad+C\left(\sum_{i=1}^2 \delta_i^2\exp(-C \delta_i t) + \frac{1}{t^2}  \right) \int_{\mathbb{T}^2}\int_\R \eta (U|\widetilde{U})\,dx_1\,dx'.
	\end{aligned}
\end{align}
\subsection{Estimate of the remaining part $\mathcal{R}_2$}\label{sec:4.6}
Substituting \eqref{est-I1} into \eqref{est} and using Young's inequality
\beq\label{2young}
\sum_{i=1}^2\left(\dot{X}_i\sum_{j=3}^6Y_{ij}\right)\le \sum_{i=1}^2\frac{\delta_i}{4M}|\dot{X}_i|^2 +\sum_{i=1}^2\frac{C}{\delta_i}\sum_{j=3}^6|Y_{ij}|^2,
\eeq
we have
\begin{align}
	\begin{aligned}\label{est-2}
	&\frac{d}{dt}\int_{\mathbb{T}^2}\int_{\bbr}a\rho\eta(U|\widetilde{U})\,dx_1\,dx' \\
	&\leq-C_1\mathcal{G}^S + \mathcal{K}_1+\mathcal{K}_2 \\
	&\quad +C\left(\sum_{i=1}^2 \delta_i^2 \exp(-C \delta_i t) + \frac{1}{t^2}  \right) \int_{\mathbb{T}^2}\int_\R \eta (U|\widetilde{U})\,dx_1 \,dx' \\
	&\quad  -\sum_{i=1}^2\frac{\delta_i}{4M}|\dot{X}_i|^2 + \sum_{i=1}^2 \frac{C}{\delta_i}\sum_{j=3}^6|Y_{ij}|^2 +\sum_{i=3}^8\mathcal{B}_i -\mathcal{G}_1-\mathcal{G}_3-\frac{1}{4}\mathcal{D},
	\end{aligned}
\end{align}
where
\begin{align*}
&\mathcal{G}^S:=\sum_{i=1}^2 \int_{\mathbb{T}^2}\int_\R |(\widetilde{v}_i)^{X_i}_{x_1}| |\Phi_i(p(v)-p(\widetilde{v})) |^2 \,dx_1\,dx' \\
&\mathcal{K}_1:=\sum_{i=1}^2 \underbrace{C \int_{\mathbb{T}^2}\int_\R |(a_i)_{x_1}^{X_i}| |p(v)-p(\widetilde{v})|^3 \,dx_1\,dx'}_{=: \mathcal{K}_{i1}} \\
 &\mathcal{K}_2:=\sum_{i=1}^2 \underbrace{C \int_{\mathbb{T}^2}\int_\R |(a_i)_{x_1}^{X_i}| |\tv-\tvi^{X_i}| |p(v)-p(\widetilde{v})|^2\,dx_1 \,dx'}_{=: \mathcal{K}_{i2}}. 
\end{align*}
In what follows,  to control the remaining terms,  we will use the good terms $\cG_1$, $\cG^S$ and the diffusion term $\cD$.\\
\noindent $\bullet$ {\bf (Estimate of $\mathcal{K}_1$):}
 To control the term $\mathcal{K}_{i1}$ by the good terms above, we localize it via $\Phi_i$ as
\begin{equation}\label{Ki1}
\mathcal{K}_{i1}
\leq C \int_{\mathbb{T}^2}  \int_{\mathbb{R}} \frac{\nu_i}{\delta_i} |(\widetilde{v}_i)_{x_1}^{X_1}| \Phi_i |p(v)-p(\widetilde{v})|^3	\, dx_1 \, dx' + C \int_{\mathbb{T}^2}  \int_{\mathbb{R}} \frac{\nu_i}{\delta_i} |(\widetilde{v}_i)_{x_1}^{X_1}| (1-\Phi_i) |p(v)-p(\widetilde{v})|^3	\, dx_1 \, dx'.
\end{equation}
Let $w=p(v)-p(\widetilde{v})$. Using the Lemma \ref{infty interpolation inequality lemma} and $\delta_i\ll\nu_i\le C\sqrt{\delta_i}$, the first term is controlled by the good terms as
\begin{align*}
&\int_{\mathbb{T}^2}  \int_{\mathbb{R}} \frac{\nu_i}{\delta_i} |(\widetilde{v}_i)_{x_1}^{X_i}| \Phi_i |p(v)-p(\widetilde{v})|^3	\, dx_1 \, dx'\\
&\le C\frac{\nu_i}{\delta_i}\lVert w \rVert_{L^\infty}\int_{\mathbb{T}^2}\int_{\mathbb{R}}|(\widetilde{v}_i)_{x_1}^{X_i}| |w_i|  dx_1 dx' \\
&\leq C\frac{\nu_i}{\delta_i} \left(\lVert w \rVert_{L^2}\lVert \partial_{x_1}w \rVert_{L^2}+\lVert \nabla_xw \rVert_{L^2}\lVert \nabla_x^2w \rVert_{L^2}\right) \sqrt{\int_{\mathbb{T}^2}  \int_{\mathbb{R}} |(\widetilde{v}_i)_{x_1}^{X_i}|  |w_i|^2	\, dx_1 \, dx'} \sqrt{\int_{\mathbb{T}^2}  \int_{\mathbb{R}} |(\widetilde{v}_i)_{x_1}^{X_i}|	\, dx_1 \, dx'}  \\
&\le C\frac{\nu_i}{\sqrt{\delta_i}} \left(\lVert w \rVert_{L^2}+\lVert \nabla_x^2w \rVert_{L^2}\right) \lVert \nabla_xw \rVert_{L^2} \sqrt{\int_{\mathbb{T}^2}  \int_{\mathbb{R}}  |(\widetilde{v}_i)_{x_1}^{X_i}|  |w_i|^2	\, dx_1 \, dx'}\\
&\le C\varepsilon_1 \lVert \nabla_xw \rVert_{L^2} \sqrt{\mathcal{G}^S} \\
&\le C\varepsilon_1 (\mathcal{D}+C_1 \mathcal{G}^S ).
\end{align*}
Using Lemma \ref{shock interaction lemma2}, the second term in \eqref{Ki1} is estimated as
\begin{align*}
C \sum\limits_{i=1}^2 \int_{\mathbb{T}^2}  \int_{\mathbb{R}} \frac{\nu_i}{\delta_i} |(\widetilde{v}_i)_{x_1}^{X_i}| (1-\Phi_i) |p(v)-p(\widetilde{v})|^3	\, dx_1 \, dx'  &\leq C \sum\limits_{i=1}^2 \varepsilon_1 \delta_i \nu_i e^{-C\delta_it} \int_{\mathbb{T}^2}  \int_{\mathbb{R}} |p(v)-p(\widetilde{v})|^2 \, dx_1 \, dx'\\
&\le C \sum\limits_{i=1}^2 \varepsilon_1 \delta_i \nu_i e^{-C\delta_it} \int_{\mathbb{T}^2}  \int_{\mathbb{R}} \eta(U|\widetilde{U}) \, dx_1 \, dx'.
\end{align*}
Combining the estimates above for each $i=1,2$, we get
\begin{align*}
\mathcal{K}_1\leq  \varepsilon_1(\mathcal{D}+C_1 \mathcal{G}^S ) +  C \sum\limits_{i=1}^2 \varepsilon_1 \delta_i \nu_i e^{-C\delta_it} \int_{\mathbb{T}^2}  \int_{\mathbb{R}} \eta(U|\widetilde{U}) \, dx_1 \, dx'.
\end{align*}
\noindent $\bullet$ {\bf (Estimate of $\mathcal{K}_2$):}
To estimate $\mathcal{K}_{2}$, we use Lemma \ref{shock interaction lemma1} for each $i=1, 2$ to get
\begin{align}
\begin{aligned} \label{estimateK2}
\mathcal{K}_{i2}
&=\frac{\nu_i}{\delta_i}  \int_{\mathbb{T}^2}  \int_{\mathbb{R}} |(\widetilde{v}_i)_{x_1}^{X_i}| |\widetilde{v}-\widetilde{v}_i^{X_i}| |p(v)-p(\widetilde{v})|^2\, dx_1 \, dx'\\
&\le C\delta_1\delta_2 e^{-C\min\left\{ \delta_1, \delta_2 \right\}t} \sum\limits_{i=1}^2 \nu_i \int_{\mathbb{T}^2}  \int_{\mathbb{R}} \eta(U|\widetilde{U}) \, dx_1 \, dx'.
\end{aligned}
\end{align}
which gives the estimate on $\mathcal{K}_2$:
\[\mathcal{K}_2\le C\nu\delta_1\delta_2\exp(-C\min\{\delta_1,\delta_2\}t)\int_{\mathbb{T}^2}\int_{\R}\eta(U|\tU)\,dx_1 \,dx'.\]
\noindent $\bullet$ {\bf (Estimate of Special Term):}
Especially, as in $C\mathcal{K}_1$, we get also
\begin{align} \label{specialterm}
\sum\limits_{i=1}^2 \int_{\mathbb{T}^2}  \int_{\mathbb{R}} |(\widetilde{v}_i)_{x_1}^{X_i}| |(p(v)-p(\widetilde{v}))|^2	\, dx_1 \, dx'\le \mathcal{G}^S+C\sum\limits_{i=1}^2 \delta_i^2e^{-C\delta_it} \int_{\mathbb{T}^2}  \int_{\mathbb{R}} \eta(U|\widetilde{U}) \, dx_1 \, dx'.    
\end{align}
This estimate is used in later calculations. 
Since we get this similar way of calculation of $\mathcal{K}_1$, we write this here in advance.

\noindent $\bullet$ {\bf (Estimate of $\frac{C}{\delta_i} |Y_{ij}|^2$ for $i=1,2, j=3,\ldots, 6$):}
Using \eqref{viscous-shock-h}$_2$ and \eqref{inta}, we estimate $Y_{i3}$ as
\begin{align*}
|Y_{i3}|&\le C \int_{\mathbb{T}^2} \int_{\mathbb{R}} \frac{\delta_i}{\nu_i} (a_i)_{x_1}^{X_i} \left| h_1-\widetilde{h}-\frac{p(v)-p(\widetilde{v})}{\sigma_i^*} \right| \, dx_1\, dx' \\
&\le C\frac{\delta_i}{\sqrt{\nu_i}}\sqrt{\int_{\mathbb{T}^2} \int_{\mathbb{R}} \frac{\sigma_i^*}{2} (a_i)_{x_1}^{X_i} \left| h_1-\widetilde{h}-\frac{p(v)-p(\widetilde{v})}{\sigma_i^*} \right|^2 \, dx_1\, dx'}\leq C \frac{\delta_i}{\sqrt{\nu_i}} \sqrt{\cG_1},
\end{align*}
which gives
\[
\frac{C}{\delta_i} |Y_{i3}|^2 \le\frac{C\delta_i}{\nu_i}\mathcal{G}_{1}.
\]
Using Lemma \ref{shock interaction lemma1}, we control $Y_{i4}$ as
\begin{align*}
|Y_{i4}| 
&\le C \int_{\mathbb{T}^2} \int_{\mathbb{R}} 
\left| \widetilde{v}-\widetilde{v}_i^{X_i}\right|  \left| p(v)-p(\widetilde{v}) \right|  \left| (\widetilde{v}_i)_{x_1}^{X_i} \right| \, dx_1\, dx'\\
&C \sqrt{\int_{\mathbb{T}^2} \int_{\mathbb{R}} 
\left| \widetilde{v}-\widetilde{v}_i^{X_i}\right|^2  \left| (\widetilde{v}_i)_{x_1}^{X_i} \right|^2 \, dx_1\, dx'} \sqrt{\int_{\mathbb{T}^2} \int_{\mathbb{R}} 
\eta(U|\widetilde{U})\, dx_1\, dx'}  \\
&\le  C\sqrt{\delta_i} \delta_1 \delta_2 \exp(-C \min \left\{\delta_1, \delta_2\right\}) t)  \sqrt{\int_{\mathbb{T}^2}\int_\R \eta(U |\widetilde{U}) \, dx_1 \, dx'},
\end{align*}
which gives
\[
\frac{C\left| Y_{i,4} \right|^2}{\delta_i}\le C\delta_1^2\delta_2^2 \exp(-C \min\left\{ \delta_1, \delta_2 \right\}t)\int_{\mathbb{T}^2} \int_{\mathbb{R}} 
\eta(U|\widetilde{U})\, dx_1\, dx'.
\]
For $Y_{i5}$,  we first estimate $\bold{h}-\widetilde{\bold{h}}$ in terms of $\bold{u}-\widetilde{\bold{u}}$ and $v-\tv$ as follows.  Using $(\tv)_{x_1}=(\tv_1)_{x_1}^{X_1}+(\tv_2)_{x_1}^{X_2}$ and $C^{-1}\le v, \tv^{X_1}_1, \tv^{X_2}_2 \le C$,  we have
\begin{align}
\begin{aligned}\label{h-est}
||\bold{h}-\widetilde{\bold{h}}||_{L^2} &\leq C||\bold{u}-\widetilde{\bold{u}}||_{L^2} + C||\nabla_x(v-\widetilde{v})||_{L^2}\le C\varepsilon_1,\\
&||p(v)-p(\widetilde{v})||_{L^2}\le C||v-\widetilde{v}||_{L^2}\le \varepsilon_1.\\
\end{aligned}
\end{align}
Hence, we estimate $Y_{i5}$ as  
\begin{align*}
|Y_{i5}|
&\le C \sqrt{\int_{\mathbb{T}^2} \int_{\mathbb{R}} 
\frac{\sigma_i^*}{2} (a_i)_{x_1}^{X_i} \left| h_1-\widetilde{h}-\frac{p(v)-p(\widetilde{v})}{\sigma_i^*}\right|^2 \, dx_1\, dx'} \cdot ||(a_i)_{x_1}^{X_i}||_{L^\infty}^\frac{1}{2} \cdot C\varepsilon_1 \\
&\le C\varepsilon_1(\nu_i \delta_i)^{\frac{1}{2}}\sqrt{\int_{\mathbb{T}^2} \int_{\mathbb{R}} 
\frac{\sigma_i^*}{2} (a_i)_{x_1}^{X_i} \left| h_1-\widetilde{h}-\frac{p(v)-p(\widetilde{v})}{\sigma_i^*}\right|^2 \, dx_1\, dx'}\\
&\leq C\nu_i\delta_i (\varepsilon_1+C\delta_1\delta_2) \sqrt{\mathcal{G}_{1}},
\end{align*}
which gives
\[
\frac{C}{\delta_i} |Y_{i5}|^2 \le C\nu_i (\varepsilon_1+C\delta_1\delta_2)^2 \mathcal{G}_{1}.
\]
Lastly, $\nu_i\le C\sqrt{\delta_i}$, by Lemma \ref{lem-rel-quant} and Lemma \ref{shock interaction lemma1}, we estimate $Y_{i6}$ as
\begin{align*}
\frac{C}{\delta_i} |Y_{i6}|^2 &\le \frac{C}{\delta_i}\Bigg(\int_{\mathbb{T}^2}  \int_{\mathbb{R}} |p(v)-p(\widetilde{v})|^2 |(a_i)_{x_1}^{X_i}| \, dx_1 \, dx'\Bigg)^2+\frac{C}{\delta_i}\Bigg(\int_{\mathbb{T}^2}  \int_{\mathbb{R}} \frac{h_2^2+h_3^2}{2}|(a_i)_{x_1}^{X_i}| \, dx_1 \, dx'\Bigg)^2 \\
&\le \frac{C}{\nu_i\delta_i}\Bigg(\int_{\mathbb{T}^2}  \int_{\mathbb{R}} |p(v)-p(\widetilde{v})|^2 |(a_i)_{x_1}^{X_i}| \, dx_1 \, dx'\Bigg)\Bigg(\int_{\mathbb{T}^2}  \int_{\mathbb{R}} |p(v)-p(\widetilde{v})|^2 |(\widetilde{v}_i)_{x_1}^{X_i}| \, dx_1 \, dx'\Bigg)+C\frac{\mathcal{G}_{i3}^2}{\delta_i}\\
&\le \frac{C \varepsilon_1^2 ||(a_i)_{x_1}^{X_i}||_{L^{\infty}}}{\nu_i\delta_i}\Bigg(\cG^S_i+C\sum\limits_{i=1}^2 \delta_i^2e^{-C\delta_it} \int_{\mathbb{T}^2}  \int_{\mathbb{R}} \eta(U|\widetilde{U}) \, dx_1 \, dx'\Bigg)+\frac{C\mathcal{G}_3}{\delta_i} ||\bold{h}-\widetilde{\bold{h}}||_{L^2}^2 ||a_i^{X_i}||_{L^{\infty}}\\
&\le C \varepsilon_1^2  \Bigg(\mathcal{G}^S+\mathcal{G}_3+\sum\limits_{i=1}^2 \delta_i^2e^{-C\delta_it} \int_{\mathbb{T}^2}  \int_{\mathbb{R}} \eta(U|\widetilde{U}) \, dx_1 \, dx'\Bigg).
\end{align*}
Combining all the estimates for $Y_{ij}$ and using the smallness of the parameters, we conclude that
\begin{align*}
	\sum\limits_{i=1}^2	\frac{C}{\delta_i} \sum\limits_{j=3}^6 |Y_{i,j}|^2 \le& \frac{C_1\mathcal{G}^S+\mathcal{G}_1+\mathcal{G}_3}{10}+C\delta_1^2\delta_2^2 e^{-C \min\left\{ \delta_1, \delta_2 \right\}t}\int_{\mathbb{T}^2} \int_{\mathbb{R}} 
\eta(U|\widetilde{U})\, dx_1\, dx'\\ 
&+C\sum\limits_{i=1}^2 \varepsilon_1^2\delta_i^2e^{-C\delta_it} \int_{\mathbb{T}^2}  \int_{\mathbb{R}} \eta(U|\widetilde{U}) \, dx_1 \, dx'.
\end{align*}
\noindent $\bullet$ {\bf (Estimate of $\mathcal{B}_3$): }
To compute $\mathcal{B}_3=\sum\limits_{i=1}^2  \int_{\mathbb{T}^2}  \int_{\mathbb{R}} ap'(\widetilde{v})(v-\widetilde{v})F_i (\widetilde{v}_i)_{x_1}^{X_i}-a(h_1-\widetilde{h})F_i(\widetilde{h}_{i})_{x_1}^{X_i} \, dx_1 \, dx'$,
we first find from $\eqref{presentF}$ and the assumption \eqref{perturbation_small} on that 
\begin{align} \label{Fbound}
||F_i||_{L^{\infty}} \le C\left(||\rho-\widetilde{\rho}||_{L^{\infty}}+||\widetilde{\rho}-\widetilde{\rho_i}^{X_i}||_{L^{\infty}}+||u-\widetilde{u}||_{L^{\infty}}+||\widetilde{u}-\widetilde{u_i}^{X_i}||_{L^{\infty}}\right) \le C(\varepsilon_1+\delta_0). 
\end{align}
Also, using the definition of $\rho$ and $h$ from \eqref{presentF}, we get 
\begin{align*}
&F_i=\sigma_i^*\frac{v-\widetilde{v}}{v}+\sigma_i^*\frac{\widetilde{v}-\widetilde{v}_i^{X_i}}{v}+\frac{h_1-\widetilde{h}}{v}+\frac{\widetilde{h}-\widetilde{h}_{i}^{X_i}}{v}\\
&\quad \ \ +(2\mu+\lambda)\frac{\partial_{x_1}(v-\widetilde{v})}{v}+(2\mu+\lambda)\frac{\partial_{x_1}(\widetilde{v}-\widetilde{v}_i^{X_i})}{v}. 
\end{align*}
Thus, by the equality above and the fact that $-\sigma_i^*(\widetilde{h}_{i})_{x_1}^{X_i}+p(\widetilde{v}_i)_{x_1}^{X_i}=0$, we have
\begin{align}
\begin{aligned} \label{expansionB3}
&F_iap'(\widetilde{v})(v-\widetilde{v})(\widetilde{v}_i)_{x_1}^{X_i}-F_ia(h_1-\widetilde{h})(\widetilde{h}_{i})_{x_1}^{X_i} \\
=&\frac{ap'(\widetilde{v})(\widetilde{v}_i)_{x_1}^{X_i}\sigma_i^*}{v}|v-\widetilde{v}|^2-\frac{ap'(\widetilde{v_i}^{X_i})(\widetilde{v}_i)_{x_1}^{X_i}}{v\sigma_i^*}|h_1-\widetilde{h}|^2 \\
&+\frac{ap'(\widetilde{v})(\widetilde{v}_i)_{x_1}^{X_i}\sigma_i^*}{v}(v-\widetilde{v})(\widetilde{v}-\widetilde{v}_i^{X_i})+\frac{ap'(\widetilde{v})(\widetilde{v}_i)_{x_1}^{X_i}}{v}(v-\widetilde{v})(\widetilde{h}-\widetilde{h}_{i}^{X_i}) \\
&-\frac{ap'(\widetilde{v_i}^{X_i})(\widetilde{v}_i)_{x_1}^{X_i}}{v}(h_1-\widetilde{h})(\widetilde{v}-\widetilde{v}_i^{X_i})-\frac{ap'(\widetilde{v_i}^{X_i})(\widetilde{v}_i)_{x_1}^{X_i}}{v \sigma_i^*}(h_1-\widetilde{h})(\widetilde{h}-\widetilde{h}_{i}^{X_i})\\
&+\frac{(2\mu+\lambda)ap'(\widetilde{v})(\widetilde{v}_i)_{x_1}^{X_i}\partial_{x_1}(v-\widetilde{v})(v-\widetilde{v})}{v}-\frac{(2\mu+\lambda)ap'(\widetilde{v_i}^{X_i})(\widetilde{v}_i)_{x_1}^{X_i}\partial_{x_1}(v-\widetilde{v})(h_1-\widetilde{h})}{\sigma_i^*v}\\ 
&+\frac{(2\mu+\lambda)ap'(\widetilde{v})(\widetilde{v}_i)_{x_1}^{X_i}\partial_{x_1}(\widetilde{v}-\widetilde{v_i}^{X_i})(v-\widetilde{v})}{v}-\frac{(2\mu+\lambda)ap'(\widetilde{v_i}^{X_i})(\widetilde{v}_i)_{x_1}^{X_i}\partial_{x_1}(\widetilde{v}-\widetilde{v_i}^{X_i})(h_1-\widetilde{h})}{\sigma_i^*v}.
\end{aligned}
\end{align}
Note that the two quadratic terms above in first equality right hand side. 

Using the Mean-Value theorem and the fact that $p'(\widetilde{v})<0<(\widetilde{v_i})_{x_1}^{X_i}\sigma_i^*$,
\begin{align*}
\frac{ap'(\widetilde{v})(\widetilde{v}_i)_{x_1}^{X_i}\sigma_i^*}{v}|v-\widetilde{v}|^2&=\frac{ap'(\widetilde{v})(\widetilde{v}_i)_{x_1}^{X_i}\sigma_i^*}{v}\frac{|p(v)-p(\widetilde{v})|^2}{|p'(c)|^2}\\
&\le -(1-C(\delta_0+\delta_i+\nu))\frac{\sigma_m|(\widetilde{v}_i)_{x_1}^{X_i}| |p(v)-p(\widetilde{v})|^2}{v_m |p'(v_m)|}.
\end{align*}
Using Young's inequality and the fact that $p'(\widetilde{v_i}^{X_i})<0<\frac{(\widetilde{v_i})_{x_1}^{X_i}}{\sigma_i^*}$,
\begin{align*}
-\frac{ap'(\widetilde{v_i}^{X_i})(\widetilde{v}_i)_{x_1}^{X_i}}{v\sigma_i^*}|h_1-\widetilde{h}|^2&\le (1+\boldsymbol{\kappa}) |\frac{ap'(\widetilde{v}_i^{X_i})}{ (\sigma_i^*)^3 v}(\widetilde{v}_i)_{x_1}^{X_i}||p(v)-p(\widetilde{v})|^2\\
&\quad +C(1+\frac{1}{\boldsymbol{\kappa}})|p'(\widetilde{v}_i^{X_i})||(\widetilde{v}_i)_{x_1}^{X_i}|\left| h_1-\widetilde{h}-\frac{p(v)-p(\widetilde{v})}{\sigma_i^*} \right|^2 \\
&\le (1+\boldsymbol{\kappa}+C(\delta_0+\varepsilon_1+\nu))|\frac{p'(v_m)(\widetilde{v}_i)_{x_1}^{X_i}}{\sigma_m^3v_m}||p(v)-p(\widetilde{v})|^2\\
&\quad +C(1+\frac{1}{\boldsymbol{\kappa}})|\frac{\delta_i}{\nu_i}(a_i)_{x_1}^{X_i}|\left| h_1-\widetilde{h}-\frac{p(v)-p(\widetilde{v})}{\sigma_i^*} \right|^2.
\end{align*}
Thus, using $\sigma_m=\sqrt{-p'(v_m)}$ and choosing for enough small $\kappa$, $\delta_0$, and $\varepsilon_1$,
\begin{align*}
 & \quad \ -\frac{ap'(\widetilde{v})\delta_i \sigma_i^* (a_i)_{x_1}^{X_i}}{v p'(\widetilde{v}_i^{X_i})\nu_i}|v-\widetilde{v}|^2+\frac{ap'(\widetilde{v})
\delta_i(a_i)_{x_1}^{X_i}}{v p'(\widetilde{v}_i^{X_i}) \sigma_i^* \nu_i}|h_1-\widetilde{h}|^2 \\ & \le (\boldsymbol{\kappa}+C(\delta_0+\varepsilon_1+\delta_i+\nu))\frac{|(\widetilde{v}_i)_{x_1}^{X_i}\sigma_m| |p(v)-p(\widetilde{v})|^2}{ v_m |p'(v_m)|}\\
&\quad +C(1+\frac{1}{\boldsymbol{\kappa}})\frac{\delta_i}{\nu_i}|(a_i)_{x_1}^{X_i}|\left| h_1-\widetilde{h}-\frac{p(v)-p(\widetilde{v})}{\sigma_i^*} \right|^2 \\
&\le \frac{C_1|(\widetilde{v}_i^{X_i})_{x_1}| |p(v)-p(\widetilde{v})|^2  }{64}+C\frac{\delta_i}{\nu_i}|(a_i)_{x_1}^{X_i}|\left| h_1-\widetilde{h}-\frac{p(v)-p(\widetilde{v})}{\sigma_i^*} \right|^2.
\end{align*}
Hence, applying results above with the facts $\partial_{x_1}(v-\widetilde{v})<\varepsilon_1<C$ and $\partial_{x_1}(\widetilde{v}-\widetilde{v_i}^{X_i})<\delta_0<C$ to the last line in \eqref{expansionB3} for each $i=1,2$, 
\begin{align}
\begin{aligned} \label{B3first}
\mathcal{B}_3\le&\frac{C_1}{80}\sum\limits_{i=1}^2 \int_{\mathbb{T}^2} \int_{\mathbb{R}} |(\widetilde{v}_i)_{x_1}^{X_i}| |p(v)-p(\widetilde{v})|^2 \, dx_1 \, dx'+\frac{1}{80}\sum\limits_{i=1}^2 \int_{\mathbb{T}^2} \int_{\mathbb{R}} |\sigma_i^*| |(a_i)_{x_1}^{X_i}|\left| h_1-\widetilde{h}-\frac{p(v)-p(\widetilde{v})}{\sigma_i^*} \right|^2 \, dx_1 \, dx' \\
&+ C\sum\limits_{i=1}^2 \int_{\mathbb{T}^2} \int_{\mathbb{R}} |(\widetilde{v}_i)_{x_1}^{X_i}| |v-\widetilde{v}| (|\widetilde{v}-\widetilde{v}_i^{X_i}|+|\widetilde{h}-\widetilde{h}_{i}^{X_i}|) \, dx_1 \, dx' \\
&+ C\sum\limits_{i=1}^2 \int_{\mathbb{T}^2} \int_{\mathbb{R}} |(\widetilde{v}_i)_{x_1}^{X_i}| |h_1-\widetilde{h}| (|\widetilde{v}-\widetilde{v}_i^{X_i}|+|\widetilde{h}-\widetilde{h}_{i}^{X_i}|) \, dx_1 \, dx'\\
&+ C\sum\limits_{i=1}^2 \int_{\mathbb{T}^2} \int_{\mathbb{R}} |(\widetilde{v}_i)_{x_1}^{X_i}| |\partial_{x_1}(v-\widetilde{v})| \left| v-\widetilde{v}- \frac{h_1-\widetilde{h}}{\sigma_i^*} \right| \, dx_1 \, dx' \\
&+ C\sum\limits_{i=1}^2 \int_{\mathbb{T}^2} \int_{\mathbb{R}} |(\widetilde{v}_i)_{x_1}^{X_i}| |\partial_{x_1}(\widetilde{v}-\widetilde{v}_i^{X_i})| \left| v-\widetilde{v}- \frac{h_1-\widetilde{h}}{\sigma_i^*} \right| \, dx_1 \, dx'. 
\end{aligned}
\end{align}
Note that 
\begin{align*} 
\left| v-\widetilde{v}- \frac{h_1-\widetilde{h}}{\sigma_i^*} \right| &\le \left| p(v)-p(\widetilde{v})- \frac{h_1-\widetilde{h}}{\sigma_i^*} \right|+\left| v-\widetilde{v}\right|+\left| p(v)-p(\widetilde{v}) \right|   \\
&\le C\left( \left| h_1-\widetilde{h}-\frac{p(v)-p(\widetilde{v})}{\sigma_i^*} \right|+\left| p(v)-p(\widetilde{v}) \right| \right). \\
\end{align*}
In addition, by the definition of $p(v)=v^{-\gamma}$,
\begin{align*}
\left| \partial_{x_1}(v-\widetilde{v}) \right|&=\left| \frac{\partial_{x_1}(p(v)-p(\widetilde{v}))}{p'(v)}+p(\widetilde{v})_{x_1}\left( \frac{1}{p'(v)}-\frac{1}{p'(\widetilde{v})} \right) \right| \\
&\le C\left( \left| \nabla_x(p(v)-p(\widetilde{v})) \right|+|(\widetilde{v})_{x_1}||p(v)-p(\widetilde{v})| \right).
\end{align*}
Then
\begin{align}
\begin{aligned} \label{B3second}
&\sum\limits_{i=1}^2 \int_{\mathbb{T}^2} \int_{\mathbb{R}} |(\widetilde{v}_i)_{x_1}^{X_i}| |\partial_{x_1}(v-\widetilde{v})| \left| v-\widetilde{v}- \frac{h_1-\widetilde{h}}{\sigma_i^*} \right| \, dx_1 \, dx' \\
\le & C\sum\limits_{i=1}^2 \int_{\mathbb{T}^2} \int_{\mathbb{R}} |(\widetilde{v}_i)_{x_1}^{X_i}| (\left| \nabla_x(p(v)-p(\widetilde{v})) \right|+ |(\widetilde{v})_{x_1}| |p(v)-p(\widetilde{v})|) (\left| h_1-\widetilde{h}-\frac{p(v)-p(\widetilde{v})}{\sigma_i^*} \right|+|p(v)-p(\widetilde{v})|)\, dx_1 \, dx' \\
\le &C\sum\limits_{i=1}^2 \int_{\mathbb{T}^2} \int_{\mathbb{R}} |(\widetilde{v}_i)_{x_1}^{X_i}| \left| \nabla_x(p(v)-p(\widetilde{v})) \right|\left| h_1-\widetilde{h}-\frac{p(v)-p(\widetilde{v})}{\sigma_i^*} \right|\, dx_1 \, dx' \\
&+C\sum\limits_{i=1}^2 \int_{\mathbb{T}^2} \int_{\mathbb{R}} |(\widetilde{v}_i)_{x_1}^{X_i}| \left| \nabla_x(p(v)-p(\widetilde{v})) \right|\left|p(v)-p(\widetilde{v}) \right|\, dx_1 \, dx' \\
&+C(\delta_1^2+\delta_2^2)\sum\limits_{i=1}^2 \int_{\mathbb{T}^2} \int_{\mathbb{R}} |(\widetilde{v}_i)_{x_1}^{X_i}| \left| p(v)-p(\widetilde{v}) \right| \left| h_1-\widetilde{h}-\frac{p(v)-p(\widetilde{v})}{\sigma_i^*} \right|\, dx_1 \, dx' \\
&+C(\delta_1^2+\delta_2^2)\sum\limits_{i=1}^2 \int_{\mathbb{T}^2} \int_{\mathbb{R}} |(\widetilde{v}_i)_{x_1}^{X_i}| \left| p(v)-p(\widetilde{v}) \right|^2 \, dx_1 \, dx'.
\end{aligned}
\end{align}
By Young's inequality,
\begin{align*}
&C\sum\limits_{i=1}^2 \int_{\mathbb{T}^2} \int_{\mathbb{R}} |(\widetilde{v}_i)_{x_1}^{X_i}| \left| \nabla_x(p(v)-p(\widetilde{v})) \right|\left| h_1-\widetilde{h}-\frac{p(v)-p(\widetilde{v})}{\sigma_i^*} \right|\, dx_1 \, dx' \\ 
\le&C\delta_0 \int_{\mathbb{T}^2} \int_{\mathbb{R}} \left| \nabla_x(p(v)-p(\widetilde{v})) \right|^2 \, dx_1 \, dx'+C\sum\limits_{i=1}^2 \int_{\mathbb{T}^2} \int_{\mathbb{R}} |(\widetilde{v}_i)_{x_1}^{X_i}| \left| h_1-\widetilde{h}-\frac{p(v)-p(\widetilde{v})}{\sigma_i^*} \right|^2 \, dx_1 \, dx' \\
\le& C \sqrt{\delta_0}(\mathcal{G}_1+\mathcal{D}).
\end{align*}
Apply the same way to other terms in \eqref{B3second} and to last line in \eqref{B3first}.

Also apply \eqref{specialterm} to first line and Lemma \ref{shock interaction lemma1} to second and third line in \eqref{B3first}. 

Then using the smallness of $\delta_0$, we get 
\[\mathcal{B}_3\le \frac{1}{40}(C_1 \mathcal{G}^S+\mathcal{G}_1+\mathcal{D})+C\sum\limits_{i=1}^2 \delta_i^2e^{-C\delta_it} \int_{\mathbb{T}^2}  \int_{\mathbb{R}} \eta(U|\widetilde{U}) \, dx_1 \, dx'+C\varepsilon_1\delta_1\delta_2e^{-C\min\left\{ \delta_1,\delta_2 \right\}t}.\]
\noindent $\bullet$ {\bf(Estimate of $\mathcal{B}_4$):}
Once again using \eqref{presentF} we have,
\begin{align*}
\mathcal{B}_4&=\sum\limits_{i=1}^2 \int_{\mathbb{T}^2}  \int_{\mathbb{R}}  \left( Q(v|\widetilde{v})+\frac{|\bold{h}-\widetilde{\bold{h}}|^2}{2} \right) F_i (a_i)_{x_1}^{X_i}\, dx_1 \, dx' \\
&\le C\sum\limits_{i=1}^2 \int_{\mathbb{T}^2} \int_{\mathbb{R}} |F_i| |(a_i)_{x_1}^{X_i}| \left( |p(v)-p(\widetilde{v})|^2+\left| h_1-\widetilde{h}-\frac{p(v)-p(\widetilde{v})}{\sigma_i^*} \right|^2+\frac{h_2^2+h_3^2}{2} \right)  \, dx_1 \, dx'\\
&\le C\sum\limits_{i=1}^2 \int_{\mathbb{T}^2} \int_{\mathbb{R}} |(a_i)_{x_1}^{X_i}| |p(v)-p(\widetilde{v})|^3 \, dx_1 \, dx'+ C\sum\limits_{i=1}^2 \int_{\mathbb{T}^2} \int_{\mathbb{R}} |\widetilde{v}-\widetilde{v}_i^{X_i}| |(a_i)_{x_1}^{X_i}| |p(v)-p(\widetilde{v})|^2 \, dx_1 \, dx'\\ 
&+C\sum\limits_{i=1}^2 \int_{\mathbb{T}^2} \int_{\mathbb{R}} |(a_i)_{x_1}^{X_i}| |h_1-\widetilde{h}| |p(v)-p(\widetilde{v})|^2 \, dx_1 \, dx'+C\sum\limits_{i=1}^2 \int_{\mathbb{T}^2} \int_{\mathbb{R}} |(a_i)_{x_1}^{X_i}| |\widetilde{h}-\widetilde{h}_{i}^{X_i}| |p(v)-p(\widetilde{v})|^2 \, dx_1 \, dx' \\
&+C\sum\limits_{i=1}^2 \int_{\mathbb{T}^2} \int_{\mathbb{R}} |(a_i)_{x_1}^{X_i}| |\partial_{x_1}(v-\widetilde{v})| |p(v)-p(\widetilde{v})|^2 \, dx_1 \, dx'\\
&+C\sum\limits_{i=1}^2 \int_{\mathbb{T}^2} \int_{\mathbb{R}} |(a_i)_{x_1}^{X_i}| |\partial_{x_1}(\widetilde{v}-\widetilde{v}_i^{X_i})| |p(v)-p(\widetilde{v})|^2 \, dx_1 \, dx'\\
&+C(\varepsilon_1+\delta_0)(\mathcal{G}_1+\mathcal{G}_3).
\end{align*}
Note that for the right hand side, the first and second terms already were estimated in \eqref{Ki1} and \eqref{estimateK2}.

On the other hand, for third term, note that
\begin{align*}
&\sum\limits_{i=1}^2 \int_{\mathbb{T}^2} \int_{\mathbb{R}} |(a_i)_{x_1}^{X_i}| |h_1-\widetilde{h}| |p(v)-p(\widetilde{v})|^2 \, dx_1 \, dx' \\
\le &C\underbrace{\sum\limits_{i=1}^2 \int_{\mathbb{T}^2} \int_{\mathbb{R}} |(a_i)_{x_1}^{X_i}| \left| h_1-\widetilde{h}-\frac{p(v)-p(\widetilde{v})}{\sigma_i^*} \right| |p(v)-p(\widetilde{v})|^2 \, dx_1 \, dx'}_{:=\mathcal{B}_{4,3,1}}\\
&+ \underbrace{C\sum\limits_{i=1}^2 \int_{\mathbb{T}^2} \int_{\mathbb{R}} |(a_i)_{x_1}^{X_i}| |p(v)-p(\widetilde{v})|^3 \, dx_1 \, dx'}_{=\mathcal{K}_1 \ in \ section \ 4.5}.
\end{align*}
Already we got the estimate on $\mathcal{K}_1$ in section 4.5 :
\begin{equation*}
\mathcal{K}_1\leq  \varepsilon_1(\mathcal{D}+C_1 \mathcal{G}^S ) +  C \sum\limits_{i=1}^2 \varepsilon_1 \delta_i \nu_i e^{-C\delta_it} \int_{\mathbb{T}^2}  \int_{\mathbb{R}} \eta(U|\widetilde{U}) \, dx_1 \, dx'.    
\end{equation*}
To estimate $\mathcal{B}_{4,3,1}$, use Lemma \ref{infty interpolation inequality lemma} as
\begin{align*}
\mathcal{B}_{4,3,1}\le &C\sum\limits_{i=1}^2 ||p(v)-p(\widetilde{v})||^2_{L^\infty} \sqrt{\mathcal{G}_1} \sqrt{\frac{\nu_i}{\delta_i}} \sqrt{\int_{\mathbb{T}^2} \int_{\mathbb{R}} |(v_i)_{x_1}^{X_i}| \, dx_1 \, dx'} \\ 
\le &C\sum\limits_{i=1}^2\left(||p(v)-p(\widetilde{v})||_{L^2}+||\nabla_x^2(p(v)-p(\widetilde{v}))||_{L^2}\right)||\nabla_x(p(v)-p(\widetilde{v})||_{L^2} \sqrt{\mathcal{G}_1} \sqrt{\nu_i}\\
\le &C\varepsilon_1\sum\limits_{i=1}^2 \sqrt{\mathcal{D}}\sqrt{\mathcal{G}_1} \\
\le &\varepsilon_1(\mathcal{D}+\mathcal{G}_1).\\
\end{align*}
For 4th term, we can easily get
\begin{align*}
&\sum\limits_{i=1}^2 \int_{\mathbb{T}^2} \int_{\mathbb{R}} |(a_i)_{x_1}^{X_i}| |\widetilde{h}-\widetilde{h}_{i}^{X_i}| |p(v)-p(\widetilde{v})|^2 \, dx_1 \, dx'\\
\le &C\sum\limits_{i=1}^2 \int_{\mathbb{T}^2} \int_{\mathbb{R}} |(a_i)_{x_1}^{X_i}| |\widetilde{v}-\widetilde{v}_{i}^{X_i}| |p(v)-p(\widetilde{v})|^2 \, dx_1 \, dx'\\
=&C\mathcal{K}_2 \ (in \ section \ 4.5)\\
\le& C \delta_1 \delta_2 e^{-C\min\left\{ \delta_1,\delta_2 \right\}t} \sum\limits_{i=1}^2 \nu_i \int_{\mathbb{T}^2} \int_{\mathbb{R}} \eta(U|\widetilde{U}) \, dx_1 \, dx'.    
\end{align*}
Finally, for 5th and 6th terms, as in $\mathcal{B}_3$, we get
\begin{align*}
 &\sum\limits_{i=1}^2 \int_{\mathbb{T}^2} \int_{\mathbb{R}} |(a_i)_{x_1}^{X_i}| |\partial_{x_1}(v-\widetilde{v})| |p(v)-p(\widetilde{v})|^2 \, dx_1 \, dx' \\
\le &\varepsilon_1(\mathcal{D}+C_1\mathcal{G}^S)+C\sum\limits_{i=1}^2 \varepsilon_1 \nu_i \delta_i e^{-C\delta_it} \int_{\mathbb{T}^2} \int_{\mathbb{R}} \eta(U|\widetilde{U}) \, dx_1 \, dx', \\ 
&\sum\limits_{i=1}^2 \int_{\mathbb{T}^2} \int_{\mathbb{R}} |(a_i)_{x_1}^{X_i}| |\partial_{x_1}(\widetilde{v}-\widetilde{v_i}^{X_i})| |p(v)-p(\widetilde{v})|^2 \, dx_1 \, dx' \\
\le &\delta_0\mathcal{G}^S + C \delta_0^2 \sum\limits_{i=1}^2 \nu_i \delta_i e^{-C\delta_it}.
\end{align*}
Thus, combining all the estimates above and using the smallness of $\delta_0$ and $\varepsilon_1$, we conclude that
\begin{align*}
\mathcal{B}_4\le &\frac{1}{40}\left(\mathcal{G}_1 + \mathcal{D} + C_1\mathcal{G}^S  + \mathcal{G}_3\right) \\ 
&+C\left( \sum\limits_{i=1}^2 \varepsilon_1 \nu_i \delta_i e^{-C\delta_it}+\delta_1\delta_2 e^{-C\min\left\{ \delta_1, \delta_2  \right\}t} \right) \int_{\mathbb{T}^2} \int_{\mathbb{R}} \eta(U|\widetilde{U}) \, dx_1 \, dx'+C \delta_0^2 \sum\limits_{i=1}^2 \nu_i \delta_i e^{-C\delta_it}.
\end{align*}
\noindent $\bullet${\bf (Estimate for $\mathcal{B}_5$, $\mathcal{B}_6$, and $\mathcal{B}_7$):}
As above, we can easily get some estimates of $\mathcal{B}_5$, $\mathcal{B}_6$, and $\mathcal{B}_7$.
\begin{align*}
\mathcal{B}_5\le& C \int_{\mathbb{T}^2} \int_{\mathbb{R}} |(\widetilde{v})_{x_1}| |\nabla_x(p(v)-p(\widetilde{v}))| |p(v)-p(\widetilde{v})| \, dx_1 \, dx' \\
\le& C \sum\limits_{i=1}^2 \int_{\mathbb{T}^2} \int_{\mathbb{R}} |(\widetilde{v}_i)_{x_1}^{X_i}|^{\frac{1}{2}} |\nabla_x(p(v)-p(\widetilde{v}))|^2 \, dx_1 \, dx'+ C \sum\limits_{i=1}^2 \int_{\mathbb{T}^2} \int_{\mathbb{R}} |(\widetilde{v}_i)_{x_1}^{X_i}|^{\frac{3}{2}} |p(v)-p(\widetilde{v})|^2 \, dx_1 \, dx'\\
 \le& \frac{1}{40}(\mathcal{D}+C_1\mathcal{G}^S)+C\sum\limits_{i=1}^2 \delta_i^2e^{-C\delta_it} \int_{\mathbb{T}^2}  \int_{\mathbb{R}} \eta(U|\widetilde{U}) \, dx_1 \, dx',\\
\mathcal{B}_6\le& C\sum\limits_{i=1}^2  \int_{\mathbb{T
}^2}  \int_{\mathbb{R}} |\frac{\nu_i(\widetilde{v}_i)_{x_1}^{X_i}}{\delta_i}| |p(v)-p(\widetilde{v})| |\partial_{x_1}(p(v)-p(\widetilde{v}))| \, dx_1 \, dx' \\
\le& C\sum\limits_{i=1}^2 \frac{\nu_i^2}{\delta_i^2} \int_{\mathbb{T
}^2}  \int_{\mathbb{R}} |(\widetilde{v}_i)_{x_1}^{X_i}|^{\frac{7}{4}} |p(v)-p(\widetilde{v})|^2   \, dx_1 \, dx'+C\sum\limits_{i=1}^2 \int_{\mathbb{T
}^2}  \int_{\mathbb{R}} |(\widetilde{v}_i)_{x_1}^{X_i}|^{\frac{1}{4}} |\partial_{x_1}(p(v)-p(\widetilde{v}))|^2   \, dx_1 \, dx'\\
\le& C\sum\limits_{i=1}^2 \sqrt{\delta_i}  \int_{\mathbb{T
}^2}  \int_{\mathbb{R}} |(\widetilde{v}_i)_{x_1}^{X_i}| |p(v)-p(\widetilde{v})|^2   \, dx_1 \, dx'+C\sqrt{\delta_i}\sum\limits_{i=1}^2 \int_{\mathbb{T
}^2}  \int_{\mathbb{R}} |\nabla_x(p(v)-p(\widetilde{v}))|^2   \, dx_1 \, dx' \\
\le& \frac{1}{40}(\mathcal{D}+C_1\mathcal{G}^S)+C\sum\limits_{i=1}^2 \delta_i^2 e^{-C\delta_it}\int_{\mathbb{T}^2} \int_{\mathbb{R}} \eta(U|\widetilde{U}) \, dx_1 \, dx', \\
\mathcal{B}_7\le& C\sum\limits_{i=1}^2 \frac{\nu_i}{\delta_i} \int_{\mathbb{T}^2} \int_{\mathbb{R}} |(\widetilde{v}_i)_{x_1}^{X_i}|^2 |p(v)-p(\widetilde{v})|^2 \, dx_1 \, dx' \\
\le& C\sum\limits_{i=1}^2 \delta_i \nu_i \int_{\mathbb{T}^2} \int_{\mathbb{R}} |(\widetilde{v}_i)_{x_1}^{X_i}| |p(v)-p(\widetilde{v})|^2 \, dx_1 \, dx'\\
\le& \frac{C_1}{40}\mathcal{G}^S+C\sum\limits_{i=1}^2 \delta_i^2 e^{-C\delta_it}\int_{\mathbb{T}^2} \int_{\mathbb{R}} \eta(U|\widetilde{U}) \, dx_1 \, dx'.
\end{align*}
Thus, combining the estimates above, we have
\[\mathcal{B}_5+\mathcal{B}_6+\mathcal{B}_7\le \frac{3}{40}(C_1\mathcal{G}^S+\mathcal{D})+C\sum\limits_{i=1}^2 \delta_i^2 e^{-C\delta_it}\int_{\mathbb{T}^2} \int_{\mathbb{R}} \eta(U|\widetilde{U}) \, dx_1 \, dx'.\]
\noindent $\bullet${\bf (Estimate for $\mathcal{B}_8$): }
Note that
\[R^*=\frac{2\mu+\lambda}{v}(\nabla_x \bold{u} \nabla_x v- div_x\bold{u} \nabla_x v)-\mu \nabla_x \times \nabla_x \times \bold{u}+\nabla_x( -p(\widetilde{v})+p(\widetilde{v}_1^{X_1})+p(\widetilde{v}_2^{X_2}) ).\]
Decompose $\mathcal{B}_8$ as
\begin{align*}
\mathcal{B}_8=&(2\mu+\lambda)\int_{\mathbb{T}^2} \int_{\mathbb{R}} a(\bold{h}-\widetilde{\bold{h}})\cdot \frac{(\nabla_x \bold{u}  \nabla_x v - div_x\bold{u} \nabla_x v)}{v} \, dx_1 \, dx'\\
&-\mu \int_{\mathbb{T}^2} \int_{\mathbb{R}} a(\bold{h}-\widetilde{\bold{h}})\cdot \nabla_x \times  \nabla_x \times \bold{u} \, dx_1 \, dx'\\
&+ \int_{\mathbb{T}^2} \int_{\mathbb{R}} a(\bold{h}-\widetilde{\bold{h}}) \cdot \nabla_x( -p(\widetilde{v})+p(\widetilde{v}_1^{X_1})+p(\widetilde{v}_2^{X_2}) ) \, dx_1 \, dx' \\
=:&\mathcal{B}_{8,1}+\mathcal{B}_{8,2}+\mathcal{B}_{8,3}
\end{align*}
 To get some bound of $\mathcal{B}_{8,1}$, decompose this into 
 \begin{align*}
 \mathcal{B}_{8,1}&= (2\mu+\lambda)\int_{\mathbb{T}^2}\int_{\mathbb{R}} a(h_1-\widetilde{h}) \frac{\partial_{x_1}\bold{u}\cdot \nabla v-div\bold{u}\partial_{x_1}v}{v} dx_1dx'\\
 & \quad +(2\mu+\lambda)\int_{\mathbb{T}^2}\int_{\mathbb{R}} ah_2 \frac{\partial_{x_2}\bold{u}\cdot \nabla v-div\bold{u}\partial_{x_2}v}{v} dx_1dx'+(2\mu+\lambda)\int_{\mathbb{T}^2}\int_{\mathbb{R}} ah_3 \frac{\partial_{x_3}\bold{u}\cdot \nabla v-div\bold{u}\partial_{x_3}v}{v} dx_1dx' \\
 & =:\mathcal{B}_{8,1,1}+\mathcal{B}_{8,1,2}+\mathcal{B}_{8,1,3}.
 \end{align*}

\noindent First, to estimate $\mathcal{B}_{8,1,1}$, let $\bold{u}':=(u_2,u_3)$, $\nabla_{x'}:=(\partial_{x_2}, \partial_{x_3})$, and $\nabla_{x'}\cdot \bold{u}':=\partial_{x_2}u_2+\partial_{x_3}u_3$.

\noindent Then
\begin{equation*}
\begin{aligned}
& \left|\partial_{x_1}\bold{u}\cdot \nabla_x v-div_x\bold{u}\partial_{x_1} v\right| \\    
= &  \left|\partial_{x_1} \bold{u}'\cdot \nabla_{x'}v-\nabla_{x'} \cdot \bold{u}' \partial_{x_1}v\right|\\ 
= &  \left|\partial_{x_1}\bold{u}'\cdot \frac{\nabla_{x'}p(v)}{-\gamma p(v)^{1+\frac{1}{\gamma}}}-\nabla_{x'}\cdot \bold{u}' \frac{\partial_{x_1}(p(v)-p(\widetilde{v}))}{-\gamma p(v)^{1+\frac{1}{\gamma}}}-\nabla_{x'}\cdot \bold{u}' \frac{\partial_{x_1}p(\widetilde{v})}{-\gamma p(v)^{1+\frac{1}{\gamma}}}\right| \\
\le & C\left(\left|\nabla_x(\bold{u}-\widetilde{\bold{u}})\right|\left|\nabla_x(p(v)-p(\widetilde{v})\right|+\left|\nabla_x(\bold{u}-\widetilde{\bold{u}})\right|\left|(\widetilde{v})_{x_1}\right|\right),
\end{aligned}    
\end{equation*}
where we used $\widetilde{\bold{u}}=(\widetilde{u},0,0)$.

\noindent Thus, 
\begin{align}
\begin{aligned}\label{B8.1.1}
\mathcal{B}_{8,1,1}\le& C\int_{\mathbb{T}^2} \int_{\mathbb{R}} |h_1-\widetilde{h}| (|\nabla_x(\bold{u}-\widetilde{\bold{u}})| |\nabla_x(p(v)-p(\widetilde{v}))|+|\nabla_x(\bold{u}-\widetilde{\bold{u}})| |(\widetilde{v})_{x_1}| )\, dx_1 \, dx' \\
\le& C\int_{\mathbb{T}^2} \int_{\mathbb{R}} |h_1-\widetilde{h}| |\nabla_x(\bold{u}-\widetilde{\bold{u}})| |\nabla_x(p(v)-p(\widetilde{v}))|\, dx_1 \, dx' \\
&+C\sum\limits_{i=1}^2 \frac{\delta_i}{\nu_i} \int_{\mathbb{T}^2} \int_{\mathbb{R}} \left|h_1-\widetilde{h}-\frac{p(v)-p(\widetilde{v})}{\sigma_i^*}\right| |\nabla_x(\bold{u}-\widetilde{\bold{u}})| |(a_i)_{x_1}^{X_i}| \, dx_1 \, dx'\\
&+C\sum\limits_{i=1}^2 \int_{\mathbb{T}^2} \int_{\mathbb{R}} |p(v)-p(\widetilde{v})| |\nabla_x(\bold{u}-\widetilde{\bold{u}})| |(\widetilde{v}_i)_{x_1}^{X_i}| \, dx_1 \, dx'\\
\le & \underbrace{C||h_1-\widetilde{h}||_{L^3} ||\nabla_x(\bold{u}-\widetilde{\bold{u}})||_{L^6} \sqrt{\mathcal{D}}}_{:=\mathcal{S}}\\
&+C\sum\limits_{i=1}^2 \int_{\mathbb{T}^2} \int_{\mathbb{R}} \left|h_1-\widetilde{h}-\frac{p(v)-p(\widetilde{v})}{\sigma_i^*}\right|^2  |(a_i)_{x_1}^{X_i}|^{2} \, dx_1 \, dx'+C\sum\limits_{i=1}^2 \frac{\delta_i^2}{\nu_i^2} \int_{\mathbb{T}^2} \int_{\mathbb{R}} |\nabla_x(\bold{u}-\widetilde{\bold{u}})|^2  \, dx_1 \, dx'\\
&+C\sum\limits_{i=1}^2 \int_{\mathbb{T}^2} \int_{\mathbb{R}} |p(v)-p(\widetilde{v})|^2  |(\widetilde{v}_i)_{x_1}^{X_i}|^{\frac{3}{2}} \, dx_1 \, dx'+C\sum\limits_{i=1}^2 \int_{\mathbb{T}^2} \int_{\mathbb{R}} |\nabla_x(\bold{u}-\widetilde{\bold{u}})|^2 |(\widetilde{v}_i)_{x_1}^{X_i}|^{\frac{1}{2}} \, dx_1 \, dx'.
\end{aligned}
\end{align}
To estimate $\mathcal{S}$, we will use the following inequalities : for any $f : \mathbb{T}^2\times \mathbb{R} \longrightarrow \mathbb{R}(\text{or} \ \mathbb{R}^3)$ belonging to $H^1$, it holds from Gagliardo-Nirenberg interpolation inequality that
\begin{align} \label{GNIAPP1}
\lVert f \rVert_{L^3}\le C\sqrt{ \lVert f \rVert_{L^6}}\sqrt{ \lVert f \rVert_{L^2}}.   
\end{align}
On the other hand, using Gagliardo-Nirenberg inequality, we have
\begin{align} \label{GNIAPP2}
\lVert f \rVert_{L^6}\le C\lVert \nabla_xf \rVert_{L^2}.    
\end{align}
Combining \eqref{GNIAPP1} and \eqref{GNIAPP2}, we get
\begin{align} \label{GNIAPP3}
\lVert f \rVert_{L^3}\le C \lVert f \rVert_{H^1}.    
\end{align}
Using \eqref{perturbation_small} and applying \eqref{GNIAPP3} to $h_1-\widetilde{h}$ and \eqref{GNIAPP2} to $\nabla_x(\bold{u}-\widetilde{\bold{u}})$, we have
\[\mathcal{S}\le C\varepsilon_1\left(||\nabla_x(\bold{u}-\widetilde{\bold{u}})||_{H^1}^2+\mathcal{D}\right).\]
Therefore, from \eqref{B8.1.1},
\begin{align*}
\mathcal{B}_{8,1,1}\le&  C\varepsilon_1(||\nabla_x(\bold{u}-\widetilde{\bold{u}})||_{H^1}^2+\mathcal{D})+C\sum\limits_{i=1}^2 \delta_i^{\frac{3}{2}}\mathcal{G}_1+C\sum\limits_{i=1}^2 \delta_i ||\nabla_x(\bold{u}-\widetilde{\bold{u}})||_{L^2}^2\\
& \ \ + C\sum\limits_{i=1}\delta_i\mathcal{G}^S+C\sum\limits_{i=1}^2 \delta_i^2 e^{-C\delta_it}\int_{\mathbb{T}^2} \int_{\mathbb{R}} \eta(U|\widetilde{U}) \, dx_1 \, dx'\\
\le &(\delta_0+\varepsilon_1)\left(C_1\mathcal{G}^S+\mathcal{G}_1+\mathcal{D}+||\nabla_x(\bold{u}-\widetilde{\bold{u}})||_{H^1}^2\right)+C\sum\limits_{i=1}^2 \delta_i^2 e^{-C\delta_it}\int_{\mathbb{T}^2} \int_{\mathbb{R}} \eta(U|\widetilde{U}) \, dx_1 \, dx'.    
\end{align*}
Likewise, we get 
\begin{align*}
 \mathcal{B}_{8,1,2}+\mathcal{B}_{8,1,3} &\le C \left(||h_2||_{L^3}+||h_3||_{L^3}\right)||\nabla_x(\bold{u}-\widetilde{\bold{u}})||_{L^6}\sqrt{\mathcal{D}}+C\sum\limits_{i=1}^2\delta_i^{\frac{5}{4}}\sqrt{\mathcal{G}_3}\left(||\nabla_x(\bold{u}-\widetilde{\bold{u}})||_{L^2}+\sqrt{\mathcal{D}}\right)\\
&\le C(\delta_0+\varepsilon_1)\left(\mathcal{G}_3+\mathcal{D}+||\nabla_x(\bold{u}-\widetilde{\bold{u}})||_{H^1}^2\right).   
\end{align*}
Altogether,
\[
\mathcal{B}_{8,1}\le (\delta_0+\varepsilon_1)\left(C_1\mathcal{G}^S+\mathcal{G}_1+\mathcal{G}_3+\mathcal{D}+||\nabla_x(\bold{u}-\widetilde{\bold{u}})||_{H^1}^2\right)+C\sum\limits_{i=1}^2\delta_i^2 e^{-C\delta_it}\int_{\mathbb{T}^2} \int_{\mathbb{R}} \eta(U|\widetilde{U}) \, dx_1 \, dx'.    
\]
To estimates $\mathcal{B}_{8,2}$, we decompose it again.
\begin{align*}
\mathcal{B}_{8,2}&=-\mu \int_{\mathbb{T}^2} \int_{\mathbb{R}} a(\bold{u}-\widetilde{\bold{u}})\cdot \nabla_x \times \nabla_x \times \bold{u} \, dx_1 \, dx'+\mu(2\mu+\lambda)\int_{\mathbb{T}^2} \int_{\mathbb{R}} a\nabla_x(v-\widetilde{v})\cdot \nabla_x \times \nabla_x \times \bold{u} \, dx_1 \, dx'\\
&=:\mathcal{B}_{8,2,1}+\mathcal{B}_{8,2,2}.
\end{align*}
Using the integration by parts and the fact that $\nabla_x \times \widetilde{u}\equiv 0$,
\begin{align} 
\begin{aligned}\label{B8.2.1}
\mathcal{B}_{8,2,1}&=-\mu \int_{\mathbb{T}^2} \int_{\mathbb{R}}\nabla_x\times \left(a(\bold{u}-\widetilde{\bold{u}})\right)\cdot \nabla_x \times \bold{u} \, dx_1 \, dx'  \\
&=-\mu \int_{\mathbb{T}^2} \int_{\mathbb{R}}\nabla_x\times \left(a(\bold{u}-\widetilde{\bold{u}})\right)\cdot \nabla_x \times \left(\bold{u}-\widetilde{\bold{u}}\right) \, dx_1 \, dx' \\
&=-\mu \int_{\mathbb{T}^2} \int_{\mathbb{R}}\left[\nabla_xa\times\left(\bold{u}-\widetilde{\bold{u}}\right)+a\nabla_x\times \left(\bold{u}-\widetilde{\bold{u}}\right) \right]\cdot \nabla_x \times \left(\bold{u}-\widetilde{\bold{u}}\right) \, dx_1 \, dx' \\
&=-\mu\int_{\mathbb{T}^2}\int_{\mathbb{R}} a|\nabla_x \times (\bold{u}-\widetilde{\bold{u}})|^2 \, dx_1 \, dx'\\
& \ \ \ -\mu\sum\limits_{i=1}^2 \int_{\mathbb{T}^2} \int_{\mathbb{R}} (a_i)_{x_1}^{X_i} (u_2(\partial_{x_1}u_2-\partial_{x_2}u_1)-u_3(\partial_{x_3}u_1-\partial_{x_1}u_3))\, dx_1 \, dx' \\
&\le -\mu\int_{\mathbb{T}^2}\int_{\mathbb{R}} a|\nabla_x \times (\bold{u}-\widetilde{\bold{u}})|^2 \, dx_1 \, dx' \\
& \ \ \ +\mu\sum\limits_{i=1}^2 \sqrt{\nu_i\delta_i}\int_{\mathbb{T}^2} \int_{\mathbb{R}} \left|(a_i)^{X_i}_{x_1}\right| (u_2^2+u_3^2)\, dx_1 \, dx' \\
& \ \ \ +\mu\sum\limits_{i=1}^2 \sqrt{\nu_i\delta_i}\int_{\mathbb{T}^2} \int_{\mathbb{R}} ( |\partial_{x_1}u_2-\partial_{x_2}u_1|^2+|\partial_{x_3}u_1-\partial_{x_1}u_3|^2)\, dx_1 \, dx'.\\
\end{aligned}
\end{align}
Substituting $\bold{u}=\bold{h}+(2\mu+\lambda)\nabla_xv$ into the second term of the r.h.s. in \eqref{B8.2.1}, we have
\begin{align*}
\mathcal{B}_{8,2,1}\le&-\mu\int_{\mathbb{T}^2}\int_{\mathbb{R}} a|\nabla_x \times (\bold{u}-\widetilde{\bold{u}})|^2 \, dx_1 \, dx' \\
&+C\sum\limits_{i=1}^2 \sqrt{\nu_i\delta_i}\int_{\mathbb{T}^2} \int_{\mathbb{R}} \left|(a_i)^{X_i}_{x_1}\right| (h_2^2+h_3^2)\, dx_1 \, dx'+C\sum\limits_{i=1}^2 \sqrt{\nu_i\delta_i}\int_{\mathbb{T}^2} \int_{\mathbb{R}}  |\nabla_{x'} v|^2\, dx_1 \, dx'\\
&+C\sum\limits_{i=1}^2 \sqrt{\nu_i\delta_i} \int_{\mathbb{T}^2} \int_{\mathbb{R}} |\nabla_x \times (\bold{u}-\widetilde{\bold{u}})|^2\, dx_1 \, dx' \\
\le&-\frac{3\mu}{4}\int_{\mathbb{T}^2}\int_{\mathbb{R}} a|\nabla_x \times (\bold{u}-\widetilde{\bold{u}})|^2 \, dx_1 \, dx'+C\sqrt{\delta_0}(\mathcal{G}_3+\mathcal{D}).       
\end{align*}
Similarly, we have
\begin{align*}
\mathcal{B}_{8,2,2}&=-\mu(2\mu+\lambda)\int_{\mathbb{T}^2} \int_{\mathbb{R}} (v-\widetilde{v}) \nabla_x a \cdot \nabla_x \times \nabla_x \times \bold{u} \, dx_1 \, dx' \\
&=-\mu(2\mu+\lambda)\sum\limits_{i=1}^2\int_{\mathbb{T}^2} \int_{\mathbb{R}} (v-\widetilde{v}) (a_i)_{x_1}^{X_i} (\partial_{x_2}(\partial_{x_1}u_2-\partial_{x_2}u_1)-\partial_{x_3}(\partial_{x_3}u_1-\partial_{x_1}u_3)) \, dx_1 \, dx'\\
&=\mu(2\mu+\lambda)\sum\limits_{i=1}^2\int_{\mathbb{T}^2} \int_{\mathbb{R}} (a_i)_{x_1}^{X_i} \partial_{x_2}(v-\widetilde{v})(\partial_{x_1}u_2-\partial_{x_2}u_1)\, dx_1 \, dx' \\
&\quad -\mu(2\mu+\lambda)\sum\limits_{i=1}^2\int_{\mathbb{T}^2} \int_{\mathbb{R}} (a_i)_{x_1}^{X_i} \partial_{x_3}(v-\widetilde{v})(\partial_{x_3}u_1-\partial_{x_1}u_3)\, dx_1 \, dx' \\
&\le C\sum\limits_{i=1}^2 \nu_i \delta_i\int_{\mathbb{T}^2}\int_{\mathbb{R}} |\nabla_x \times (\bold{u}-\widetilde{\bold{u}})|^2 \, dx_1 \, dx' +  C\sum\limits_{i=1}^2 \nu_i \delta_i\int_{\mathbb{T}^2}\int_{\mathbb{R}} |\nabla_{x'} (v-\widetilde{v})|^2 \, dx_1 \, dx' \\
&\le \frac{\mu}{4}\int_{\mathbb{T}^2}\int_{\mathbb{R}} a|\nabla_x \times (\bold{u}-\widetilde{\bold{u}})|^2 \, dx_1 \, dx'+C\delta_0\mathcal{D}. \\
\end{align*}
Therefore, we have
\[\mathcal{B}_{8,2}\le -\frac{\mu}{2}\int_{\mathbb{T}^2}\int_{\mathbb{R}} a|\nabla_x \times (\bold{u}-\widetilde{\bold{u}})|^2 \, dx_1 \, dx'+\sqrt{\delta_0}(\mathcal{G}_3+\mathcal{D}).\]
Finally to estimate $\mathcal{B}_3$, note that 
\[\mathcal{B}_{8,3}=-\int_{\mathbb{T}^2} \int_{\mathbb{R}} a(h_1-\widetilde{h})\partial_{x_1}(p(\widetilde{v})-p(\widetilde{v}_1^{X_1})-p(\widetilde{v}_2^{X_2})) \, dx_1 \, dx'.\]
Since $\partial_{x_1}(p(\widetilde{v})-p(\widetilde{v}_1^{X_1})-p(\widetilde{v}_2^{X_2}))=(p'(\widetilde{v})-p'(\widetilde{v}_1^{X_1}))(\widetilde{v}_1)_{x_1}^{X_1}+(p'(\widetilde{v})-p'(\widetilde{v}_2^{X_2}))(\widetilde{v}_2)_{x_1}^{X_2}$, we use Lemma \ref{shock interaction lemma1} to have
\begin{align*}
\mathcal{B}_{8,3}&\le C\varepsilon_1 \sum\limits_{i=1}^2 \int_{\mathbb{T}^2} 	\int_{\mathbb{R}} |\widetilde{v}-\widetilde{v}_i^{X_i}| |(\widetilde{v}_i)_{x_1}^{X_i}| \, dx_1 \, dx' \\
&\le C\varepsilon_1 \delta_1 \delta_2 e^{-C\min\left\{ \delta_1,\delta_2 \right\}t}.
\end{align*}
Combining the estimates above, we have
\begin{align*}
\mathcal{B}_8&\le \frac{1}{40}(C_1\mathcal{G}^S+\mathcal{G}_1+\mathcal{G}_3+\mathcal{D})+C(\delta_0+\varepsilon_1)||\nabla_x(\bold{u}-\widetilde{\bold{u}})||^2_{H^1} \\ 
& \quad +C\sum\limits_{i=1}^2\delta_i^2 e^{-C\delta_it}\int_{\mathbb{T}^2} \int_{\mathbb{R}} \eta(U|\widetilde{U}) \, dx_1 \, dx'+C\varepsilon_1 \delta_1\delta_2 e^{-C\min\left\{ \delta_1,\delta_2 \right\}t}.
\end{align*}
In conclusion, we summarize that \eqref{est-2} implies
\begin{align}
\begin{aligned} \label{4.36con}
	\frac{d}{dt}&\int_{\mathbb{T}^2}\int_{\R}a\rho\eta(U|\tU)\,dx_1 \,dx'+\sum_{i=1}^2\frac{\delta_i}{4M}|\dot{X}_i|^2+\frac{1}{2}\mathcal{G}_1+\frac{1}{2}\mathcal{G}_3 +\frac{C_1}{2}\mathcal{G}^S+\frac{1}{10}\mathcal{D}\\
	&\le C\left(\sum_{i=1}^2(\delta_i \exp(-C\delta_i t)+ \varepsilon_1\nu_i\delta_i\exp(-C\delta_it)+\delta_1\delta_2\exp(-C\min(\delta_1,\delta_2)t)+\frac{1}{t^2}\right)\int_{\mathbb{T}^2}\int_{\R}a\rho\eta(U|\tU)\,dx_1 \,dx'\\
	&\quad +C(\delta_0+\varepsilon_1)\delta_1\delta_2\exp(-C\min\left(\delta_1, \delta_2\right)t)+C(\delta_0+\varepsilon_1)||\nabla_x(\bold{u}-\widetilde{\bold{u}})||^2_{H^1}.
\end{aligned}
\end{align}

\subsection{Estimate in small time} \label{sec:sest}
Notice that the estimate  \eqref{est-I1} on $\mathcal{R}_1$ has the coefficient $\frac{1}{t^2}$, which is not integrable near $t=0$. Hence, to have the desired result, we would find a rough estimate for a short time $t\le 1$ and then we return to the preceding right-hand side $\mathcal{R}$ in \eqref{est}:
\begin{align*}
	\mathcal{R}=-\sum_{i=1}^2\frac{\delta_i}{M}|\dot{X}_i|^2 + \sum_{i=1}^2\left(\dot{X_i}\sum_{j=3}^6Y_{ij}\right) +\sum_{i=1}^8\mathcal{B}_i -\mathcal{G}_1-\mathcal{G}_2-\mathcal{G}_3-\mathcal{D}.
\end{align*}
Using Young's inequality \eqref{2young}, we first get
\begin{align*}
	\mathcal{R}+\sum_{i=1}^2\frac{\delta_i}{4M}|\dot{X}_i|^2 + \mathcal{G}_1 +  \mathcal{G}_3+ \mathcal{D} + \mathcal{G}^S \le  \sum_{i=1}^2\left(\frac{C}{\delta_i}\sum_{j=3}^6|Y_{ij}|^2\right)+\sum_{i=1}^8 \mathcal{B}_i+ \mathcal{G}^S.
\end{align*}
In addition, by \eqref{perturbation_small} and Lemma \ref{lem:shock-est} with \eqref{inta}, we get
\begin{align*}
\sum_{j=3}^6|Y_{ij}| &\le  C  \norm{(\tv_i)_{x_1}^{X_i}}_{L^2}  (\|h-\widetilde{h}\|_{L^2}+\norm{\pv-p(\tv)}_{L^2}) \\
& \quad \ \ +  \norm{(a_i)_{x_1}^{X_i}}_{L^\infty}  (\|h-\widetilde{h}\|_{L^2}^2+\norm{\pv-p(\tv)}_{L^2}^2) \\
&\le C \delta_i \varepsilon_1,
\end{align*}
which yields
\[
 \sum_{i=1}^2\left(\frac{C}{\delta_i}\sum_{j=3}^6|Y_{ij}|^2\right) \le C\varepsilon_1^2  \sum_{i=1}^2 \delta_i.
\]
Similarly, we get
\[
\sum_{i=1}^8 \mathcal{B}_i \le  C\sum_{i=1}^2 \Big(\norm{(a_i)_{x_1}^{X_i}}_{L^\infty} +  \norm{(\tv_i)_{x_1}^{X_i}}_{L^\infty} \Big) \left(\norm{v-\tv}_{H^1}^2+\norm{\bold{u}-\widetilde{\bold{u}}}_{H^1}^2\right) \le C \varepsilon_1^2  \sum_{i=1}^2\delta_i,
\]	
and
\[
 \mathcal{G}^S \le C\sum_{i=1}^2 \norm{(\tv_i)_{x_1}^{X_i}}_{L^\infty} \norm{\pv-p(\tv)}_{L^2}^2 \le C\varepsilon_1^2\sum_{i=1}^2 \delta_i^2.
\]
Thus, the estimates above provide a rough bound: for any $\delta_1, \delta_2\in (0,\delta_0)$,
\beq\label{sest}
\mathcal{R}+\sum_{i=1}^2\frac{\delta_i}{4M}|\dot{X}_i|^2 + \mathcal{G}_1 + \mathcal{G}_3+\mathcal{D} + \mathcal{G}^S \le C\delta_0,\quad t>0.
\eeq

\subsection{Proof of Lemma 4.1} We here finish the proof of Lemma 4.1. First of all, from \eqref{sest} with \eqref{est}, we have a rough estimate for $t\le 1$ as follows: 
\begin{align*}
\frac{d}{dt}\int_{\mathbb{T}^2}\int_{\R}a\rho\eta(U|\tU)\,dx_1\,dx'+\sum_{i=1}^2\frac{\delta_i}{4M}|\dot{X}_i|^2 + \mathcal{G}_1 + \mathcal{G}_3+\mathcal{D} + \mathcal{G}^S\le C\delta_0,
\end{align*}
which yields
\begin{align}
 \begin{aligned} \label{est-short-time}
	&\int_{\mathbb{T}^2}\int_{\R} a\rho\eta(U|\tU)\,d x_1\,dx'\Bigg|_{t=1} +\int_0^1\bigg(\sum_{i=1}^2\frac{\delta_i}{4M}|\dot{X}_i|^2 + \mathcal{G}_1 + \mathcal{G}_3+\mathcal{D} + \mathcal{G}^S \bigg) dt \\
 \le& \int_{\mathbb{T}^2}\int_{\R} a\eta(U|\tU)\,d x_1 \, dx'\Bigg|_{t=0}+C\delta_0.
 \end{aligned}   
\end{align}
On the other hand, for $t\ge 1$, we apply Gr\"onwall inequality to \eqref{4.36con} to have that
\begin{align}
	\begin{aligned}\label{est-long-time}
	\int_{\mathbb{T}^2}\int_{\R}a\rho&\eta(U(t,x)|\tU(t,x))\,dx_1\,dx'+\int_1^t\left(\sum_{i=1}^2\frac{\delta_i}{4M}|\dot{X}_i|^2+\frac{1}{2}\mathcal{G}_1 +\frac{1}{2}\mathcal{G}_3+\frac{C_1}{2}\mathcal{G}^S+\frac{1}{8}\mathcal{D}\right)\,ds\\
	&\le C\int_{\mathbb{T}^2}\int_{\R}a\rho\eta(U|\tU)\,dx_1 \,dx'\Bigg|_{t=1}+C \delta_0+C(\delta_0+\varepsilon_1)\int_{1}^{t}||\nabla_x(\bold{u}-\widetilde{\bold{u}})||^2_{H^1}\, ds.
	\end{aligned}
\end{align}
In the end, combining the estimates \eqref{est-short-time} and \eqref{est-long-time}, we conclude that
\begin{align*}
	\int_{\mathbb{T}^2}\int_{\R}a\rho&\eta(U(t,x)|\tU(t,x))\,dx_1 \,dx'+\int_0^t \left(\sum_{i=1}^2 \delta_i|\dot{X}_i|^2 + \mathcal{G}_1+\mathcal{G}_3+\mathcal{G}^S+\mathcal{D}\right)\,ds\\
	&\le C\int_{\mathbb{T}^2}\int_{\R}a(0,x)\rho(0,x)\eta(U_0(x)|\tU(0,x))\,dx_1\,dx' +C(\delta_0+\varepsilon_1)\int_{1}^{t}||\nabla_x(\bold{u}-\widetilde{\bold{u}})||^2_{H^1}\, ds+C\delta_0 ,
\end{align*}
which together with $\frac{1}{2}\le a\le \frac{3}{2}$, $\mathbf{G}_1(U) \sim \mathcal{G}_1(U)$, $\bold{G}_3(U) \sim \mathcal{G}_3(U)$, and $\mathbf{D}(U) \sim \mathcal{D}(U)$, completes the proof of Lemma \ref{lemma4.1}.

\section{Proof of Proposition \ref{prop:main}}\label{sec:5}
\setcounter{equation}{0}

In this section, we present the higher order estimates as in the following lemmas. Then the following lemmas and Lemma \ref{lemma4.1} complete the proof of Proposition \ref{prop:main}. The proofs for the higher order estimates follows the typical energy method. So, we only present the statements for the lemmas here, and postpone those proofs in Appendix.

\begin{lemma} \label{lem:1st}
Under the assumptions of Proposition 3.2, there exists positive constant $C$, which is independent of $\nu_i$, $\delta_i$, $\varepsilon_1$, and $T$ such that for all $t\in[0,T]$, it holds
\begin{equation}
\begin{aligned} \label{lemma5.1}
&\lVert v-\widetilde{v} \rVert_{H^1}^2+\lVert \bold{u}-\widetilde{\bold{u}} \rVert_{L^2}^2+\int_{0}^{t} \left( \sum\limits_{i=1}^2\delta_i|\dot{X}_i(\tau)|^2+\mathbf{G}_1(\tau)+\mathbf{G}_3(\tau)+\mathcal{G}^S(\tau)+\mathbf{D}(\tau) \right) \, d\tau+\int_{0}^t \lVert \nabla_x(\bold{u}-\widetilde{\bold{u}})(\tau) \rVert_{L^2}^2 \, d\tau \\
&\le C	\left( \lVert v_0-\widetilde{v}(0,\cdot) \rVert_{H^1}^2 + \lVert \bold{u}_0-\widetilde{\bold{u}}(0,\cdot) \rVert_{L^2}^2 \right)+C(\delta_0+\varepsilon_1)\int_{0}^t \lVert \nabla_x^2(\bold{u}-\widetilde{\bold{u}})(\tau) \rVert_{L^2}^2  \, d\tau+C\delta_0^{\frac{1}{2}}
\end{aligned}    
\end{equation}
\end{lemma}

\begin{lemma} \label{lem:2nd}
Under the assumptions of Proposition \ref{prop:main}, there exists positive constant $C$, which is independent of $\nu_i$, $\delta_i$, $\varepsilon_1$, 
and $T$, such that for all $t\in[0,T]$, it holds
\begin{align} 
\lVert \nabla_x(\bold{u}-\widetilde{\bold{u}}) \rVert_{L^2}^2+\int_{0}^t\lVert \nabla_x^2(u-\widetilde{u})(\tau) \rVert_{L^2}^2 \, d\tau \le C\left( \lVert v_0-\widetilde{v}(0,\cdot) \rVert_{H^1}^2 + \lVert \bold{u}_0-\widetilde{\bold{u}}(0,\cdot) \rVert_{H^1}^2 \right)+C\delta_0^{\frac{1}{2}}.    
\end{align}
\end{lemma}

\begin{lemma} \label{lem:3rd}
 Under the hypotheses of Proposition \ref{prop:main}, there exists positive constant $C$, which is independent of $\nu_i$, $\delta_i$, $\varepsilon_1$, and $T$, such that for all $t\in\left[ 0, t \right]$, it holds
\begin{equation}
\begin{aligned} \label{lemma5.3eq}
&\lVert \nabla_x^2(v-\widetilde{v})(t) \rVert^2+\int_{0}^{t} \lVert \nabla_x^2(v-\widetilde{v})(\tau) \rVert^2 \, d\tau \\
\le &C\left( \lVert v_0-\widetilde{v}(0,\cdot) \rVert_{H^2}^2 + \lVert \bold{u}_0-\widetilde{\bold{u}}(0,\cdot) \rVert_{H^1}^2\right)+C(\delta_0+\varepsilon_1)\int_{0}^{t} \lVert \nabla_x^3(\bold{u}-\widetilde{\bold{u}})(\tau) \rVert^2 \, d\tau+C\delta_0^{\frac{1}{2}}.
\end{aligned}    
\end{equation}
\end{lemma}

\begin{lemma} \label{lem:4th}
Under the assumptions of Proposition \ref{prop:main}, there exists positive constant $C$, which is independent of $\nu_i$, $\delta_i$, $\varepsilon_1$,
and $T$, such that for all $t\in\left[ 0, t \right]$, it holds
\begin{equation} \label{5.17}
\lVert \nabla_x^2(\bold{u}-\widetilde{\bold{u}})(t) \rVert^2+\int_{0}^{t} \lVert \nabla_x^3(\bold{u}-\widetilde{\bold{u}})(\tau) \rVert^2 \, d\tau \le C\left( \lVert v_0-\widetilde{v}(0,\cdot) \rVert_{H^2}^2 + \lVert \bold{u}_0-\widetilde{\bold{u}}(0,\cdot) \rVert_{H^2}^2\right)+C\delta_0^{\frac{1}{2}}.   
\end{equation} 
\end{lemma}

\begin{appendix}
\setcounter{equation}{0}
\section{Proof of Lemma \ref{lem:rel-ent}}  \label{app:rel}
First, it follows from the equations that
\begin{align*}
    &\partial_t(a\rho Q(v|\widetilde{v}))=-\rho Q(v|\widetilde{v})\sum_{i=1}^2 (\sigma_i+\dot{X}_i(t))(a_i)_{x_1}^{X_i}+a\rho_t Q(v|\widetilde{v})+a\rho[-v_t(p(v)-p(\widetilde{v}))+p'(\widetilde{v})(v-\widetilde{v})(\widetilde{v})_t]\\
&\qquad \qquad \qquad \  =-\rho Q(v|\widetilde{v})\sum_{i=1}^2 \dot{X}_i(t)(a_i)_{x_1}^{X_i}+a\rho_tQ(v|\widetilde{v})-\rho Q(v|\widetilde{v})\sum\limits_{i=1}^2 \sigma_i(a_i)_{x_1}^{X_i} \\
&\qquad \qquad \qquad \quad \ -a(p(v)-p(\widetilde{v}))\rho (v-\widetilde{v})_t + \rho a p(v|\widetilde{v})\sum_{i=1}^2(\sigma_i+\dot{X}_i(t))(\widetilde{v_i})_{x_1}^{X_i},\\
& div_x(a\rho \bold{u} Q(v|\widetilde{v}))=\rho u_1 Q(v|\widetilde{v})\sum\limits_{i=1}^2 (a_i)_{x_1}^{X_i}+aQ(v|\widetilde{v})div_x(\rho \bold{u}) \\ 
&\qquad \qquad \qquad \qquad \  -a\rho u_1 p(v|\widetilde{v})\sum\limits_{i=1}^2(\widetilde{v}_i)_{x_1}^{X_i}-a(p(v)-p(\widetilde{v}))\rho \nabla_x(v-\widetilde{v})\cdot \bold{u}. \\
\end{align*}
These imply that
\begin{align}
\begin{aligned} \label{rvt equation1}
-a(p(v)-p(\widetilde{v}))\rho (v-\widetilde{v})_t &=\partial_t(a\rho Q(v|\widetilde{v}))-aQ(v|\widetilde{v})\partial_t\rho\\
&\qquad+\rho Q(v|\widetilde{v})\sum\limits_{i=1}^2( \sigma_i+\dot{X_i}(t))(a_i)_{x_1}^{X_i} \\
&\qquad- \rho a p(v|\widetilde{v})\sum_{i=1}^2(\sigma_i+\dot{X}_i(t))(\widetilde{v_i})_{x_1}^{X_i}, \\
-a(p(v)-p(\widetilde{v}))\rho \bold{u}\cdot \nabla_x(v-\widetilde{v})&=div_x(a\rho\bold{u}Q(v|\widetilde{v}))-aQ(v|\widetilde{v})div_x(\rho \bold{u})\\
&\qquad -\rho u_1Q(v|\widetilde{v})\sum\limits_{i=1}^2(a_i)_{x_1}^{X_i} \\
&\qquad + a\rho u_1p(v|\widetilde{v})\sum\limits_{i=1}^2(\widetilde{v}_i)_{x_1}^{X_i}.
\end{aligned}
\end{align}
Thus, multiplying \eqref{mixed equations}$_1$ by $-a(p(v)-p(\widetilde{v}))$ and using the definition of $\sigma^{\ast}_i$, it holds
\begin{equation}
\begin{aligned} \label{before diffusion transform}
&\partial_t(a\rho Q(v|\widetilde{v}))+div_x(a\rho \bold{u} Q(v|\widetilde{v}))+a(p(v)-p(\widetilde{v}))div_x(\bold{h}-\widetilde{\bold{h}})\\
=&-\rho Q(v|\widetilde{v})\sum\limits_{i=1}^2\dot{X}_i(t)(a_i)_{x_1}^{X_i}-a\rho p'(\widetilde{v})(v-\widetilde{v})\sum\limits_{i=1}^2\dot{X}_i(t)(\widetilde{v}_i)_{x_1}^{X_i} \\ 
&+Q(v|\widetilde{v})\sum\limits_{i=1}^2 F_i (a_i)_{x_1}^{X_i}-Q(v|\widetilde{v})\sum\limits_{i=1}^2 \sigma_i^* (a_i)_{x_1}^{X_i} \\ 
&-ap(v|\widetilde{v})\sum\limits_{i=1}^2 F_i (\widetilde{v}_i)_{x_1}^{X_i}+ap(v|\widetilde{v})\sum\limits_{i=1}^2 \sigma_i^* (\widetilde{v}_i)_{x_1}^{X_i}\\ 
&+a(p(v)-p(\widetilde{v}))\sum\limits_{i=1}^2 F_i(\widetilde{v}_i)_{x_1}^{X_i}\\ 
&-(2\mu+\lambda)a(p(v)-p(\widetilde{v}))\Delta_x(v-\widetilde{v}).
\end{aligned}
\end{equation}
Note that 
\begin{equation*}
\nabla_xv=\frac{\nabla_xp(v)}{-\gamma p(v)^{1+\frac{1}{\gamma}}} \ \ and \ \ \nabla_x\widetilde{v}=\frac{\nabla_xp(\widetilde{v})}{-\gamma p(\widetilde{v})^{1+\frac{1}{\gamma}}}.
\end{equation*}
To get the desired diffusion term, we manipulate the last term in \eqref{before diffusion transform} as
\begin{align*}
&-(2\mu+\lambda)a(p(v)-p(\widetilde{v}))\Delta_x(v-\widetilde{v})\\
=&-(2\mu+\lambda)div_x(a(p(v)-p(\widetilde{v}))\nabla_x(v-\widetilde{v}))\\
\ &+(2\mu+\lambda)\nabla_x(a(p(v)-p(\widetilde{v})))\cdot \nabla_x(v-\widetilde{v})\\
=&-(2\mu+\lambda)div_x(a(p(v)-p(\widetilde{v}))\nabla_x(v-\widetilde{v}))\\
\ &+(2\mu+\lambda)a\nabla_x(p(v)-p(\widetilde{v}))\cdot \left(\frac{\nabla_xp(v)}{-\gamma p(v)^{1+\frac{1}{\gamma}}}-\frac{\nabla_xp(\widetilde{v})}{-\gamma p(\widetilde{v})^{1+\frac{1}{\gamma}}}\right)\\
\ &+(2\mu+\lambda)\sum\limits_{i=1}^2(a_i)_{x_1}^{X_i}(p(v)-p(\widetilde{v}))\left(\frac{\partial_{x_1}(p(v))}{-\gamma p(v)^{1+\frac{1}{\gamma}}}-\frac{\partial_{x_1}p(\widetilde{v})}{-\gamma p(\widetilde{v})^{1+\frac{1}{\gamma}}}\right)\\
=&-(2\mu+\lambda)div_x(a(p(v)-p(\widetilde{v}))\nabla_x(v-\widetilde{v}))\\
&-(2\mu+\lambda)a\frac{|\nabla_x(p(v)-p(\widetilde{v}))|^2}{\gamma p(v)^{1+\frac{1}{\gamma}}}\\ 
&-(2\mu+\lambda)a\partial_{x_1}(p(v)-p(\widetilde{v})) \partial_{x_1}p(\widetilde{v})\left( \frac{1}{\gamma p(v)^{1+\frac{1}{\gamma}}}- \frac{1}{\gamma p(\widetilde{v})^{{1+\frac{1}{\gamma}}}}\right)\\
&-(2\mu+\lambda)(p(v)-p(\widetilde{v}))\frac{\partial_{x_1}(p(v)-p(\widetilde{v}))}{\gamma p(v)^{1+\frac{1}{\gamma}}}\sum\limits_{i=1}^2 (a_i)_{x_1}^{X_i}\\ 
&-(2\mu+\lambda)(p(v)-p(\widetilde{v}))\partial_{x_1}p(\widetilde{v})\left( \frac{1}{\gamma p(v)^{1+\frac{1}{\gamma}}}- \frac{1}{\gamma p(\widetilde{v})^{{1+\frac{1}{\gamma}}}}\right) \sum\limits_{i=1}^2 (a_i)_{x_1}^{X_i}.  
\end{align*}
On the other hand, note that
\begin{align*}
&\bullet \ \partial_t\left( a\rho \frac{|\bold{h}-\widetilde{\bold{h}}|^2}{2} \right)=-\rho \frac{|\bold{h}-\widetilde{\bold{h}}|^2}{2}\sum\limits_{i=1}^2(\sigma_i+\dot{X}_i(t))(a_i)_{x_1}^{X_i}+a\rho_t \frac{|\bold{h}-\widetilde{\bold{h}}|^2}{2}+a\rho (\bold{h}-\widetilde{\bold{h}})(\bold{h}-\widetilde{\bold{h}})_t \\
&\bullet \ div_x\left( a\rho \bold{u} \frac{|\bold{h}-\widetilde{\bold{h}}|^2}{2} \right)=\rho u_1 \frac{|\bold{h}-\widetilde{\bold{h}}|^2}{2}\sum\limits_{i=1}^2(a_i)_{x_1}^{X_i}+a\frac{|\bold{h}-\widetilde{\bold{h}}|^2}{2}div_x(\rho \bold{u})\\ 
& \qquad \qquad \qquad \qquad \qquad \ \ +a\rho u_1 (\bold{h}-\widetilde{\bold{h}}) \cdot \partial_{x_1}(\bold{h}-\widetilde{\bold{h}})+a\rho u_2 (\bold{h}-\widetilde{\bold{h}}) \cdot \partial_{x_2}(\bold{h}-\widetilde{\bold{h}})+a\rho u_3 (\bold{h}-\widetilde{\bold{h}}) \cdot \partial_{x_3}(\bold{h}-\widetilde{\bold{h}}) \\ 
& \qquad \qquad \qquad \qquad \quad \ \ =\rho u_1 \frac{|\bold{h}-\widetilde{\bold{h}}|^2}{2}\sum\limits_{i=1}^2(a_i)_{x_1}^{X_i}+adiv_x(\rho \bold{u})\frac{|\bold{h}-\widetilde{\bold{h}}|^2}{2}+a(\bold{h}-\widetilde{\bold{h}})\cdot \rho\bold{u} \nabla_x(\bold{h}-\widetilde{\bold{h}})\\ 
&\bullet \ div_x(a(p(v)-p(\widetilde{v}))(\bold{h}-\widetilde{\bold{h}}))-a(p(v)-p(\widetilde{v}))div_x(\bold{h}-\widetilde{\bold{h}})\\
&\quad =\partial_{x_1}(a(p(v)-p(\widetilde{v})))(h_1-\widetilde{h})+\partial_{x_2}(a(p(v)-p(\widetilde{v})))h_2+\partial_{x_3}(a(p(v)-p(\widetilde{v})))h_3 \\ 
& \quad =(p(v)-p(\widetilde{v}))(h_1-\widetilde{h})\sum\limits_{i=1}^2(a_i)_{x_1}^{X_i}+a(\bold{h}-\widetilde{\bold{h}}) \cdot \nabla_x(p(v)-p(\widetilde{v})).
\end{align*}

\noindent Multiplying \eqref{mixed equations}$_2$ by $a(\bold{h}-\bold{\widetilde{h}})$ and using these equalities above, we have
\begin{equation}
\begin{aligned} \label{4.4.2 change equation}
&\partial_t\left( a\rho \frac{|\bold{h}-\widetilde{\bold{h}}|^2}{2} \right)+div_x\left( a\rho \bold{u} \frac{|\bold{h}-\widetilde{\bold{h}}|^2}{2} \right)\\
+&div_x(a(p(v)-p(\widetilde{v}))(\bold{h}-\widetilde{\bold{h}}))-a(p(v)-p(\widetilde{v}))div_x(\bold{h}-\widetilde{\bold{h}})\\
&=-\rho \frac{|\bold{h}-\widetilde{\bold{h}}|^2}{2}\sum\limits_{i=1}^2 \dot{X}_i(t)(a_i)_{x_1}^{X_i}+\rho \frac{|\bold{h}-\widetilde{\bold{h}}|^2}{2}\sum\limits_{i=1}^2F_i(a_i)_{x_1}^{X_i}-\rho \frac{|\bold{h}-\widetilde{\bold{h}}|^2}{2}\sum\limits_{i=1}^2\sigma_i^*(a_i)_{x_1}^{X_i}\\ 
&\qquad +a\rho (h_1-\widetilde{h})\sum\limits_{i=1}^2\dot{X}_i(t) (\widetilde{h}_{i})_{x_1}^{X_i}-a(h_1-\widetilde{h})\sum\limits_{i=1}^2F_i (\widetilde{h}_{i})_{x_1}^{X_i}+a(\bold{h}-\widetilde{\bold{h}})\cdot R^* \\ 
&\qquad+\rho u_1 \frac{|\bold{h}-\widetilde{\bold{h}}|^2}{2}\sum\limits_{i=1}^2 (a_i)_{x_1}^{X_i}+(p(v)-p(\widetilde{v}))(h_1-\widetilde{h})\sum\limits_{i=1}^2 (a_i)_{x_1}^{X_i}. \\
\end{aligned}
\end{equation}
Adding \eqref{before diffusion transform} and \eqref{4.4.2 change equation} together and then integrating the result by parts over $\Omega$, we get 
\begin{align}
\begin{aligned}
&\frac{d}{dt}\int_{\mathbb{T}^2} \int_{\mathbb{R}} a\rho\left( Q(v|\widetilde{v})+\frac{|\bold{h}-\widetilde{\bold{h}}|^2}{2} \right) \, dx_1 \, dx' \\ 
&=\sum\limits_{i=1}^2 \int_{\mathbb{T}^2} \int_{\mathbb{R}} \dot{X}_i(t)\rho\left(-Q(v|\widetilde{v})(a_i)_{x_1}^{X_i} -\frac{|\bold{h}-\widetilde{\bold{h}}|^2}{2}(a_i)_{x_1}^{X_i}-ap'(\widetilde{v})(v-\widetilde{v})(\widetilde{v}_i)_{x_1}^{X_i}+a(h_1-\widetilde{h})(\widetilde{h}_{i})_{x_1}^{X_i}   \right) \, dx_1 \, dx' \\ 
&\quad +\sum\limits_{i=1}^2 \int_{\mathbb{T}^2} \int_{\mathbb{R}} (a_i)_{x_1}^{X_i}(p(v)-p(\widetilde{v}))(h_1-\widetilde{h}) \, dx_1 \, dx'+\sum\limits_{i=1}^2 \int_{\mathbb{T}^2} \int_{\mathbb{R}} \sigma_i^*ap(v|\widetilde{v})  (\widetilde{v}_i)_{x_1}^{X_i} \, dx_1 \, dx' \\
& \quad +\sum\limits_{i=1}^2\int_{\mathbb{T}^2} \int_{\mathbb{R}}ap'(\widetilde{v})(v-\widetilde{v})F_i (\widetilde{v}_i)_{x_1}^{X_i} \, dx_1 \, dx' - \sum\limits_{i=1}^2 \int_{\mathbb{T}^2} \int_{\mathbb{R}} a(h_1-\widetilde{h})F_i(\widetilde{h}_{i})_{x_1}^{X_i} \, dx_1 \, dx' \\ 
& \quad +\sum\limits_{i=1}^2 \int_{\mathbb{T}^2} \int_{\mathbb{R}} (Q(v|\widetilde{v})+\frac{|\bold{h}-\widetilde{\bold{h}}|^2}{2}) F_i (a_i)_{x_1}^{X_i} \, dx_1 \, dx' \\ 
& \quad -(2\mu+\lambda)\int_{\mathbb{T}^2} \int_{\mathbb{R}} a\partial_{x_1}(p(v)-p(\widetilde{v})) \partial_{x_1}p(\widetilde{v})\left( \frac{1}{\gamma p(v)^{1+\frac{1}{\gamma}}}- \frac{1}{\gamma p(\widetilde{v})^{{1+\frac{1}{\gamma}}}}\right)\\
&\quad -(2\mu+\lambda)\sum\limits_{i=1}^2 \int_{\mathbb{T}^2} \int_{\mathbb{R}} (a_i)_{x_1}^{X_i}(p(v)-p(\widetilde{v}))\frac{\partial_{x_1}(p(v)-p(\widetilde{v}))}{\gamma p(v)^{1+\frac{1}{\gamma}}} \, dx_1 \, dx' \\
&\quad -(2\mu+\lambda)\sum\limits_{i=1}^2\int_{\mathbb{T}^2} \int_{\mathbb{R}} (a_i)_{x_1}^{X_i}(p(v)-p(\widetilde{v}))\partial_{x_1}p(\widetilde{v})\left( \frac{1}{\gamma p(v)^{1+\frac{1}{\gamma}}}- \frac{1}{\gamma p(\widetilde{v})^{{1+\frac{1}{\gamma}}}}\right) \, dx_1 \, dx'  \\ 
&\quad +\int_{\mathbb{T}^2} \int_{\mathbb{R}} a(\bold{h}-\widetilde{\bold{h}})\cdot R^* \, dx_1 \, dx'- \sum\limits_{i=1}^2 \int_{\mathbb{T}^2} \int_{\mathbb{R}} \sigma_i^*Q(v|\widetilde{v}) (a_i)_{x_1}^{X_i} \, dx_1 \, dx' \\
& \quad -(2\mu+\lambda)\int_{\mathbb{T}^2} \int_{\mathbb{R} }a\frac{|\nabla_x(p(v)-p(\widetilde{v}))|^2}{\gamma p(v)^{1+\frac{1}{\gamma}}} \, dx_1 \, dx', \\
\end{aligned}
\end{align}
which complete the proof of Lemma 4.2.

\section{Proof of Lemma \ref{lem:1st}}
Recall 
\begin{equation}    
\begin{cases} \label{lemma5.1eq}
\rho(v-\widetilde{v})_t+\rho \bold{u}\cdot \nabla_x(v-\widetilde{v})-\rho \sum\limits_{i=1}^2\dot{X}_i(t)(\widetilde{v}_i)_{x_1}^{X_i}+\sum\limits_{i=1}^2F_i(\widetilde{v}_i)_{x_1}^{X_i}=div_x(\bold{u}-\widetilde{\bold{u}}), \\
\quad \rho(\bold{u}-\widetilde{\bold{u}})_t+\rho \bold{u} \nabla_x(\bold{u}-\widetilde{\bold{u}})+\nabla_x(p(v)-p(\widetilde{v}))-\rho\sum\limits_{i=1}^2\dot{X}_i(t)(\widetilde{\bold{u}}_i)_{x_1}+\sum\limits_{i=1}^2 F_i (\widetilde{\bold{u}}_i)_{x_1}\\
=\mu\Delta_x(\bold{u}-\widetilde{\bold{u}})+(\mu+\lambda)\nabla_x div_x(\bold{u}-\widetilde{\bold{u}})-\nabla_x(p(\widetilde{v})-p(\widetilde{v}_1^{X_1})-p(\widetilde{v}_2^{X_2}) ).
\end{cases}
\end{equation}
Multiplying \eqref{lemma5.1eq}$_1$ by $-(p(v)-p(\widetilde{v}))$ and using the similar way as in Section 4.1, we get
\begin{equation*}
\begin{aligned}
&\partial_t(\rho Q(v|\widetilde{v}))+div_x(\rho \bold{u} Q(v|\widetilde{v}))+(p(v)-p(\widetilde{v}))div_x(\bold{u}-\widetilde{\bold{u}}) \\
=&-\rho \sum\limits_{i=1}^2 \rho p'(\widetilde{v})(v-\widetilde{v})\dot{X}_i(\widetilde{v}_i)_{x_1}^{X_i}+\sum\limits_{i=1}^2 p(v|\widetilde{v})\sigma_i^* (\widetilde{v}_i)_{x_1}^{X_i}+\sum\limits_{i=1}^2p'(\widetilde{v})(v-\widetilde{v})F_i(\widetilde{v}_i)_{x_1}^{X_i}.   
\end{aligned}
\end{equation*}
Multiplying \eqref{lemma5.1eq}$_2$ by $\bold{u}-\widetilde{\bold{u}}$ and using the similar way as in Section 4.1, we get 
\begin{equation*}
\begin{aligned}
&\partial_t\left( \rho \frac{|\bold{u}-\widetilde{\bold{u}}|^2}{2} \right)+div_x\left( \rho \bold{u} \frac{|\bold{u}-\widetilde{\bold{u}}|^2}{2} \right)+div_x\left[ (p(v)-p(\widetilde{v}))(\bold{u}-\widetilde{\bold{u}}) \right]-(p(v)-p(\widetilde{v}))div_x(\bold{u}-\widetilde{\bold{u}})\\
=&\rho(u_1-\widetilde{u})\sum\limits_{i=1}^2 \dot{X}_i(\widetilde{u}_{i})_{x_1}^{X_i}-(u_1-\widetilde{u})\sum\limits_{i=1}^2F_i(\widetilde{u}_{i})_{x_1}^{X_i}-\mu \left| \nabla_x(\bold{u}-\widetilde{\bold{u}}) \right|^2-(\mu+\lambda)\left| div_x(\bold{u}-\widetilde{\bold{u}}) \right|^2\\
&+\mu div_x\left[ \nabla_x(p(v)-p(\widetilde{v}))\cdot(\bold{u}-\widetilde{\bold{u}}) \right]+(\mu+\lambda)div\left[ div_x(\bold{u}-\widetilde{\bold{u}})(\bold{u}-\widetilde{\bold{u}}) \right]\\
&-\nabla_x( p(\widetilde{v})-p(\widetilde{v}_1^{X_1})_{x_1}-p(\widetilde{v}_2^{X_2})_{x_1} )\cdot(\bold{u}-\widetilde{\bold{u}}).\\
\end{aligned}    
\end{equation*}
Adding two equalities above and integrating that over $\Omega$, we get
\begin{equation}
\begin{aligned}
&\frac{d}{dt}\int_{\mathbb{T}^2} \int_{\mathbb{R}} \rho\left( Q(v|\widetilde{v})+\frac{\left| \bold{u}-\widetilde{\bold{u}} \right|^2}{2} \right)  \, dx_1\, dx'\\
+&\int_{\mathbb{T}^2} \int_{\mathbb{R}} \mu\left| \nabla_x(\bold{u}-\widetilde{\bold{u}}) \right|^2+(\mu+\lambda)\left| div_x(\bold{u}-\widetilde{\bold{u}}) \right|^2  \, dx_1\, dx'\\
&=\sum\limits_{i=1}^2 \dot{X}_i(t)Z_i(t)+\sum\limits_{i=1}^4I_i(t),
\end{aligned}
\end{equation}
where
\begin{equation*}
\begin{aligned}
Z_i(t)&=-\int_{\mathbb{T}^2}\int_{\mathbb{R}} \rho p'(\widetilde{v})(v-\widetilde{v})(\widetilde{v}_i)_{x_1}^{X_i}  \, dx_1 \, dx'+\int_{\mathbb{T}^2}\int_{\mathbb{R}} \rho (u_1-\widetilde{u})(\widetilde{u}_{i})_{x_1}^{X_i}  \, dx_1 \, dx'\\
&=:Z_{1,i}(t)+Z_{2,i}(t),\\
I_{1,i}(t)&=\sigma_i^* \int_{\mathbb{T}^2} \int_{\mathbb{R}} p(v|\widetilde{v})(\widetilde{v}_i)_{x_1}^{X_i}  \, dx_1 \, dx',\\ 
I_{2,i}(t)&=\sigma_i^* \int_{\mathbb{T}^2} \int_{\mathbb{R}} F_ip'(\widetilde{v})(v-\widetilde{v})(\widetilde{v}_i)_{x_1}^{X_i}  \, dx_1 \, dx',\\
I_{3,i}(t)&=-\int_{\mathbb{T}^2} \int_{\mathbb{R}} F_i(u_1-\widetilde{u})(\widetilde{u}_{i})_{x_1}^{X_i}  \, dx_1 \, dx', \\
I_4(t)&=-\int_{\mathbb{T}^2} \int_{\mathbb{R}} \nabla_x( p(\widetilde{v})-p(\widetilde{v}_1^{X_1})-p(\widetilde{v}_2^{X_2}) )\cdot(\bold{u}-\widetilde{\bold{u}})  \, dx_1 \, dx'.
\end{aligned}
\end{equation*}
First, by \eqref{specialterm} we get the estimates of $Z_i(t)$.
\begin{equation*}
\begin{aligned}
\left| Z_{1,i}(t) \right|&\le C \sqrt{\int_{\mathbb{T}^2} \int_{\mathbb{R}}	|(\widetilde{v}_i)_{x_1}^{X_i}|  \, dx_1 \, dx'} \sqrt{\int_{\mathbb{T}^2} \int_{\mathbb{R}}	\left| v-\widetilde{v} \right|^2  |(\widetilde{v}_i)_{x_1}^{X_i}|  \, dx_1 \, dx'} \\
&\le C\sqrt{\delta_i}	\sqrt{\left( \mathcal{G}^S+\sum\limits_{i=1}^2 \delta_i^2e^{-C\delta_it} \int_{\mathbb{T}^2}  \int_{\mathbb{R}} \eta(U|\widetilde{U}) \, dx_1 \, dx' \right)}\\
\Longrightarrow |\dot{X}_i| \left| Z_{1,i} \right| &\le \frac{\delta_i}{8}|\dot{X}_i|^2+C\mathcal{G}^S+C\sum\limits_{i=1}^2 \delta_i^2e^{-C\delta_it} \int_{\mathbb{T}^2}  \int_{\mathbb{R}} \eta(U|\widetilde{U}) \, dx_1 \, dx',
\end{aligned}    
\end{equation*}
\begin{equation*}
\begin{aligned}
\left| Z_{2,i}(t) \right|
&\le C \int_{\mathbb{T}^2} \int_{\mathbb{R}}			\Bigg(\left| h_1-\widetilde{h}-\frac{p(v)-p(\widetilde{v})}{\sigma_i^*} \right|+ \left| p(v)-p(\widetilde{v}) \right| \\
&\qquad \qquad \qquad \quad +\left| \partial_{x_1}(p(v)-p(\widetilde{v})) \right|+ 	|(\widetilde{v}_i^{X_i})_{x_1}|	\left| v-\widetilde{v} \right|\Bigg) | (\widetilde{v}_i)_{x_1}^{X_i}|  \, dx_1 \, dx'\\
&\le C\sqrt{\int_{\mathbb{T}^2} \int_{\mathbb{R}}	 | (\widetilde{v}_i)_{x_1}^{X_i} |  \, dx_1 \, dx'}  \sqrt{\int_{\mathbb{T}^2} \int_{\mathbb{R}}		\left| h_1-\widetilde{h}-\frac{p(v)-p(\widetilde{v})}{\sigma_i^*} \right|^2 |(a_i)_{x_1}^{X_i}|  \, dx_1 \, dx'}\\
& \quad +C\sqrt{\int_{\mathbb{T}^2} \int_{\mathbb{R}}	 | (\widetilde{v}_i)_{x_1}^{X_i} |  \, dx_1 \, dx'}  \sqrt{\int_{\mathbb{T}^2} \int_{\mathbb{R}}	\left|	p(v)-p(\widetilde{v})\right|^2 | (\widetilde{v}_i)_{x_1}^{X_i} |  \, dx_1 \, dx'}\\
& \quad +C\sqrt{\int_{\mathbb{T}^2} \int_{\mathbb{R}}	 | (\widetilde{v}_i)_{x_1}^{X_i} |^2  \, dx_1 \, dx'}  \sqrt{\int_{\mathbb{T}^2} \int_{\mathbb{R}}	\left|	\nabla_x(p(v)-p(\widetilde{v}))\right|^2   \, dx_1 \, dx'}\\
&\le C\sqrt{\delta_i}(\sqrt{\bold{G}_1}+\sqrt{ \mathcal{G}^S+\sum\limits_{i=1}^2 \delta_i^2e^{-C\delta_it} \int_{\mathbb{T}^2}  \int_{\mathbb{R}} \eta(U|\widetilde{U}) \, dx_1 \, dx' }+\delta_i\sqrt{\bold{D}})\\
\Longrightarrow |\dot{X}_i| \left| Z_{2,i} \right| &\le \frac{\delta_i}{8}|\dot{X}_i|^2+C\bold{G}_1+C\mathcal{G}^S+C\sum\limits_{i=1}^2 \delta_i^2e^{-C\delta_it} \int_{\mathbb{T}^2}  \int_{\mathbb{R}} \eta(U|\widetilde{U}) \, dx_1 \, dx'+C\delta_i^2\bold{D}. 
\end{aligned}    
\end{equation*}
Also, we estimate $I$-terms. 

By \eqref{specialterm},
\begin{equation*}
\begin{aligned}
|I_{1,i}|&\le C\int_{\mathbb{T}^2} \int_{\mathbb{R}}	\left| p(v)-p(\widetilde{v}) \right|^2  | (\widetilde{v}_i)_{x_1}^{X_i} |  \, dx_1 \, dx' \\ 
&\le C\mathcal{G}^S+C\sum\limits_{i=1}^2 \delta_i^2e^{-C\delta_it} \int_{\mathbb{T}^2}  \int_{\mathbb{R}} \eta(U|\widetilde{U}) \, dx_1 \, dx'. 
\end{aligned}    
\end{equation*}
As in $\mathcal{B}_4$,
\begin{equation*}
\begin{aligned}
|I_{2,i}|&\le C\int_{\mathbb{T}^2} \int_{\mathbb{R}}	|F_i| \left| p(v)-p(\widetilde{v}) \right|  | (\widetilde{v}_i)_{x_1}^{X_i} |  \, dx_1 \, dx' \\
&\le C\mathcal{G}^S+C\frac{\delta_i}{\nu_i}(\bold{G}_1+\bold{G}_3+\bold{D})+C\sum\limits_{i=1}^2\delta_i^2e^{-C\delta_it} \int_{\mathbb{T}^2}  \int_{\mathbb{R}} \eta(U|\widetilde{U}) \, dx_1 \, dx'.    
\end{aligned}    
\end{equation*}
Using $u_1-\widetilde{u}=h_1-\widetilde{h}+(2\mu+\lambda)\partial_{x_1}(v-\widetilde{v})$, as in $\mathcal{B}_4$, we get
\begin{equation*}
\begin{aligned}
|I_{3,i}|\le C\mathcal{G}^S+C\frac{\delta_i}{\nu_i}(\bold{G}_1+\bold{G}_3+\bold{D})+C\sum\limits_{i=1}^2\delta_i^2e^{-C\delta_it} \int_{\mathbb{T}^2}  \int_{\mathbb{R}} \eta(U|\widetilde{U}) \, dx_1 \, dx'.    
\end{aligned}    
\end{equation*}
Finally, to estimate $I_4$, note that 
\begin{align*}
\nabla_x( p(\widetilde{v})-p(\widetilde{v}_1^{X_1})-p(\widetilde{v}_2^{X_2}) )\cdot(\bold{u}-\widetilde{\bold{u}})=\partial_{x_1}(p(\widetilde{v})-p(\widetilde{v}_1^{X_1})-p(\widetilde{v}_2^{X_2}))(u_1-\widetilde{u}),    
\end{align*}
\begin{align*}
\partial_{x_1}(p(\widetilde{v})-p(\widetilde{v}_1^{X_1})-p(\widetilde{v}_2^{X_2}))&=p'(\widetilde{v})((\widetilde{v}_1)_{x_1}^{X_1}+(\widetilde{v}_2)_{x_1}^{X_2})-p'(\widetilde{v}_1^{X_1})(\widetilde{v}_1)_{x_1}^{X_1}-p'(\widetilde{v}_2^{X_2})(\widetilde{v}_2)_{x_1}^{X_2}\\
&=(\widetilde{v}_1)_{x_1}^{X_1}(p'(\widetilde{v})-p'(\widetilde{v}_1^{X_1}))+(\widetilde{v}_2)_{x_1}^{X_2}(p'(\widetilde{v})-p'(\widetilde{v}_2^{X_2})). 
\end{align*}
Then as in $\mathcal{B}_3$, we get
\begin{equation*}
\begin{aligned}    
|I_4|&\le C\sum\limits_{i=1}^2 \int_{\mathbb{T}^2}  \int_{\mathbb{R}} | (\widetilde{v}_i)_{x_1}^{X_i} | | \widetilde{v}-\widetilde{v}_i^{X_i} | \left(| h_1-\widetilde{h}|+\left| \partial_{x_1}(v-\widetilde{v}) \right|  \right)\, dx_1 \, dx'\\
&\le C\delta_1\delta_2e^{-C\min\left\{ \delta_1,\delta_2 \right\}t} + C\delta_0(\bold{G}_1+\bold{D}+\mathcal{G}^S)+C\sum\limits_{i=1}^2\delta_0\delta_i^2e^{-C\delta_it} \int_{\mathbb{T}^2}  \int_{\mathbb{R}} \eta(U|\widetilde{U}) \, dx_1 \, dx'.    
\end{aligned}    
\end{equation*}
Therefore,
\begin{align} 
\begin{aligned}\label{lemma5.1ieq1}
&\frac{d}{dt}\int_{\mathbb{T}^2} \int_{\mathbb{R}} \rho\left( Q(v|\widetilde{v})+\frac{\left| \bold{u}-\widetilde{\bold{u}} \right|^2}{2} \right)  \, dx_1\, dx' + \mu\int_{\mathbb{T}^2} \int_{\mathbb{R}} \left| \nabla_x(\bold{u}-\widetilde{\bold{u}}) \right|^2  \, dx_1\, dx' \\
\le &  \sum\limits_{i=1}^2\frac{\delta_i}{4}|\dot{X}_i|^2+\mathcal{C}_1\left(\bold{G}_1+\bold{G}_3+\mathcal{G}^S+\bold{D}\right)\\
& \ \ +C\delta_1\delta_2e^{-C\min\left\{ \delta_1,\delta_2 \right\}t}+C\delta_0\sum\limits_{i=1}^2\delta_ie^{-C\delta_it} \int_{\mathbb{T}^2}  \int_{\mathbb{R}} \eta(U|\widetilde{U}) dx_1 dx',
\end{aligned}
\end{align}
where $\mathcal{C}_1>0$ is some constant.

Integrate \eqref{lemma5.1ieq1} over $\left[0, t\right]$ for any $t\le T$ and then multiply the result by $\frac{1}{2\max\left\{ 1, \mathcal{C}_1 \right\}}$.

Then by the smallness of the parameters, we have
\begin{align} 
\begin{aligned}\label{lemma5.1ieq2}
&\frac{1}{2\max\left\{ 1, \mathcal{C}_1 \right\}}\int_{\mathbb{T}^2} \int_{\mathbb{R}} \rho\left( Q(v|\widetilde{v})+\frac{\left| \bold{u}-\widetilde{\bold{u}} \right|^2}{2} \right)  \, dx_1\, dx' + \frac{\mu}{2\max\left\{ 1, \mathcal{C}_1 \right\}}\int_{0}^{t}\lVert \nabla_x\left(\bold{u}-\bold{\widetilde{u}}\right) \rVert^2_{L^2(\mathbb{R})}\, d\tau \\
\le &  \frac{1}{2\max\left\{ 1, \mathcal{C}_1 \right\}}\int_{\mathbb{T}^2}\int_{\mathbb{R}}\left( Q(v_0(x) | \widetilde{v}(0,x))+\frac{|\bold{u}_0(x)-\widetilde{\bold{u}}(0,x)|^2}{2} \right) \,dx_1 \, dx' \\ 
& \ + \frac{1}{2}\int_{0}^t \left(\sum\limits_{i=1}^2{\delta_i}|\dot{X}_i|^2+\bold{G}_1+\bold{G}_3+\mathcal{G}^S+\bold{D}\right) \, d\tau+\frac{\sqrt{\max\left\{ \delta_1, \delta_2 \right\}}}{2}+\frac{\sqrt{\delta_0}\varepsilon_1}{2}.\\
\end{aligned}
\end{align}
Adding \eqref{lemma5.1ieq2} to Lemma \ref{lemma4.1}, we complete the proof of the Lemma \ref{lemma5.1}.

\section{Proof of Lemma \ref{lem:2nd}}
Multiply \eqref{lemma5.1eq}$_2$ by $-v\Delta_x(\bold{u}-\widetilde{\bold{u}})$ and then integrate that over $\Omega$.

Using the integration by parts, it holds
\begin{equation}
\begin{aligned} \label{lemma5.2eq}
&\frac{d}{dt}\int_{\mathbb{T}^2} \int_{\mathbb{R}} \frac{\left| \nabla_x(\bold{u}-\widetilde{\bold{u}}) \right|^2}{2} \, dx_1 \, dx'\\
&+\mu\int_{\mathbb{T}^2} \int_{\mathbb{R}} v\left| \Delta_x(\bold{u}-\widetilde{\bold{u}}) \right|^2 \, dx_1 \, dx+(\mu+\lambda)\int_{\mathbb{T}^2} \int_{\mathbb{R}} \nabla_x div_x(\bold{u}-\widetilde{\bold{u}}) \cdot v\Delta_x(\bold{u}-\widetilde{\bold{u}})  \, dx_1 \, dx'\\
=&-\int_{\mathbb{T}^2} \int_{\mathbb{R}} 	\left( \left( \nabla_x \bold{u} \nabla_x(\bold{u}-\widetilde{\bold{u}})   \right) : \nabla_x(\bold{u}-\widetilde{\bold{u}}) - div_x\bold{u}\frac{	\left| \nabla_x(\bold{u}-\widetilde{\bold{u}}) \right|^2}{2} \right)  \, dx_1 \, dx'\\
&-\sum\limits_{i=1}^2\dot{X}_i(t)\int_{\mathbb{T}^2} \int_{\mathbb{R}} 	(\widetilde{u}_{i})_{x_1}^{X_i}\Delta_x(u_1-\widetilde{u})  \, dx_1 \, dx'+\sum\limits_{i=1}^2\int_{\mathbb{T}^2} \int_{\mathbb{R}} 	vF_i(\widetilde{u}_{i})_{x_1}^{X_i}\Delta_x(u_1-\widetilde{u})  \, dx_1 \, dx'\\
&+\int_{\mathbb{T}^2} \int_{\mathbb{R}} v\Delta_x(\bold{u}-\widetilde{\bold{u}})\cdot \nabla_x(p(v)-p(\widetilde{v}))	 \, dx_1 \, dx'.   
\end{aligned}    
\end{equation}
Note that 
\begin{equation*}
\begin{aligned}
&  (\mu+\lambda)\int_{\mathbb{T}^2} \int_{\mathbb{R}} \nabla_x div_x(\bold{u}-\widetilde{\bold{u}}) \cdot v\Delta_x(\bold{u}-\widetilde{\bold{u}})  \, dx_1 \, dx'\\
= \ &(\mu+\lambda)\int_{\mathbb{T}^2} \int_{\mathbb{R}}  v\left| \nabla_x div_x(\bold{u}-\widetilde{\bold{u}}) \right|^2 \, dx_1 \, dx' \\
&  +(\mu+\lambda)\int_{\mathbb{T}^2} \int_{\mathbb{R}}  v\nabla_x div_x(\bold{u}-\widetilde{\bold{u}})\cdot\left[ \Delta_x(\bold{u}-\widetilde{\bold{u}})- \nabla_x div_x(\bold{u}-\widetilde{\bold{u}})\right] \, dx_1 \, dx'.
\end{aligned}    
\end{equation*}
For the above last term, by integration by parts and by the formula \[\nabla_x \cdot \left[ \Delta_x(\bold{u}-\widetilde{\bold{u}})-\nabla_x div_x(\bold{u}-\widetilde{\bold{u}})\right]=0,\]
it holds
\begin{equation*}
\begin{aligned}
&(\mu+\lambda)\int_{\mathbb{T}^2} \int_{\mathbb{R}}  v\nabla_x div_x(\bold{u}-\widetilde{\bold{u}})\left[ \Delta_x(\bold{u}-\widetilde{\bold{u}})- \nabla_x div_x(\bold{u}-\widetilde{\bold{u}})\right] \, dx_1 \, dx'\\
=&-(\mu+\lambda)\int_{\mathbb{T}^2} \int_{\mathbb{R}}  div_x(\bold{u}-\widetilde{\bold{u}}) \nabla_x \cdot (v\left[ \Delta_x(\bold{u}-\widetilde{\bold{u}})- \nabla_x div_x(\bold{u}-\widetilde{\bold{u}})\right]) \, dx_1 \, dx'\\
=&-(\mu+\lambda)\int_{\mathbb{T}^2} \int_{\mathbb{R}}  div_x(\bold{u}-\widetilde{\bold{u}}) \nabla_x v \cdot (\left[ \Delta_x(\bold{u}-\widetilde{\bold{u}})- \nabla_x div_x(\bold{u}-\widetilde{\bold{u}})\right]) \, dx_1 \, dx'. 
\end{aligned}    
\end{equation*}
Substituting this into \eqref{lemma5.2eq}, we get
\begin{equation*}
\begin{aligned}
&\frac{d}{dt}\int_{\mathbb{T}^2} \int_{\mathbb{R}} \frac{\left| \nabla_x(\bold{u}-\widetilde{\bold{u}}) \right|^2}{2} \, dx_1 \, dx'\\
&+\underbrace{\mu\int_{\mathbb{T}^2} \int_{\mathbb{R}} v\left| \Delta_x(\bold{u}-\widetilde{\bold{u}}) \right|^2 \, dx_1 \, dx+(\mu+\lambda)\int_{\mathbb{T}^2} \int_{\mathbb{R}}  v\left| \nabla_x div_x(\bold{u}-\widetilde{\bold{u}}) \right|^2 \, dx_1 \, dx'}_{:=\mathcal{D}_2}  \\
=&-\int_{\mathbb{T}^2} \int_{\mathbb{R}} 	\left( \left( \nabla_x \bold{u} \nabla_x(\bold{u}-\widetilde{\bold{u}})   \right) : \nabla_x(\bold{u}-\widetilde{\bold{u}}) - div_x\bold{u}\frac{	\left| \nabla_x(\bold{u}-\widetilde{\bold{u}}) \right|^2}{2} \right)  \, dx_1 \, dx' \\
&-\sum\limits_{i=1}^2\dot{X}_i(t)\int_{\mathbb{T}^2} \int_{\mathbb{R}} 	(\widetilde{u}_{i})_{x_1}^{X_i}\Delta_x(u_1-\widetilde{u})  \, dx_1 \, dx'+\sum\limits_{i=1}^2\int_{\mathbb{T}^2} \int_{\mathbb{R}} 	vF_i(\widetilde{u}_{i})_{x_1}^{X_i}\Delta_x(u_1-\widetilde{u})  \, dx_1 \, dx'\\
&+\int_{\mathbb{T}^2} \int_{\mathbb{R}} v\Delta_x(\bold{u}-\widetilde{\bold{u}})\cdot \nabla_x(p(v)-p(\widetilde{v}))	 \, dx_1 \, dx'\\
&+(\mu+\lambda)\int_{\mathbb{T}^2} \int_{\mathbb{R}}  div_x(\bold{u}-\widetilde{\bold{u}}) \nabla_x v \cdot (\left[ \Delta_x(\bold{u}-\widetilde{\bold{u}})- \nabla_x div_x(\bold{u}-\widetilde{\bold{u}})\right]) \, dx_1 \, dx' \\ 
:=&\sum\limits_{i=1}^5J_i(t). \\
\end{aligned}    
\end{equation*}
Recall $\varepsilon_1=\sup\limits_{0\le t\le T}\lVert (v-\widetilde{v}), \bold{u}-\widetilde{\bold{u}}) \rVert_{H^2}$.

Using \eqref{GNIAPP3},
\begin{equation}
\begin{aligned} \label{GNIAPP4}
\lVert \nabla_x(\bold{u}-\widetilde{\bold{u}}) \rVert_{L^3}&\le C\lVert \nabla_x(\bold{u}-\widetilde{\bold{u}}) \rVert_{H^1}\le C\varepsilon_1, \\
\lVert \nabla_x(v-\widetilde{v}) \rVert_{L^3}&\le C\lVert \nabla_x(v-\widetilde{v}) \rVert_{H^1}\le C\varepsilon_1.    
\end{aligned}    
\end{equation}
From this, we get
\begin{equation*}
\begin{aligned}
J_1(t)&\le C\lVert \nabla_x(\bold{u}-\widetilde{\bold{u}}) \rVert_{L^3}\lVert \nabla_x(\bold{u}-\widetilde{\bold{u}}) \rVert_{L^6}\lVert \nabla_x(\bold{u}-\widetilde{\bold{u}}) \rVert_{L^2}+C\lVert \nabla_x \widetilde{\bold{u}} \rVert_{L^\infty}\lVert \nabla_x(\bold{u}-\widetilde{\bold{u}}) \rVert_{L^2}^2 \\
&\le C\varepsilon_1 \lVert \nabla_x(\bold{u}-\widetilde{\bold{u}}) \rVert_{H^1}\lVert \nabla_x(\bold{u}-\widetilde{\bold{u}}) \rVert_{L^2} + C(\delta_1^2+\delta_2^2)\lVert \nabla_x(\bold{u}-\widetilde{\bold{u}}) \rVert_{L^2}^2 \\
&\le C(\varepsilon_1+\delta_1^2+\delta_2^2)\left( \mathcal{D}_2(t)+\lVert \nabla_x(\bold{u}-\widetilde{\bold{u}}) \rVert_{L^2}^2 \right).
\end{aligned}    
\end{equation*}
By Young's inequality,
\begin{equation*}
\begin{aligned}
J_2(t)&\le \sum\limits_{i=1}^2 | \dot{X}_i(t)| \lVert (\widetilde{u}_{i})_{x_1}^{X_i} \rVert_{L^2}\sqrt{\mathcal{D}_2}\le C\sum\limits_{i=1}^2 | \dot{X}_i(t) |\delta_i^{\frac{3}{2}}\sqrt{\mathcal{D}_2}\le C\sum\limits_{i=1}^2 \delta_i^2 | \dot{X}_i(t) |^2+C\sum\limits_{i=1}^2 \delta_i\mathcal{D}_2(t).   
\end{aligned}    
\end{equation*}
As in $J_2$ and $\mathcal{B}_4$, we get
\begin{equation*}
\begin{aligned}
J_3(t)&\le C\sum\limits_{i=1}^2 \delta_i\sqrt{\frac{\delta_i}{\nu_i}}(\sqrt{\bold{G}_1}+\sqrt{\bold{G}_3}+\sqrt{\mathcal{G}^S}+\sqrt{\bold{D}})\sqrt{\mathcal{D}_2}\\
&\le C(\delta_1+\delta_2) (\bold{G}_1+\bold{G}_3+\mathcal{G}^S+\bold{D}+\mathcal{D}_2).
\end{aligned}    
\end{equation*}
By Young's inequality,
\begin{equation*}
\begin{aligned}
J_4(t)&\le C\int_{\mathbb{T}^2} \int_{\mathbb{R}} v\left| \Delta_x(\bold{u}-\widetilde{\bold{u}}) \right| \left|  \nabla_x(p(v)-p(\widetilde{v})) \right|	 \, dx_1 \, dx' \\
&\le \frac{\mu}{5}\int_{\mathbb{T}^2} \int_{\mathbb{R}} v\left| \Delta_x(\bold{u}-\widetilde{\bold{u}}) \right|^2	 \, dx_1 \, dx'+C\lVert \nabla_x(p(v)-p(\widetilde{v})) \rVert_{L^2}^2 \\
&\le \frac{1}{5}\mathcal{D}_2+C\bold{D}.    
\end{aligned}    
\end{equation*}
Finally, by Holder's inequality and Gagliardo-Nirenberg interpolation inequality,
\begin{equation*}
\begin{aligned}
J_5\le& C\int_{\mathbb{T}^2} \int_{\mathbb{R}} \left| div_x(\bold{u}-\widetilde{\bold{u}}) \right|\left| \nabla_x(v-\widetilde{v}) \right|\left( \left| \Delta_x(\bold{u}-\widetilde{\bold{u}}) \right|+\left| \nabla_x div_x(\bold{u}-\widetilde{\bold{u}}) \right|  \right)	 \, dx_1 \, dx'\\
&+C\int_{\mathbb{T}^2} \int_{\mathbb{R}} \left| div_x(\bold{u}-\widetilde{\bold{u}}) \right|\left| (\widetilde{v})_{x_1} \right|\left( \left| \Delta_x(\bold{u}-\widetilde{\bold{u}}) \right|+\left| \nabla_x div_x(\bold{u}-\widetilde{\bold{u}}) \right|  \right)	 \, dx_1 \, dx'\\
\le& C\lVert div_x(\bold{u}-\widetilde{\bold{u}}) \rVert_{L^6}\lVert \nabla_x(v-\widetilde{v}) \rVert_{L^3}\left( \lVert \Delta_x(\bold{u}-\widetilde{\bold{u}}) \rVert_{L^2}+\lVert \nabla_x div_x(\bold{u}-\widetilde{\bold{u}}) \rVert_{L^2} \right)\\
&+C\lVert div_x(\bold{u}-\widetilde{\bold{u}}) \rVert_{L^2}\lVert( \widetilde{v})_{x_1} \rVert_{L^\infty}\left( \lVert \Delta_x(\bold{u}-\widetilde{\bold{u}}) \rVert_{L^2}+\lVert \nabla_x div_x(\bold{u}-\widetilde{\bold{u}}) \rVert_{L^2} \right) \\
\le& C(\delta_1^2+\delta_2^2+\varepsilon_1)\sqrt{\mathcal{D}_2}\lVert div_x(\bold{u}-\widetilde{\bold{u}}) \rVert_{H^1} \\
\le &C(\delta_1^2+\delta_2^2+\varepsilon_1)(\mathcal{D}_2+\lVert \nabla_x(\bold{u}-\widetilde{\bold{u}}) \rVert_{L^2}^2).  
\end{aligned}    
\end{equation*}
Therefore, from the estimates above, we have
\begin{equation*}
\begin{aligned}
&\frac{d}{dt}\lVert \nabla_x(\bold{u}-\widetilde{\bold{u}}) \rVert_{L^2}^2+\mathcal{D}_2(t)\\
\le&\mathcal{C}_2\left(\sum\limits_{i=1}^2\delta_i^2| \dot{X}_i(t) |^2+\bold{G}_1+\bold{G}_3+\mathcal{G}^S+\bold{D}+\lVert \nabla_x(\bold{u}-\widetilde{\bold{u}}) \rVert_{L^2}^2)\right),
\end{aligned}    
\end{equation*}
where $\mathcal{C}_2>0$ is some constant.

Integrate the inequality above over $\left[ 0,t \right]$ for any $t\le T$ and then multiply the result by $\frac{1}{2\max\left\{ 1, \mathcal{C}_2 \right\}}$.

Then we have
\begin{align} 
\begin{aligned}\label{lemma5.2ieq2}
&\frac{1}{2\max\left\{ 1, \mathcal{C}_2 \right\}}\lVert \nabla_x(\bold{u}-\widetilde{\bold{u}}) \rVert_{L^2}^2 + \frac{1}{2\max\left\{ 1, \mathcal{C}_2 \right\}}\int_{0}^{t}\mathcal{D}_2\, d\tau \\
\le &  \frac{1}{2\max\left\{ 1, \mathcal{C}_2 \right\}}\lVert \nabla_{x}\left(\bold{u}_0-\widetilde{\bold{u}}(0,\cdot)\right) \rVert_{L^2}^2 \\ 
& \ + \frac{1}{2}\int_{0}^t \left(\sum\limits_{i=1}^2{\delta_i}|\dot{X}_i|^2+\bold{G}_1+\bold{G}_3+\mathcal{G}^S+\bold{D}\right) \, d\tau+\frac{1}{2}\int_{0}^{t} \lVert \nabla_x(\bold{u}-\widetilde{\bold{u}}) \rVert_{L^2}^2  \, d\tau.\\
\end{aligned}
\end{align}
Using  $\lVert \nabla_x^2(\bold{u}-\widetilde{\bold{u}})(t) \rVert_{L^2}^2 \sim \mathcal{D}_2(t)$ and adding \eqref{lemma5.2ieq2} to \eqref{lemma5.1}, we finish the proof of the Lemma \ref{lem:2nd}.

\section{Proof of Lemma \ref{lem:3rd}}
We set $\phi:=v-\widetilde{v}$, $\psi:=\bold{u}-\widetilde{\bold{u}}$ for simplicity, and substitute this into $(5.2)$ as 
\begin{equation}
\begin{cases}
\partial_t\phi+\bold{u}\cdot \nabla_x\phi- \sum\limits_{i=1}^2\dot{X}_i(t)(\widetilde{v}_i)_{x_1}^{X_i}+v\sum\limits_{i=1}^2F_i(\widetilde{v}_i)_{x_1}^{X_i}=vdiv_x\psi, \\
\partial_t\psi+\bold{u}\nabla_x\psi+vp'(v)\nabla_x\phi+v(p'(v)-p'(\widetilde{v}))\nabla_x\widetilde{v}-\sum\limits_{i=1}^2\dot{X}_i(t)(\widetilde{\bold{u}}_i)_{x_1}+v\sum\limits_{i=1}^2 F_i (\widetilde{\bold{u}}_i)_{x_1}\\
\ = \mu v \Delta_x\psi+(\mu+\lambda)v\nabla_x div_x\psi-\nabla_x( p(\widetilde{v})-p(\widetilde{v}_1^{X_1})-p(\widetilde{v}_2^{X_2})).
\end{cases}    
\end{equation}
Taking $\nabla_x \partial_{x_j}$ $(j=1,2,3)$ to $(5.9)_1$ and $\partial_{x_j}$ $(j=1,2,3)$ to $(5.9)_2$, we have
\begin{equation} \label{5.12}
\begin{cases}
\partial_t\nabla_x \partial_{x_j}\phi+\bold{u} \nabla_x(\nabla_x \partial_{x_j}\phi)- \sum\limits_{i=1}^2\dot{X}_i(t)\nabla_x \partial_{x_j}(\widetilde{v}_i)_{x_1}^{X_i}\\
+v\sum\limits_{i=1}^2F_i\nabla_x \partial_{x_j}(\widetilde{v}_i)_{x_1}^{X_i}+\nabla_x \partial_{x_j}\bold{u} \nabla_x \phi+\nabla_x \bold{u} \nabla_x\partial_{x_j} \phi
+\partial_j\bold{u} \nabla_x(\nabla_x\phi)\\+\sum\limits_{i=1}^2\nabla_x \partial_{x_j}(vF_i)(\widetilde{v}_i)_{x_1}^{X_i}+\sum\limits_{i=1}^2\nabla_x(vF_i)\partial_{x_j}(\widetilde{v}_i)_{x_1}^{X_i}+\sum\limits_{i=1}^2\partial_{x_j}(vF_i)\nabla_x(\widetilde{v}_i)_{x_1}^{X_i} \\ \  =v\nabla_x \partial_{x_j}div_x(\psi)+\nabla_x \partial_{x_j} vdiv_x\psi+\partial_{x_j}v \nabla_x div_x\psi + \nabla_x v \partial_{x_j} div_x\psi, \\
\partial_t\partial_{x_j}\psi+\bold{u} \nabla_x\partial_{x_j}\psi+\partial_{x_j}\bold{u} \nabla_x\psi+vp'(v)\nabla_x\partial_{x_j}\phi+\partial_{x_j}(vp'(v))\nabla_x\phi\\
+\partial_{x_j}(v(p'(v)-p'(\widetilde{v})))\nabla_x\widetilde{v} +v(p'(v)-p'(\widetilde{v})))\nabla_x\partial_{x_j}(\widetilde{v})\\-\sum\limits_{i=1}^2\dot{X}_i(t)\partial_{x_j}(\widetilde{\bold{u}}_i)_{x_1}+\sum\limits_{i=1}^2vF_i\partial_{x_j}(\widetilde{\bold{u}}_i)_{x_1}+\sum\limits_{i=1}^2\partial_{x_j}(vF_i)(\widetilde{\bold{u}}_i)_{x_1}\\
\ =\mu v \Delta_x\partial_{x_j}\psi+(\mu+\lambda)v\nabla_x \partial_{x_j}div_x\psi\\
\quad +\partial_{x_j}v(\mu\Delta_x\psi+(\mu+\lambda)\nabla_x div_x\psi)-\partial_{x_j}\nabla_x( p(\widetilde{v})-p(\widetilde{v}_1^{X_1})-p(\widetilde{v}_2^{X_2}) ).
\end{cases}
\end{equation}
Multiply \eqref{5.12}$_1$ by $\rho(2\mu+\lambda)\nabla_x\partial_{x_j}\phi$ and then sum $i$ from $1$ to $3$.

In addition, integrating the result over $\Omega$, we have
\begin{equation} \label{5.13}
\begin{aligned}
&(2\mu+\lambda)\frac{d}{dt}\int_{\mathbb{T}^2} \int_{\mathbb{R}} \rho \frac{|\nabla_x^2\phi|^2}{2}  \, dx_1 \, dx'-(2\mu+\lambda)\sum\limits_{j=1}^3\int_{\mathbb{T}^2} \int_{\mathbb{R}} \nabla_x\partial_{x_j}\phi\cdot \nabla_x\partial_{x_j} div_x\psi \, dx_1 \, dx'\\
=&(2\mu+\lambda)\sum\limits_{i=1}^2\dot{X}_i(t)\int_{\mathbb{T}^2} \int_{\mathbb{R}}  \rho \partial_{x_1}^2\phi (\widetilde{v}_i)_{x_1x_1x_1}^{X_i} \, dx_1 \, dx'-(2\mu+\lambda)\sum\limits_{i=1}^2\int_{\mathbb{T}^2} \int_{\mathbb{R}}  F_i \partial_{x_1}\phi (\widetilde{v}_i)_{x_1x_1x_1}^{X_i} \, dx_1 \, dx'\\
&-(2\mu+\lambda)\sum\limits_{j=1}^3\int_{\mathbb{T}^2} \int_{\mathbb{R}}  \rho \nabla_x \partial_{x_j}\phi \cdot 	\left[ \nabla_x\partial_j\bold{u} \nabla_x \phi +\nabla_x \bold{u} \nabla_x\partial_{x_j} \phi+\partial_{x_j}\bold{u} \nabla_x(\nabla_x\phi)  \right] \, dx_1 \, dx'\\
&-(2\mu+\lambda)\sum\limits_{i=1}^2\sum\limits_{j=1}^3	\int_{\mathbb{T}^2} \int_{\mathbb{R}} \rho\left( \nabla_x\partial_{x_j} \phi \cdot \nabla_x \partial_{x_j}(vF_i)(\widetilde{v}_i)_{x_1}^{X_i}+\partial_{x_1}\partial_{x_j}\phi \partial_{x_j}(vF_i)(\widetilde{v}_i)_{x_1x_1}^{X_i} \right) \, dx_1 \, dx'\\
&-(2\mu+\lambda)\sum\limits_{i=1}^2 \int_{\mathbb{T}^2} \int_{\mathbb{R}} \rho \nabla_x\partial_{x_1}\phi \cdot \nabla_x(vF_i)(\widetilde{v}_i)_{x_1x_1}^{X_i} \, dx_1 \, dx'\\
&+(2\mu+\lambda)\sum\limits_{j=1}^3 \int_{\mathbb{T}^2} \int_{\mathbb{R}} \rho \nabla_x\partial_{x_j}\phi \cdot 	\left[ \nabla_x \partial_{x_j} v div_x \psi + \partial_{x_j} v \nabla_x div_x \psi + \nabla_x v \partial_{x_j} div_x\psi \right]  \, dx_1 \, dx'.    
\end{aligned}    
\end{equation}
Especially, in this calculation above, we significantly use the equality below which comes from integration by parts 
\begin{equation*}
\begin{aligned}
\int_{\mathbb{T}^2} \int_{\mathbb{R}} -\rho_t \frac{\left| \nabla_x^2 \phi \right|^2}{2}+\rho \nabla_x ^2 \phi \cdot \bold{u}(\nabla_x^3 \phi) \, dx_1 \, dx'&=\int_{\mathbb{T}^2} \int_{\mathbb{R}} -\rho_t \frac{\left| \nabla_x^2 \phi \right|^2}{2}- div_x (\rho \bold{u}) \frac{\left| \nabla_x^2 \phi \right|^2}{2} \, dx_1 \, dx' \\ 
&=0.
\end{aligned}   
\end{equation*}
Multiply \eqref{5.12}$_1$ by $\rho(2\mu+\lambda)\nabla_x\partial_{x_j}\phi$ and then sum $i$ from $1$ to $3$.

In addition, integrating the result over $\Omega$, we have
\begin{equation} \label{5.14}
\begin{aligned}
&\int_{\mathbb{T}^2} \int_{\mathbb{R}} -p'(v)|\nabla_x^2 \phi|^2  \, dx_1 \, dx'+(2\mu+\lambda)\sum\limits_{j=1}^3\int_{\mathbb{T}^2} \int_{\mathbb{R}} \nabla_x\partial_{x_j}\phi \cdot \nabla_x\partial_{x_j} div_x\psi \, dx_1 \, dx'\\
=&\frac{d}{dt}\sum\limits_{j=1}^3\int_{\mathbb{T}^2} \int_{\mathbb{R}} \rho\partial_{x_j}\psi \cdot \nabla_x\partial_{x_j}\phi \, dx_1 \, dx'+\sum\limits_{j=1}^3\int_{\mathbb{T}^2} \int_{\mathbb{R}} \rho\nabla_x\partial_{x_j}\phi \cdot \partial_{x_j}\bold{u} \nabla_x \psi \, dx_1 \, dx'\\
&-\sum\limits_{j=1}^3\int_{\mathbb{T}^2} \int_{\mathbb{R}} \rho \partial_{x_j}\psi\cdot \left[ \nabla_x\partial_{x_j} \partial_t \phi +\bold{u} \nabla_x(\nabla_x\partial_{x_j}\phi) \right] \, dx_1 \, dx'\\
&+\sum\limits_{j=1}^3 \int_{\mathbb{T}^2} \int_{\mathbb{R}} \rho \partial_{x_j}(vp'(v)) \nabla_x\partial_{x_j}\phi \cdot \nabla_x\phi  \, dx_1 \, dx'-\sum\limits_{i=1}^2\dot{X}_i(t) \int_{\mathbb{T}^2} \int_{\mathbb{R}} \rho \partial_{x_1}^2 \phi (\widetilde{u}_{i})_{x_1x_1}^{X_i} \, dx_1 \, dx'\\
&+\sum\limits_{i=1}^2\sum\limits_{j=1}^3\int_{\mathbb{T}^2} \int_{\mathbb{R}} \rho \partial_{x_j}(v(p'(v)-p'(\widetilde{v})))\partial_{x_1} \partial_{x_j} \phi (\widetilde{v}_i)_{x_1}^{X_i} \, dx_1 \, dx'\\
&+\sum\limits_{i=1}^2 \int_{\mathbb{T}^2} \int_{\mathbb{R}} (p'(v)-p'(\widetilde{v}))\partial_{x_1}^2 \phi (\widetilde{v}_i)_{x_1x_1}^{X_i}  \, dx_1 \, dx'\\
&+\sum\limits_{i=1}^2\int_{\mathbb{T}^2} \int_{\mathbb{R}} F_i \partial_{x_1}^2 \phi (\widetilde{u}_{i})_{x_1x_1}^{X_i} \, dx_1 \, dx'+\sum\limits_{i=1}^2\sum\limits_{j=1}^3\int_{\mathbb{T}^2} \int_{\mathbb{R}}  \rho \partial_{x_j}(vF_i)\partial_{x_1} \partial_{x_j} \phi (\widetilde{u}_{i})_{x_1}^{X_i} \, dx_1 \, dx'\\
&-\sum\limits_{j=1}^3 \int_{\mathbb{T}^2} \int_{\mathbb{R}} \rho \partial_jv (\mu\Delta_x\psi + (\mu+\lambda)\nabla_x div_x \psi) \cdot \nabla_x\partial_{x_j} \phi \, dx_1 \, dx'\\
&+\sum\limits_{j=1}^3 \int_{\mathbb{T}^2} \int_{\mathbb{R}} \rho \nabla_x \partial_{x_j}\phi \cdot \partial_{x_j}\nabla_x( p(\widetilde{v})-p(\widetilde{v}_1^{X_1})-p(\widetilde{v}_2^{X_2})) \, dx_1 \, dx'  
\end{aligned}    
\end{equation}
As in \eqref{5.13}, we use the equality below in this calculation above.
\begin{equation*}
 \int_{\mathbb{T}^2} \int_{\mathbb{R}} \nabla_x\partial_{x_j}\phi \cdot (\Delta_x\partial_{x_j}\psi-\nabla_x\partial_{x_j}div_x\psi) \, dx_1 \, dx'=-\int_{\mathbb{T}^2} \int_{\mathbb{R}} \partial_{x_j}\phi \nabla_x \cdot (\Delta_x\partial_{x_j}\psi-\nabla_x\partial_{x_j}div_x\psi) \, dx_1 \, dx',
\end{equation*}
from $\nabla_x\cdot \left[ \Delta_x \partial_{x_j}\psi-\nabla_x div_x \partial_{x_j}\psi \right]=0$.    

Adding \eqref{5.13} and \eqref{5.14} together and then integrating the result over $\left[ 0, t \right]$ for $t\in \left[ 0, T \right]$, we get
\begin{equation} \label{5.15}
(2\mu+\lambda)\int_{\mathbb{T}^2} \int_{\mathbb{R}} \rho \frac{\left| \nabla_x^2 \phi \right|^2}{2} \, dx_1 \, dx' \Big|^{\tau=t}_{\tau=0} \  + \int_{0}^{t} \int_{\mathbb{T}^2} \int_{\mathbb{R}} \rho \frac{\left| \nabla_x^2 \phi \right|^2}{2} \, dx_1 \, dx' \, d\tau=\sum\limits_{j=1}^8 K_j(t),  
\end{equation}
where
\begin{equation*}
\begin{aligned}
K_1(t)=&\sum\limits_{j=1}^3\int_{\mathbb{T}^2} \int_{\mathbb{R}} \rho \partial_{x_j}\psi \cdot \nabla_x\partial_{x_j}\phi \, dx_1\, dx' \Big|^{\tau=t}_{\tau=0},\\
K_2(t)=&\sum\limits_{i=1}^2 \int_{0}^{t} \dot{X}_i(\tau) \int_{\mathbb{T}^2} \int_{\mathbb{R}} \rho \partial_{x_1}^2\phi \left[(2\mu+\lambda)(\widetilde{v}_i)_{x_1x_1x_1}^{X_i}- (\widetilde{u}_{i})_{x_1x_1}^{X_i} \right] \, dx_1 \, dx' \, d\tau,\\
K_3(t)=&-\sum\limits_{i=1}^2 \int_{0}^{t} \int_{\mathbb{T}^2} \int_{\mathbb{R}} F_i \partial_{x_1}^2\phi \left[(2\mu+\lambda)(\widetilde{v}_i)_{x_1x_1x_1}^{X_i}- (\widetilde{u}_{i})_{x_1x_1}^{X_i} \right] \, dx_1 \, dx' \, d\tau,\\
K_4(t)=&-(2\mu+\lambda)\sum\limits_{i=1}^2\sum\limits_{j=1}^3 \int_{0}^{t} \int_{\mathbb{T}^2} \int_{\mathbb{R}} \rho\left( \nabla_x\partial_{x_j} \phi \cdot \nabla_x \partial_{x_j}(vF_i)(\widetilde{v}_i)_{x_1}^{X_i}+\partial_{x_1}\partial_{x_j}\phi \partial_{x_j}(vF_i)(\widetilde{v}_i)_{x_1x_1}^{X_i} \right)  \, dx_1 \, dx' \, d\tau\\
&-(2\mu+\lambda) \sum\limits_{i=1}^2 \int_{0}^{t} \int_{\mathbb{T}^2} \int_{\mathbb{R}}  \rho \nabla_x\partial_{x_1}\phi \cdot \nabla_x(vF_i)(\widetilde{v}_i)_{x_1x_1}^{X_i} \, dx_1 \, dx' \, d\tau \\
&+\sum\limits_{i=1}^2\sum\limits_{j=1}^3 \int_{0}^{t} \int_{\mathbb{T}^2} \int_{\mathbb{R}}  \rho \partial_{x_j}(vF_i)\partial_{x_1} \partial_{x_j} \phi (\widetilde{u}_{i})_{x_1}^{X_i}  \, dx_1 \, dx' \, d\tau,\\
K_5(t)=&-\sum\limits_{j=1}^3 \int_{0}^{t} \int_{\mathbb{T}^2} \int_{\mathbb{R}} \rho \partial_{x_j}\psi\cdot \left[ \nabla_x\partial_{x_j} \partial_t \phi +\bold{u}\nabla_x(\nabla_x\partial_{x_j}\phi) \right] \, dx_1 \, dx' \, d\tau,\\
K_6(t)=&-(2\mu+\lambda)\sum\limits_{j=1}^3 \int_{0}^{t} \int_{\mathbb{T}^2} \int_{\mathbb{R}}  \rho \nabla_x \partial_{x_j}\phi \cdot 	\left[ \nabla_x\partial_{x_j}\bold{u} \nabla_x \phi +\nabla_x \bold{u} \nabla_x\partial_{x_j} \phi+\partial_{x_j}\bold{u} \nabla_x(\nabla_x\phi)  \right] \, dx_1 \, dx' \, d\tau\\
&+\sum\limits_{j=1}^3 \int_{0}^{t} \int_{\mathbb{T}^2} \int_{\mathbb{R}} \rho\nabla_x\partial_{x_j}\phi \cdot \left[ \partial_{x_j}\bold{u}  \nabla_x \psi +  \partial_{x_j}(vp'(v)) \nabla_x\phi \right]  \, dx_1 \, dx' \, d\tau, \\
K_7(t)=&\sum\limits_{i=1}^2 \sum\limits_{j=1}^3 \int_{0}^{t} \int_{\mathbb{T}^2} \int_{\mathbb{R}}  \rho \partial_{x_j}(v(p'(v)-p'(\widetilde{v})))\partial_{x_1} \partial_{x_j} \phi (\widetilde{v}_i)_{x_1}^{X_i}  \, dx_1 \, dx' \, d\tau\\
&+\sum\limits_{i=1}^2\int_{0}^{t} \int_{\mathbb{T}^2} \int_{\mathbb{R}}  (p'(v)-p'(\widetilde{v}))\partial_{x_1}^2 \phi (\widetilde{v}_i)_{x_1x_1}^{X_i}  \, dx_1 \, dx' \, d\tau,\\
K_8(t)=&(2\mu+\lambda)\sum\limits_{j=1}^3 \int_{0}^{t} \int_{\mathbb{T}^2} \int_{\mathbb{R}} \rho \nabla_x \partial_{x_j}\phi \cdot 	\left[ \nabla_x \partial_{x_j} v div_x \psi + \partial_{x_j} v \nabla_x div_x \psi + \nabla_x v \partial_{x_j} div_x\psi \right]  \, dx_1 \, dx' \, d\tau\\
&-\sum\limits_{j=1}^3 \int_{0}^{t} \int_{\mathbb{T}^2} \int_{\mathbb{R}} \rho \partial_{x_j}v \nabla_x \partial_{x_j} \phi \cdot \left[ \mu\Delta_x\psi + (\mu+\lambda)\nabla_x div_x \psi \right] \, dx_1 \, dx' \, d\tau,\\
K_9(t)=&\sum\limits_{j=1}^3 \int_{0}^{t} \int_{\mathbb{T}^2} \int_{\mathbb{R}} \rho \nabla_x \partial_{x_j}\phi \cdot \partial_{x_j}\nabla_x( p(\widetilde{v})-p(\widetilde{v}_1^{X_1})-p(\widetilde{v}_2^{X_2}))  \, dx_1 \, dx' \, d\tau. \\
\end{aligned}    
\end{equation*}
By Young's inequality, 
\begin{equation*}
K_1(t)\le \frac{2\mu+\lambda}{16}\lVert (\sqrt{\rho} \nabla_x^2\phi)(t) \rVert_{L^2}^2+C\left( \lVert ( \nabla_x\psi)(t) \rVert_{L^2}^2 + \lVert ( \nabla_x^2\phi)(0) \rVert_{L^2}^2 + \lVert ( \nabla_x\psi)(0) \rVert_{L^2}^2\right).    
\end{equation*}
By Lemma \ref{lem:shock-est}, 
\begin{equation*}
K_2(t)\le C\sum\limits_{i=1}^2 \delta_i \int_{0}^{t} | \dot{X}_i(\tau) |  \lVert \partial_{x_1}^2\phi \rVert_{L^2} \lVert (\widetilde{v}_i)_{x_1}^{X_i} \rVert_{L^2}  \, d\tau\le C\sum\limits_{i=1}^2 \delta_i^2 \int_{0}^{t} | \dot{X}_i(\tau) |^2  \lVert \sqrt{|p'(v)|} \partial_{x_1}^2\phi \rVert_{L^2}^2  \, d\tau.
\end{equation*}
Using Lemma \ref{lem:shock-est}, as in $\mathcal{B}_4$, we get  
\begin{equation*}
\begin{aligned}
K_3(t)\le &C\sum\limits_{i=1}^2 \delta_i^2 \int_{0}^{t} (\sqrt{\bold{G}_1}+\sqrt{\bold{G}_3}+\sqrt{\mathcal{G}^S}+\sqrt{\bold{D}}\\
& \qquad \qquad \qquad +C\left( \varepsilon_1 \nu_i \delta_i e^{-C\delta_it}+\delta_1\delta_2 e^{-C\min\left\{ \delta_1, \delta_2  \right\}t} \right) \int_{\mathbb{T}^2} \int_{\mathbb{R}} \eta(U|\widetilde{U}) \, dx_1 \, dx' \lVert \partial_{x_1}^2 \phi \rVert_{L^2}) \, d\tau\\
\le &C\sum\limits_{i=1}^2 \delta_i^2 \int_{0}^{t} (\bold{G}_1+\bold{G}_3+\mathcal{G}^S+\bold{D}+ \lVert \sqrt{|p'(v)|}\partial_{x_1}^2 \phi \rVert_{L^2}^2) d\tau+C\sum\limits_{i=1}^2\delta_i^2\varepsilon_1.    
\end{aligned}
\end{equation*}
To estimate $K_4(t)$, note that 
\begin{equation}
\begin{aligned} \label{presentvF}
vF_i&=\sigma_i^*(v-\widetilde{v})+\sigma_i^*(\widetilde{v}-\widetilde{v}_i^{X_i})+(u_1-\widetilde{u})+(\widetilde{u}-\widetilde{u}_{i}^{X_i})\\
&=\sigma_i^*\phi+\psi_1+\sigma_i^*(\widetilde{v}-\widetilde{v}_i^{X_i})+(\widetilde{u}-\widetilde{u}_{i}^{X_i}).    
\end{aligned}    
\end{equation}
From this,
\begin{equation*}
\begin{aligned}
K_4(t)\le &C\sum\limits_{i=1}^2 \int_{0}^{t} \int_{\mathbb{T}^2} \int_{\mathbb{R}} |\nabla_x^2 \phi| (\left| \nabla_x^2 \phi \right| +\left| \nabla_x^2 \psi \right|)|(\widetilde{v}_i)_{x_1}^{X_i} | \, dx_1 \, dx'  \, d\tau\\
&+C\int_{0}^{t} \int_{\mathbb{T}^2} \int_{\mathbb{R}}  \left| \nabla_x^2 \phi \right| | (\widetilde{v}_1)_{x_1}^{X_1} | | (\widetilde{v}_2)_{x_1}^{X_2} | \, dx_1 \, dx'  \, d\tau\\
&+C\sum\limits_{i=1}^2 \int_{0}^{t} \int_{\mathbb{T}^2} \int_{\mathbb{R}} \left| \nabla_x^2 \phi \right| (\left| \nabla_x \phi \right| +\left| \nabla_x \psi \right|)| (\widetilde{v}_i)_{x_1}^{X_i} |  \, dx_1 \, dx'  \, d\tau\\
\le &C\sum\limits_{i=1}^2 \delta_i^2 \int_{0}^{t} \lVert \nabla_x^2 \phi  \rVert_{L^2} (\lVert \nabla_x \phi  \rVert_{H^1}+\lVert \nabla_x \psi  \rVert_{H^1}) d\tau+C\int_{0}^{t}  \lVert \sqrt{|p'(v)|}\nabla_x^2 \phi  \rVert_{L^2}^2 \, d\tau+C\int_{0}^{t}\delta_1\delta_2 e^{-C\min\left\{ \delta_1, \delta_2 \right\}t}d\tau\\
\le &C\sum\limits_{i=1}^2 \delta_i^2 \int_{0}^{t} \left( \lVert \sqrt{|p'(v)|}\nabla_x^2 \phi  \rVert_{L^2} + \bold{D}(\tau)+ \mathcal{G}^S(\tau)+\varepsilon_1\nu_i \delta_ie^{-C\delta_it}+\lVert \nabla_x \psi  \rVert_{H^1}^2 \right) \, d\tau + C\delta_0.    
\end{aligned}    
\end{equation*}
Using the equation $(5.10)_1$, we can compute the term $K_5(t)$ as
\begin{equation*}
\begin{aligned}
K_5(t)=&-\sum\limits_{i=1}^2\int_{0}^{t}\dot{X}_i(\tau) \int_{\mathbb{T}^2} \int_{\mathbb{R}} \rho\partial_{x_1}\psi_1 (\widetilde{v}_i)_{x_1x_1x_1}^{X_i} dx_1  \, dx'  \, d\tau \\
&+\sum\limits_{i=1}^2\int_{0}^{t}\dot{X}_i(\tau) \int_{\mathbb{T}^2} \int_{\mathbb{R}} F_i \partial_{x_1}\psi_1 (\widetilde{v}_i)_{x_1x_1x_1}^{X_i} dx_1  \, dx'  \, d\tau \\
    &+\sum\limits_{j=1}^3\int_{0}^{t} \int_{\mathbb{T}^2} \int_{\mathbb{R}} \rho \partial_{x_j} \psi \cdot \left[ \nabla_x\partial_j\bold{u} \nabla_x\phi+\nabla_x \bold{u} \nabla_x\partial_{x_j}\phi +\partial_{x_j} \bold{u} \nabla_x(\nabla_x \phi)   \right] dx_1  \, dx'  \, d\tau\\
&+\sum\limits_{j=1}^3\sum\limits_{i=1}^2\int_{0}^{t} \int_{\mathbb{T}^2} \int_{\mathbb{R}} \rho \partial_j \psi \cdot \left[ \nabla_x\partial_{x_j}(vF_i)(\widetilde{v}_i)_{x_1}^{X_i}+ \nabla_x(vF_i)\partial_{x_j}(\widetilde{v}_i)_{x_1}^{X_i})+\partial_{x_j}(vF_i)\nabla_x(\widetilde{v}_i)_{x_1}^{X_i} \right] dx_1  \, dx'  \, d\tau\\
&-\sum\limits_{j=1}^3\int_{0}^{t} \int_{\mathbb{T}^2} \int_{\mathbb{R}} \rho \partial_{x_j} \psi \cdot \left[ \nabla_x\partial_{x_j}v div_x\psi + \partial_{x_j}v \nabla_x div_x\psi+\nabla_x v \partial_{x_j} div_x\psi \right] dx_1  \, dx'  \, d\tau\\
&-\sum\limits_{j=1}^3 \int_{0}^{t} \int_{\mathbb{T}^2} \int_{\mathbb{R}} \partial_{x_j} \psi \cdot \nabla_x \partial_{x_j} div_x\psi dx_1  \, dx'  \, d\tau\\
=:&\sum\limits_{j=1}^6 K_{5,j}(t).
\end{aligned}
\end{equation*}
By Lemma \ref{lem:shock-est}, we get
\begin{equation*}
K_{5,1}(t)\le C\sum\limits_{i=1}^2 \delta_i^2 \int_{0}^{t} |\dot{X}_i(\tau)|^2 d\tau + C\sum\limits_{i=1}^2 \delta_i^2 \int_{0}^{t} \lVert \nabla_x \psi \rVert_{L^2}^2 d\tau.  
\end{equation*}
As in $\mathcal{B}_4$, we get
\begin{equation*}
\begin{aligned}
K_{5,2}(t)&\le C\sum\limits_{i=1}^2 \delta_i^2 \int_{0}^{t} F_i^2+|\partial_{x_1} \psi_1|^2 \, d\tau\\
&\le C \sum\limits_{i=1}^2 \delta_i^2 \int_{0}^{t} \left( \lVert \nabla_x \psi \rVert_{L^2}^2 + \bold{G}_1 + \bold{G}_3 +\mathcal{G}^S +\bold{D} \right)  \, d\tau+C\sum\limits_{i=1}^2\delta_i^2\varepsilon_1.
\end{aligned}
\end{equation*}
By \eqref{GNIAPP4}, similarly to $\mathcal{B}_4$, we have
\begin{equation*}
\begin{aligned}
K_{5,3}(t)\le &C\int_{0}^{t}  \int_{\mathbb{T}^2} \int_{\mathbb{R}} \left| \nabla_x \psi \right|  \left(  \left| \nabla_x^2 \psi \right| \left| \nabla_x \phi \right| +   \left| \nabla_x \psi \right| \left| \nabla_x^2 \phi \right| \right) \,dx_1 \, dx' \, d\tau\\
&+C\int_{0}^{t}  \int_{\mathbb{T}^2} \int_{\mathbb{R}} \left| \nabla_x \psi \right|  \left(  \left| (\widetilde{u})_{x_1x_1} \right| \left| \partial_{x_1} \phi \right| +   \left| (\widetilde{u})_{x_1} \right| \left| \nabla_x \partial_{x_1} \phi \right| \right) \,dx_1 \, dx' \, d\tau\\
\le &C\int_{0}^{t} \lVert \nabla_x \psi \rVert_{L^6} \lVert \nabla_x(\phi, \psi) \rVert_{L^3} \lVert \nabla_x^2 (\phi, \psi) \rVert_{L^2} \, d\tau+C\sum\limits_{i=1}^2 \delta_i^2 \int_{0}^{t} \lVert \nabla_x \psi \rVert_{L^2} (\lVert \partial_{x_1} \phi \rVert_{L^2}+ \lVert \nabla_x \partial_{x_1} \phi \rVert_{L^2}) \, d\tau\\
\le &C\varepsilon_1 \int_{0}^{t} (\lVert \sqrt{|p'(v)|} \nabla_x \phi \rVert_{L^2}^2+ \lVert \nabla_x \psi \rVert_{H^1}^2) \, d\tau \\
&+C\sum\limits_{i=1}^2\delta_i^2 \int_{0}^{t} \left( \lVert \sqrt{|p'(v)|}\nabla_x \phi \rVert_{L^2}^2+ \lVert \nabla_x \psi \rVert_{L_2}^2+\bold{D}(\tau)+\mathcal{G}^S(\tau) \right) \, d\tau+C\sum\limits_{i=1}^2\delta_i^2 \varepsilon_1.    
\end{aligned}    
\end{equation*}
We use \eqref{presentvF} again. Then as in $\mathcal{B}_4$, 
\begin{equation*}
\begin{aligned}
K_{5,4}(t)&\le C\sum\limits_{i=1}^2 \delta_i^2 \int_{0}^{t} \int_{\mathbb{T}^2} \int_{\mathbb{R}} \left| \nabla_x \psi \right| \left(  \left| \nabla_x^2 \phi \right|+\left| \nabla_x^2 \psi \right|+\left| \nabla_x \phi \right|+\left| \nabla_x \psi \right|\right) \, dx_1 \, dx' \, d\tau+C\delta_i^2\int_{0}^{t}\delta_1\delta_2 e^{-C\min\left\{ \delta_1, \delta_2 \right\}t}d\tau \\ 
&\le C\sum\limits_{i=1}^2 \delta_i^2 \int_{0}^{t} \left( \lVert \sqrt{|p'(v)|} \nabla_x^2\phi \rVert_{L^2}^2+\lVert \nabla_x^2\psi \rVert_{H^1}^2+\bold{D}(\tau)+\mathcal{G}^S(\tau)+\varepsilon_1\nu_i\delta_ie^{-C\delta_it}  \right)d\tau+C\delta_0.     
\end{aligned}
\end{equation*}
The same as $K_{5,3}(t)$, we have
\begin{equation*}
\begin{aligned}
K_{5,5}(t)\le& C\int_{0}^{t} \int_{\mathbb{T}^2}  \int_{\mathbb{R}} \left| \nabla_x \psi \right| \left( \left| \nabla_x^2 \phi \right| \left| \nabla_x \psi \right|+\left| \nabla_x \phi \right|\left| \nabla_x^2 \psi \right| \right) \, dx_1 \, dx' \, d\tau \\
&+C\int_{0}^{t}  \int_{\mathbb{T}^2} \int_{\mathbb{R}} \left| \nabla_x \psi \right|  \left(  \left| (\widetilde{v})_{x_1x_1} \right| \left| \nabla_x \psi \right| +   \left| (\widetilde{v})_{x_1} \right| \left| \nabla_x^2 \psi \right| \right) \,dx_1 \, dx' \, d\tau\\
\le&C\varepsilon_1\int_{0}^{t} \left( \lVert \sqrt{|p'(v)|}\nabla_x^2 \phi \rVert_{L_2}^2 + \lVert \nabla_x \psi \rVert_{H^1}^2 \right) \, d\tau+C\sum\limits_{i=1}^2 \delta_i^2 \int_{0}^{t}  \lVert \nabla_x \psi \rVert_{H^1}^2  \, d\tau.    
\end{aligned}    
\end{equation*}
By integration by parts over $\Omega$,
\begin{equation*}
K_{5,6}(t)=\int_{0}^{t} \lVert \nabla_x div_x \psi \rVert_{L^2}^2 d\tau.    
\end{equation*}
Combining of the estimates above, we get
\begin{equation*}
\begin{aligned}
K_5(t)\le \frac{1}{8} &\int_{0}^{t} \lVert \sqrt{|p'(v)|} \nabla_x^2\phi \rVert_{L_2}^2 \, d\tau + C\int_{0}^{t} \lVert \nabla_x^2\psi \rVert^2 \, d\tau+C\sum\limits_{i=1}^2 \delta_i^2 \int_{0}^{t} | \dot{X}_i(\tau) |^2 \, d\tau\\
&+C\sum\limits_{i=1}^2 \delta_i^2 \int_{0}^{t} 	\left(\bold{G}_1(\tau)+ \bold{G}_3(\tau)+\mathcal{G}^S(\tau)+\bold{D}(\tau) \right) \, d\tau + C\sum\limits_{i=1}^2(\delta_i+\varepsilon_1) \int_{0}^{t} \lVert \nabla_x \psi \rVert_{L_2}^2 d\tau +C\delta_0.    
\end{aligned}
\end{equation*}
By Holder's inequality and \eqref{GNIAPP4},
\begin{equation*}
\begin{aligned}
K_6(t)\le &C\int_{0}^{t} \int_{\mathbb{T}^2} \int_{\mathbb{R}} \left| \nabla_x^2\phi \right| \left( \left| \nabla_x^2\psi \right| \left| \nabla_x\phi \right| + \left| \nabla_x\psi \right|\left| \nabla_x^2 \phi \right| + \left| \nabla_x \psi \right|^2 + \left| \nabla_x \phi \right|^2 \right) \, dx_1 \, dx' \, d\tau \\
&+C\int_{0}^{t} \int_{\mathbb{T}^2} \int_{\mathbb{R}} \left| \nabla_x^2\phi \right| \left( \left| \partial_1 \phi \right| \left| (\widetilde{u})_{x_1x_1} \right| + \left| \partial_1\psi \right| \left| (\widetilde{u})_{x_1} \right| + \left| \nabla_x^2\phi \right|\left| (\widetilde{u})_{x_1} \right| + \left| \nabla_x \phi \right| \left| (\widetilde{v})_{x_1} \right| \right) \, dx_1 \, dx' \, d\tau\\
\le &C\int_{0}^{t} \lVert \nabla_x^2\phi \rVert_{L^2}\left[ \lVert \nabla_x^2\psi \rVert_{L^6}\lVert \nabla_x\phi \rVert_{L^3}+\lVert \nabla_x\psi \rVert_{L^6}\lVert \nabla_x\psi \rVert_{L^3}+\lVert \nabla_x^2\phi \rVert_{L^2}\lVert \nabla_x\psi \rVert_{L^\infty} \right] \, d\tau\\
&+C\int_{0}^{t} \lVert \nabla_x^2\phi \rVert_{L^2} \lVert \nabla_x\phi \rVert_{L^3}\lVert \nabla_x\phi \rVert_{L^6} \, d\tau+C(\delta_1^2+\delta_2^2)\int_{0}^{t} \lVert \nabla_x^2\phi \rVert_{L^2} \left( \lVert \nabla_x(\phi,\psi) \rVert_{L^2}+\lVert \nabla_x^2\phi \rVert_{L^2} \right) \, d\tau\\
\le &C(\delta_0+\varepsilon_1)\int_{0}^{t} \lVert \sqrt{\left| p'(v) \right|}\nabla_x^2\phi \rVert_{L^2}^2+\bold{D}(\tau)+\mathcal{G}^S(\tau)+\lVert \sqrt{\left| p'(v) \right|}\nabla_x\psi \rVert_{H^2}^2\, d\tau+C\delta_0,    
\end{aligned}    
\end{equation*}
where we use the fact that (from Lemma \ref{infty interpolation inequality lemma})
\begin{equation*}
\lVert \nabla_x^2\phi \rVert_{L^2}^2\lVert \nabla_x\psi \rVert_{L^\infty}\le C\lVert \nabla_x^2\phi \rVert_{L^2}^2\lVert \nabla_x\psi \rVert_{H^2}\le C\varepsilon_1 \lVert \nabla_x\phi \rVert_{L^2}^2\lVert \nabla_x\psi \rVert_{H^2}\le C\varepsilon_1\left( \lVert \sqrt{	\left| p'(v) \right|}\nabla_x^2\phi \rVert_{L^2}^2 + \lVert \nabla_x\psi \rVert_{H^2}^2 \right).
\end{equation*}
As above, we get
\begin{equation*}
\begin{aligned}
K_7(t)\le &C\sum\limits_{i=1}^2 \int_{0}^{t} \int_{\mathbb{T}^2} \int_{\mathbb{R}}	\left( \left| \nabla_x\phi \right|\left| \phi \right|+\left| \phi \right|\left| (\widetilde{v}_i)_{x_1}^{X_i} \right|+\left| \nabla_x\phi \right| \right)\left|\nabla_x^2\phi \right| \left| (\widetilde{v}_i)_{x_1}^{X_i} \right| \, dx_1 \, dx' \, d\tau\\
&+C\sum\limits_{i=1}^2\delta_i \int_{0}^{t} \int_{\mathbb{T}^2} \int_{\mathbb{R}}	\left| \phi \right| \left| \partial_{x_1}^2\phi\right|\left| (\widetilde{v}_i)_{x_1}^{X_i} \right| \, dx_1 \, dx' \, d\tau \\
\le& C\sum\limits_{i=1}^2\delta_i \int_{0}^{t} \left( \lVert \sqrt{\left| p'(v) \right|} \nabla_x^2\phi \rVert_{L^2}^2+\bold{D}(\tau)+\mathcal{G}^S(\tau) \right) d\tau+C\delta_0.
\end{aligned}    
\end{equation*}
Likewise, using \eqref{GNIAPP4} we can estimate $K_8(t)$ as 
\begin{equation*}
\begin{aligned}
K_8(t)\le &\int_{0}^{t} \int_{\mathbb{T}^2} \int_{\mathbb{R}} \left| \nabla_x^2 \psi \right| \left| \nabla_x^2 \phi \right| \left( \left| \nabla_x \phi \right|+ \left| (\widetilde{v})_{x_1} \right|  \right) \, dx_1 \, dx' \, d\tau\\
&+C\int_{0}^{t} \int_{\mathbb{T}^2} \int_{\mathbb{R}} \left| \nabla_x^2 \phi \right| \left( \left| \nabla_x \psi \right| \left| \nabla_x^2 \phi \right|+  \left| \nabla_x \psi \right|\left| (\widetilde{v})_{x_1x_1} \right|  \right) \, dx_1 \, dx' \, d\tau\\
\le &C\int_{0}^{t} \left[ \lVert \nabla_x\phi \rVert_{L^3} \lVert \nabla_x^2\psi \rVert_{L^6} \lVert \nabla_x^2\phi \rVert_{L^2} +\lVert \nabla_x\psi \rVert_{L^\infty}\lVert \nabla_x^2\phi \rVert_{L^2}^2\right] \, d\tau\\
&+C\delta_0^2 \int_{0}^{t} \lVert \nabla_x^2\phi \rVert_{L^2}\left[  \lVert \nabla_x^2\psi \rVert_{L^2}+\lVert \nabla_x\psi \rVert_{L^2} \right] \, d\tau\\
\le &C(\delta_0+\varepsilon_1)\int_{0}^{t} \left[ \lVert \sqrt{\left| p'(v) \right|} \nabla_x^2\phi\rVert_{L^2}+\lVert \nabla_x\psi \rVert_{H^2}^2 \right] \, d\tau.
\end{aligned}    
\end{equation*}
Finally, to estimate $K_9(t)$, note that
\begin{equation*}
\begin{aligned}
\partial_{x_1}^2(p(\widetilde{v})-p(\widetilde{v}_1^{X_1})-p(\widetilde{v}_2^{X_2}) )=&( p''(\widetilde{v})-p''(\widetilde{v}_1^{X_1}) )| (\widetilde{v}_1')^{X_1}|^2+( p''(\widetilde{v})-p''(\widetilde{v}_2^{X_2}) )| (\widetilde{v}_2')^{X_2} |^2\\
&+( p'(\widetilde{v})-p'(\widetilde{v}_1^{X_1}))(\widetilde{v}_1'')^{X_1}+( p'(\widetilde{v})-p'(\widetilde{v}_2^{X_2}))(\widetilde{v}_2'')^{X_2}+2p''(\widetilde{v})(\widetilde{v}_1')^{X_1}(\widetilde{v}_2')^{X_2}.
\end{aligned}
\end{equation*}
Using this fact, Lemma \ref{lem:shock-est}, and \ref{shock interaction lemma1}, we get
\begin{equation*}
\begin{aligned}
K_9(t)&\le C \int_{0}^{t} \int_{\mathbb{T}^2} \int_{\mathbb{R}} \left| \nabla_x^2\phi \right| |\partial_{x_1}^2( p(\widetilde{v})-p(\widetilde{v}_1^{X_1})-p(\widetilde{v}_2^{X_2}) )| \, dx_1 \, dx' \, d\tau\\
&\le C\int_{0}^{t} \left| \nabla_x^2\phi \right| \left( \delta_1 |\widetilde{v}-\widetilde{v}_1^{X_1}| | (\widetilde{v}_1)_{x_1}^{X_1}|+\delta_2|\widetilde{v}-\widetilde{v}_2^{X_2}| |(\widetilde{v}_2)_{x_1}^{X_2}| + |(\widetilde{v}_1)_{x_1}^{X_1}| |(\widetilde{v}_2)_{x_1}^{X_2}| \right)\, d\tau\\
&\le C(\delta_1^2+\delta_2^2)\int_{0}^{t} \lVert \sqrt{\left| p'(v)\right|} \nabla_x^2\phi \rVert_{L^2}^2  \, d\tau + C\delta_0.
\end{aligned}
\end{equation*}
Combining the estimates above and using Lemma \ref{lemma4.1}, Lemma \ref{lemma5.1}, and Lemma \ref{lem:2nd}, we can obtain the desired result \eqref{lemma5.3eq}, which completes the proof of Lemma \ref{lem:3rd}.

\section{Proof of Lemma \ref{lem:4th}}
Multiply $(5.10)_2$ by $-\Delta_x\partial_{x_j}\psi$ and then sum $j$ from $1$ to $3$. 

\noindent In addition, integrating the result over $\Omega$, we get
\begin{equation}
\begin{aligned}
\frac{d}{dt}\int_{\mathbb{T}^2}\int_{\mathbb{R}} \frac{\left| \nabla_x^2 \psi \right|^2}{2} \, &dx_1  \,dx'+\underbrace{\mu \int_{\mathbb{T}^2}\int_{\mathbb{R}} v \left| \nabla_x \Delta_x \psi \right|^2 \, dx_1  \,dx'+(\mu+\lambda) \int_{\mathbb{T}^2}\int_{\mathbb{R}} v \left| \nabla_x^2 div_x\psi \right|^2 \, dx_1  \,dx'}_{\mathcal{D}_3(t)}\\
&=:\sum\limits_{i=1}^8L_i(t),   
\end{aligned}
\end{equation}
where
\begin{equation}
\begin{aligned}
L_1(t)=&-\sum\limits_{j,k=1}^3\int_{\mathbb{T}^2} \int_{\mathbb{R}} 
\partial_{x_k}\bold{u}\cdot \nabla_x\partial_{x_j} \psi \partial_{x_k} \partial_{x_j} \psi \, dx_1 \, dx'+\int_{\mathbb{T}^2} \int_{\mathbb{R}} 
div_x\bold{u}\frac{\left| \nabla_x^2\psi \right|^2}{2} \, dx_1 \, dx',\\
L_2(t)=&\sum\limits_{j=1}^3\int_{\mathbb{T}^2}\int_{\mathbb{R}} vp'(v)\Delta_x \partial_{x_j}\psi \cdot \nabla_x\partial_{x_j} \phi \, dx_1 \, dx',\\
L_3(t)=&\sum\limits_{i=1}^3 \int_{\mathbb{T}^2}\int_{\mathbb{R}} \partial_{x_j}\bold{u} \cdot \nabla_x\psi \Delta_x\partial_{x_j}\psi \, dx_1 \, dx'+\sum\limits_{i=1}^3\int_{\mathbb{T}^2}\int_{\mathbb{R}} \partial_{x_j}\left( vp'(v) \right) \Delta_x\partial_{x_j}\psi \cdot \nabla_x\phi \, dx_1 \, dx',\\
L_4(t)=&\sum\limits_{j=1}^3\int_{\mathbb{T}^2} \int_{\mathbb{R}} \partial_{x_j}(v(p'(v)-p'(\widetilde{v}))) (\widetilde{v})_{x_1} \Delta_x\partial_{x_j}\psi_1 \, dx_1 \, dx'\\
&+\int_{\mathbb{T}^2} \int_{\mathbb{R}} v(p'(v)-p'(\widetilde{v})) (\widetilde{v})_{x_1x_1} \Delta_x\partial_1\psi_1 \, dx_1 \, dx',\\
L_5(t)=&-\sum\limits_{i=1}^2\dot{X}_i(t)\int_{\mathbb{T}^2} \int_{\mathbb{R}} \Delta_x\partial_{x_1}\psi_1 (\widetilde{u}_{i})_{x_1x_1}^{X_i} \, dx_1 \, dx',\\
L_6(t)=&\sum\limits_{i=1}^2\int_{\mathbb{T}^2} \int_{\mathbb{R}} vF_i(\widetilde{u}_{i})_{x_1x_1}^{X_i} \Delta_x\partial_{x_1} \psi_1\, dx_1 \, dx',\\
L_7(t)=&\sum\limits_{i=1}^2 \sum\limits_{j=1}^3 \int_{\mathbb{T}^2} \int_{\mathbb{R}} \partial_{x_j}(vF_i)(\widetilde{u}_{i})_{x_1}^{X_i} \Delta_x \partial_{x_j}\psi_1 \, dx_1 \, dx',\\
L_8(t)=&(\mu+\lambda)\sum\limits_{j=1}^3 \int_{\mathbb{T}^2} \int_{\mathbb{R}} \left( \Delta_x\partial_{x_j}\psi-\nabla_x\partial_{x_j} div_x\psi \right)\cdot \nabla_x v \partial_{x_j} div_x\psi \, dx_1 \, dx'\\
&- \sum\limits_{j=1}^3 \int_{\mathbb{T}^2} \int_{\mathbb{R}} \partial_{x_j}v \left( \mu\Delta_x\psi + (\mu+\lambda)\nabla_x div_x\psi \right) \cdot \Delta_x \partial_{x_j} \psi \, dx_1 \, dx',\\
L_9(t)=&\sum\limits_{j=1}^3 \int_{\mathbb{T}^2} \int_{\mathbb{R}} \Delta_x \partial_{x_j} \psi \cdot \partial_{x_j} \nabla_x ( p(\widetilde{v})-p(\widetilde{v}_1^{X_1})-p(\widetilde{v}_2^{X_2}) ) \, dx_1 \, dx'.        
\end{aligned}
\end{equation}
By Holder's inequality and \eqref{GNIAPP4}, 
\begin{equation*}
\begin{aligned}
L_1(t)&\le C\lVert \nabla_x\psi \rVert_{L^3}\lVert \nabla_x^2\psi \rVert_{L^6}\lVert \nabla_x^2\psi \rVert_{L^2}+(\delta_1^2+\delta_2^2)\lVert \nabla_x^2\psi \rVert_{L^2}^2\\
&\le C(\delta_0+\varepsilon_1)\lVert \nabla_x^2\psi \rVert_{H^1}^2\\
&\le C(\delta_0+\varepsilon_1)(\mathcal{D}_3(t)+\lVert \nabla_x^2\psi \rVert_{L^2}^2).
\end{aligned}
\end{equation*}
Using Cauchy-Schwartz inequality, we get
\begin{equation*}
\begin{aligned}
L_2(t)&\le \int_{\mathbb{T}^2} \int_{\mathbb{R}} \left| \nabla_x \Delta_x \psi \right| \left| \nabla_x^2\phi \right| \, dx_1 \, dx'\\
&\le \frac{1}{8}\mathcal{D}_3(t)+C\lVert \nabla_x^2 \phi \rVert_{L^2}^2.   
\end{aligned}    
\end{equation*}
Using Lemma 2.3, we have
\begin{equation*}
\begin{aligned}
L_3(t)\le &\int_{\mathbb{T}^2} \int_{\mathbb{R}} \left| \nabla_x\psi \right|^2 \left| \nabla_x \Delta_x \psi \right| \, dx_1 \, dx'+C(\delta_1^2+\delta_2^2)\int_{\mathbb{T}^2} \int_{\mathbb{R}} \left| \nabla_x\psi \right| \left| \nabla_x \Delta_x \psi \right| \, dx_1 \, dx'\\
&+\int_{\mathbb{T}^2} \int_{\mathbb{R}} \left| \nabla_x \phi \right|^2 \left| \nabla_x \Delta_x \psi \right| \, dx_1 \, dx'+C(\delta_1^2+\delta_2^2)\int_{\mathbb{T}^2} \int_{\mathbb{R}} \left| \nabla_x \phi \right| \left| \nabla_x \Delta_x \psi \right| \, dx_1 \, dx'\\
\le &C\left[ \lVert \nabla_x\psi \rVert_{L^3}\lVert \nabla_x\psi \rVert_{L^6}+\lVert \nabla_x\phi \rVert_{L^3}\lVert \nabla_x\phi \rVert_{L^6}+(\delta_1^2+\delta_2^2)\lVert \nabla_x(\phi,\psi) \rVert_{L^2} \right]\sqrt{\mathcal{D}_3(t)}\\
\le &C(\delta_0+\varepsilon_1)(\mathcal{D}_3(t)+\lVert \nabla_x\psi \rVert_{H^1}^2+\lVert \nabla_x\phi \rVert_{H^1}^2)\\
\le &C(\delta_0+\varepsilon_1)(\mathcal{D}_3(t)+\mathcal{G}^S(t)+\bold{D}(t)+\sum\limits_{i=1}^2 \delta_i^2e^{-C\delta_it} \int_{\mathbb{T}^2}  \int_{\mathbb{R}} \eta(U|\widetilde{U}) \, dx_1 \, dx'+\lVert \nabla_x\psi \rVert_{H^1}^2+\lVert \nabla_x^2\phi \rVert_{L^2}^2).  
\end{aligned}    
\end{equation*}
As above,
\begin{equation*}
\begin{aligned}
L_4(t)&\le C\int_{\mathbb{T}^2} \int_{\mathbb{R}}  (\left| \nabla_x \phi \right|+\left| (\widetilde{v})_{x_1} \right|\left| \phi \right|)\left|( \widetilde{v})_{x_1} \right| \left| \nabla_x^3 \psi_1 \right| \, dx_1 \, dx'+C\sum\limits_{i=1}^2 \delta_i \int_{\mathbb{T}^2} \int_{\mathbb{R}} \left|\phi \right| |( \widetilde{v}_i)_{x_1}^{X_i}|\left| \nabla_x^3 \psi_1 \right| \, dx_1 \, dx'\\
&\le C\delta_0(\mathcal{D}_3(t)+\mathcal{G}^S(t)+\bold{D}(t)+\sum\limits_{i=1}^2 \delta_i^2e^{-C\delta_it} \int_{\mathbb{T}^2}  \int_{\mathbb{R}} \eta(U|\widetilde{U}) \, dx_1 \, dx'), \\
L_5(t)&\le C\sum\limits_{i=1}^2\delta_i| \dot{X}_i(t) |\sqrt{\mathcal{D}_3(t)}\lVert (\widetilde{v}_{i})_{x_1}^{X_i}\rVert_{L^2}\\
&\le C\sum\limits_{i=1}^2\delta_i^{\frac{5}{2}}| \dot{X}_i(t) |\sqrt{\mathcal{D}_3(t)}\\
&\le C\sum\limits_{i=1}^2\delta_i^{\frac{5}{2}} (\lVert \dot{X}_i(t)\rVert_{L^2}^2+\mathcal{D}_3(t)).
\end{aligned}
\end{equation*}
Using the definition of $F_i$, as in $\mathcal{B}_3$, we get
\begin{align*}
L_6(t)\le& \sum\limits_{i=1}^2 \delta_i^2\Bigg( \sqrt{\bold{G}_1}+\sqrt{\bold{G}_3}+\sqrt{\mathcal{G}^S}+\sqrt{\bold{D}}\\
&\qquad \qquad +\sqrt{\varepsilon_1 \nu_i \delta_i e^{-C\delta_it}}+\sqrt{\delta_1\delta_2 e^{-C\min\left\{ \delta_1, \delta_2  \right\}t}}+\sqrt{\delta_0^2 \nu_i \delta_i e^{-C\delta_it}} \Bigg) \lVert \Delta_x\partial_{x_1} \psi_1 \rVert_{L^2}. \\
\le& \sum\limits_{i=1}^2 \delta_i^2\left( \bold{G}_1+\bold{G}_3+\mathcal{G}^S+\bold{D}+\varepsilon_1\nu_i\delta_ie^{-C\delta_it}+\delta_1\delta_2e^{-C\min\left\{ \delta_1,\delta_2 \right\}t}+\delta_0^2\nu_i\delta_ie^{-C\delta_it}+\mathcal{D}_3 \right).
\end{align*}
As above, with \eqref{presentvF} and lemma \ref{shock interaction lemma1}, we get
\begin{align*}
L_7(t)&\le \sum\limits_{i=1}^2 \delta_i^2\left( \lVert\nabla_x(\phi,\psi_1)\rVert_{L^2}+\lVert  (\widetilde{v}_i)_{x_1}^{X_i} (\widetilde{v}-\widetilde{v}_i^{X_i}) \rVert_{L^2} \right)\lVert\Delta_x\partial_{x_i}(\phi,\psi_1)\rVert_{L^2} \\
&\le \sum\limits_{i=1}^2 \delta_i^2\left( \bold{G}_1+\bold{G}_3+\mathcal{G}^S+\bold{D}+\varepsilon_1\nu_i\delta_ie^{-C\delta_it}+\delta_1\delta_2e^{-C\min\left\{ \delta_1,\delta_2 \right\}t}+\delta_0^2\nu_i\delta_ie^{-C\delta_it}+\mathcal{D}_3 \right).
\end{align*}
Using Holder's inequality and \eqref{GNIAPP4}, we get
\begin{align*}
L_8(t)\le& C \int_{\mathbb{T}^2}\int_{\mathbb{R}}  (\left|\Delta_x\nabla_x \psi\right|+\left|\nabla_x^2div_x\psi\right|)(\left|\phi\right|+\delta_1^2+\delta_2^2)\left|\nabla_x div_x\psi\right|) \, dx_1 \, dx'\\ 
&+ C \int_{\mathbb{T}^2}\int_{\mathbb{R}}  (\left| \nabla_x\phi \right|+\delta_1^2+\delta_2^2)(\left| \Delta_x\psi \right|+\left| \nabla_x div_x\psi \right|)\left| \Delta_x\nabla_x \psi \right|  \, dx_1 \, dx' \\
\le& \sqrt{\mathcal{D}_3(t)}\lVert \nabla_x\phi \rVert_{L^3}\left(\lVert \Delta_x\psi \rVert_{L^6}+\lVert \nabla_x div_x\psi \rVert_{L^6}\right)+(\delta_1^2+\delta_2^2)\sqrt{\mathcal{D}_3(t)}(\lVert \Delta_x\psi \rVert_{L^2}+\lVert \nabla_x div_x\psi \rVert_{L^2}) \\
\le& (\delta_0+\varepsilon_1)\sqrt{\mathcal{D}_3(t)}(\lVert \Delta_x\psi \rVert_{H^1}+\lVert \nabla_x div_x\psi \rVert_{H^1}) \\
\le& (\delta_0+\varepsilon_1)(\mathcal{D}_3(t)+\lVert \Delta_x\psi \rVert_{H^1}^2+\lVert \nabla_x div_x\psi \rVert_{H^1}^2).
\end{align*}
Finally, as in $K_9$, we get
\begin{align*}
L_9(t)&\le \int_{\mathbb{T}^2}\int_{\mathbb{R}} \left| \nabla_x\Delta_x \psi \right| | \partial_{x_1}^2 (p(\widetilde{v})-p(\widetilde{v}_1^{X_1})-p(\widetilde{v}_2^{X_2})) | \, dx_1 \, dx' \\
&\le C(\delta_1^2+\delta_2^2)\mathcal{D}_3(t)+C\delta_0.
\end{align*}
Combine the estimates above and then integrate result over $\left[ 0,t \right]$ for $t\in \left[0, T\right]$. \\
\noindent Also, using Lemmas 4.1, 5.1, 5.2, and 5.3, we can get the desired result \eqref{5.17}, which completes the proof of Lemma 5.4.

\end{appendix}

\bibliographystyle{amsplain}
\bibliography{reference} 
\end{document}